\documentclass{amsart}
\usepackage{latexsym,amsxtra,amscd,ifthen}
\usepackage{newunicodechar}

\usepackage{cleveref}
\usepackage{amsfonts}
\usepackage{color,graphicx}
\usepackage{verbatim}
\usepackage{amsmath}
\usepackage{amsthm}
\usepackage{amssymb}
\usepackage{url}
\usepackage{todonotes}

\usepackage{tikz}
\usetikzlibrary{arrows,shapes,chains}

\theoremstyle{plain}

\newtheorem{theorem}{Theorem}
\newtheorem{lemma}[theorem]{Lemma}
\newtheorem{proposition}[theorem]{Proposition}
\newtheorem{construction}[theorem]{Construction}
\newtheorem{corollary}[theorem]{Corollary}

\newtheorem{convention}[theorem]{Convention}
\newtheorem{notation}[theorem]{Notation}

\numberwithin{theorem}{section}
\numberwithin{equation}{theorem}

\theoremstyle{definition}
\newtheorem{definition}[theorem]{Definition}
\newtheorem{example}[theorem]{Example}
\newtheorem{remark}[theorem]{Remark}
\newtheorem{hypothesis}[theorem]{Hypothesis}

\newtheorem*{question*}{Question}

\newcommand{\fm}{\mathfrak{m}}

\newcommand{\fp}{\mathfrak{p}}

\newcommand{\bfq}{\mathbf{q}}

\newcommand{\bfx}{\mathbf{x}}
\newcommand{\bfw}{\mathbf{w}}
\newcommand{\bfd}{\mathbf{d}}

\newcommand{\Supp}{\mathbb{S}}
\newcommand{\OO}{\mathbb{O}}

\DeclareMathOperator{\Proj}{Proj}

\DeclareMathOperator{\Aut}{Aut}
\DeclareMathOperator{\gr}{gr}

\DeclareMathOperator{\GKdim}{GKdim}
\DeclareMathOperator{\gldim}{gldim}
\DeclareMathOperator{\MaxSpec}{MaxSpec}

\DeclareMathOperator{\Reg}{Reg}
\DeclareMathOperator{\Spec}{Spec}

\newcommand\kk{{\Bbbk}}

\begin{document}

\title{Valuation method for Nambu-Poisson algebras}

\author{Hongdi Huang, Xin Tang, Xingting Wang and James J. Zhang}

\address{Huang: Department of Mathematics, Rice University,
Houston, TX 77005, USA}

\email{hh40@rice.edu}

\address{Tang: Department of Mathematics \& Computer Science,
Fayetteville State University, Fayetteville, NC 28301,
USA}

\email{xtang@uncfsu.edu}

\address{Wang: Department of Mathematics, Louisiana State University, 
Baton Rouge, Louisiana 70803, USA}
\email{xingtingwang@math.lsu.edu}

\address{Zhang: Department of Mathematics, Box 354350,
University of Washington, Seattle, Washington 98195, USA}

\email{zhang@math.washington.edu}

\begin{abstract}
Automorphism, isomorphism, and embedding problems are
investigated for a family of Nambu-Poisson algebras
(or $n$-Lie Poisson algebras) using Poisson valuations. 
\end{abstract}


%17B63 (2000-now) Poisson algebras
%17B40 (1973-now) Automorphisms, derivations, other operators for Lie algebras and superalgebras
%16W20 (1991-now) Automorphisms and endomorphisms


%\thanks{ }

%\date{\today}
\maketitle

%\tableofcontents

% \setcounter{section}{-1}
\section*{Introduction}
\label{xxsec0}

A Nambu-Poisson manifold was introduced in the work of Nambu 
in 1973 \cite{Nam}. There are several nice survey articles 
about geometric aspects of the Nambu-Poisson bracket; see 
\cite{Nak, Va}. The present paper concerns the corresponding
notion of a Nambu-Poisson algebra. An equivalent notion 
called $n$-Lie Poisson algebra was introduced by 
Umirbaev-Zhelyabin in \cite{UZ} in 2021. A Nambu-Poisson 
algebra is exactly an $n$-Lie Poisson algebra for any 
$n\geq 3$, see Section \ref{xxsec1} for definitions. Some 
examples of Nambu-Poisson algebras are given in Example
\ref{xxex1.3}. 

This paper deals with automorphism, isomorphism, and embedding 
problems for several classes of Nambu-Poisson algebras and 
Nambu-Poisson fields. We will use the Poisson 
valuation concept introduced recently by the authors \cite{HTWZ3}
to tackle these important questions. In \cite{HTWZ1, HTWZ2, HTWZ3}, we have used valuations to solve rigidity of grading, 
rigidity of filtration, automorphism, isomorphism, and 
embedding problems for Poisson algebras and fields. This paper aims to show that  
the valuation method works well for Nambu-Poisson algebras. 

Let $\Bbbk$ be a base field. For simplicity, we assume that 
$\Bbbk$ is of characteristic 0, though all definitions make 
sense without this extra assumption. Most algebraic objects 
are over $\Bbbk$. Let $n$ be a fixed integer $\geq 2$ and, in 
a large part of the paper, $n\geq 3$. 

\begin{definition}
\label{xxdef0.1}
Let $N$ be a Nambu-Poisson field (or an $n$-Lie Poisson field)
with $n$-Lie Poisson bracket $\{-,\cdots,-\}: N^{\otimes n}\to
N$. A map 
$$\nu: N\to {\mathbb Z}\cup\{\infty\}$$
is called a {\it $0$-valuation} on $N$ if, for all 
elements $a,b_1,\cdots,b_n\in N$,
\begin{enumerate}
\item[(1)]
$\nu(a)=\infty$ if and only if $a=0$,
\item[(2)]
$\nu(a)=0$ for all $a\in \Bbbk^{\times}:=\Bbbk\setminus \{0\}$,
\item[(3)]
$\nu(b_1 b_2)=\nu(b_1)+\nu(b_2)$,
\item[(4)]
$\nu(b_1+b_2)\geq \min\{\nu(b_1),\nu(b_2)\}$, with equality if 
$\nu(b_1)\neq \nu(b_2)$,
\item[(5)]
$\nu(\{b_1,b_2,\cdots,b_n\})\geq \sum_{s=1}^n \nu(b_{s})$.
\end{enumerate}
\end{definition}

In Definition \ref{xxdef1.6}, we will introduce a slightly more 
general version of a $0$-valuation, called $(\OO,w)$-valuation, 
for an ordered group $\OO$ and an element $w\in \OO$. Given a  
$0$-valuation $\nu$, we can define a filtration
${\mathbb F}^{\nu}:=\{F_i^{\nu}\mid i\in {\mathbb Z}\}$ of $N$ by
$$F_i^{\nu}:=\{a\in N\mid \nu(a)\geq i\}, \quad {\text{for all $i$}}.$$
The {\it associated graded ring} of $\nu$ is defined to be 
$$\gr_{\nu}(N)=\bigoplus_{i\in {\mathbb Z}} F_i^{\nu}/F_{i+1}^{\nu}.$$
Note that there is an induced Nambu-Poisson bracket on
$\gr_{\nu}(N)$. The degree 0 part $(\gr_{\nu}(N))_0$ is called the 
{\it residue field} of $\nu$.

In this paper, we will use the Gelfand-Kirillov dimension (or 
GK-dimension for short) instead of transcendence degree. We refer 
to the book \cite{KL} for basic definitions and properties 
related to the GK-dimension. If $B$ is a field over $\Bbbk$, then 
the GK-dimension of $B$ is equal to the transcendence degree of 
$B$ (over $\Bbbk$). If $B$ is a commutative affine algebra, then 
the GK-dimension of $B$ equals the Krull dimension of $B$. 
We introduce two more definitions before we state some results.

\begin{definition}
\label{xxdef0.2}
A $0$-valuation $\nu$ on a Nambu-Poisson field $N$ is 
{\it faithful} if
\begin{enumerate}
\item[(1)]
the image of $\nu$ is ${\mathbb Z}\cup\{\infty\}$,
\item[(2)]
the induced Nambu-Poisson bracket on $\gr_{\nu}(N)$ is nonzero, 
and
\item[(3)]
$\GKdim \gr_{\nu}(N)=\GKdim N$.
\end{enumerate}
\end{definition}

Let $Q$ be another Nambu-Poisson field. We write $N\to_{\nu} Q$ 
if $Q$ is isomorphic to the Nambu-Poisson fraction field of 
the Nambu-Poisson associated graded domain $\gr_{\nu}(N)$.

\begin{definition}
\label{xxdef0.3}
Let $N$ be a Nambu-Poisson field. 
\begin{enumerate}
\item[(1)] The {\it depth} of $N$ is defined to be
$$\bfd(N):=\sup\{n\mid N=N_0\to_{\nu_1} N_1\to_{\nu_2}
N_2 \to \cdots \to_{\nu_{n}} N_n\}$$
where each $\nu_i$ is a faithful $0$-valuation and 
$N_i\not\cong N_j$ for all $0\leq i\neq j\leq n$.
\item[(2)]
The {\it width} of $N$ is defined to be
$$\bfw(N):=\# \{ [Q \mid N\to_{\nu} Q]/\cong\}$$
where $\nu$ runs over all faithful $0$-valuations of $N$ and 
$[-]/\cong$ means the isomorphism classes. 
\end{enumerate}
\end{definition}

Here is one of the main results. Let $m\geq n\geq 2$ be two 
integers. Let $T$ denote a Nambu-Poisson torus 
$\Bbbk[x_1^{\pm 1},\cdots,x_m^{\pm 1}]$ with Nambu-Poisson 
bracket determined by $\{x_{i_1},\cdots,x_{i_n}\}=
q_{i_1\cdots i_n} x_{i_1}\cdots x_{i_n}$ for some parameters 
$q_{i_1\cdots i_n}\in \Bbbk$, for all $1\leq i_1<
\cdots< i_n\leq m$. If $A$ is a Nambu-Poisson domain, we 
use $Q(A)$ to denote the Nambu-Poisson fraction field of
$A$.  

\begin{theorem}[Theorem \ref{xxthm4.3}]
\label{xxthm0.4}
Retain the above notation. Suppose $T$ is Nambu-Poisson simple.
Let $N$ be the Nambu-Poisson fraction field $Q(T)$. Then 
$\bfd(N)=0$ and $\bfw(N)=1$. In particular, for every faithful 
$0$-valuation $\nu$, $Q(\gr_{\nu}(N))\cong N$.
\end{theorem}

The above theorem is a generalization of \cite[Theorem 0.4]{HTWZ3}.
A complete set of $0$-valuations on $N$ as in the above theorem
is given in the proof of Theorem \ref{xxthm4.3}(1). 

The valuation method is useful for many algebraic problems in 
the area of Poisson algebras and Nambu-Poisson algebras. Here 
is a partial list of topics.

\begin{enumerate}
\item[(1)]
Automorphism problem, see Theorems \ref{xxthm0.7} and 
\ref{xxthm0.11}. Computing the automorphism group is always a challenging task. In \cite{MTU2}, Makar-Limanov, Turusbekova, 
and Umirbaev computed the group of the Poisson algebra 
automorphisms for elliptic Poisson algebras. Using Poisson 
valuations, more Poisson automorphism groups  
were computed in \cite{HTWZ1, HTWZ2, HTWZ3}. One can also compute the automorphism group of 
several families of Nambu-Poisson fields using Poisson valuation. Theorems 
\ref{xxthm0.7} and \ref{xxthm0.11} are the first of this kind.
\item[(2)]
Isomorphism problem, see Theorem \ref{xxthm4.14}.
\item[(3)]
Embedding problem, see Theorem \ref{xxthm4.16} and 
Corollaries \ref{xxcor5.9},
\ref{xxcor5.10}, and \ref{xxcor5.11}.
\item[(4)]
Rigidity of grading, see Theorem \ref{xxthm0.9}.
\item[(5)]
Rigidity of filtration, see Theorem \ref{xxthm0.10}.
\item[(6)]
Dixmier problem, see Theorem \ref{xxthm0.8}.
\item[(7)]
Classification project, see \cite{HTWZ2}. Classification of 
Poisson fields of transcendence degree two is an important 
project that is connected with important projects in 
noncommutative algebra and noncommutative algebraic geometry.
\end{enumerate}

To illustrate some ideas of the valuation method, we will 
study several classes of Nambu-Poisson algebras that are 
related to the following construction. 

Let $m$ be a positive integer and let $\Bbbk^{[m]}$ denote the 
(commutative) polynomial ring $\Bbbk[t_1,\cdots, t_{m}]$. 
If $f\in \Bbbk^{[m]}$, the partial derivative 
$\frac{\partial f}{\partial t_i}$ is denoted by $f_{t_i}$.

\begin{construction}\label{xxcon0.5}
We fix an integer $n\geq 3$. Let $\Omega\in \Bbbk^{[n+1]}$, 
called a {\it potential}. 
\begin{enumerate}
\item[(1)]
First we can define a Nambu-Poisson {\rm{(}}$n$-Lie{\rm{)}} 
bracket on the polynomial ring $\Bbbk^{[n+1]}$ by
\begin{equation}
\label{E0.5.1}\tag{E0.5.1}
\{f_1,\cdots, f_n\}
%=Jac(f_1,\cdots, f_n, \Omega)
:=
\det \begin{pmatrix} 
(f_1)_{t_1} & (f_1)_{t_2} & \cdots & (f_1)_{t_n} & (f_1)_{t_{n+1}}\\
(f_2)_{t_1} & (f_2)_{t_2} & \cdots & (f_2)_{t_n} & (f_2)_{t_{n+1}}\\
\cdots      & \cdots      & \cdots & \cdots      & \cdots\\
(f_n)_{t_1} & (f_n)_{t_2} & \cdots & (f_n)_{t_n} & (f_n)_{t_{n+1}}\\
\Omega_{t_1}&\Omega_{t_2} & \cdots & \Omega_{t_n}&\Omega_{t_{n+1}}
\end{pmatrix}
\end{equation}
for all $f_i\in \Bbbk^{[n+1]}$. By \cite[p.578]{UZ}, this is a
Nambu-Poisson algebra and it is denoted by $A_{\Omega}$. {\rm{(}}In 
\cite[p.578]{UZ} it is denoted by $P_{\Omega}$, but we will use 
the symbol $P_{\Omega}$ for a different algebra in sequel.{\rm{)}}
\item[(2)]
It is easy to check that $\Omega$ is in the Nambu-Poisson 
center of $A_{\Omega}$. 
Hence $A_{\Omega}$ has a Nambu-Poisson 
factor ring $P_{\Omega-\xi}:= A_{\Omega}/(\Omega-\xi)$ where 
$\xi\in \Bbbk$. If $\xi=0$, we also write it as $P_{\Omega}$.
\item[(3)]
If $\kk^{[n+1]}$ is graded and $\Omega$ is homogeneous of degree 
equal to $\sum_{s=1}^{n+1} \deg(t_{s})$, then 
$P_{\Omega-\xi}\cong P_{\Omega-1}$ when $\xi\neq 0$. So, 
in this case, we can assume that $\xi$ is either 0 or $1$.
\end{enumerate}
Unless otherwise stated, $A_{\Omega}$, $P_{\Omega}$, 
$P_{\Omega-1}$, and $P_{\Omega-\xi}$ will be the Nambu-Poisson 
algebras defined as above.
\end{construction}

For simplicity, we consider a grading on $\Bbbk^{[n+1]}$ by 
$\deg t_i=1$ for all $i=1,\cdots, n+1$. This is called the 
original grading or Adams grading of $\Bbbk^{[n+1]}$. 
We call an element $\Omega$ in $\Bbbk^{[n+1]}$ a potential regarding Construction \ref{xxcon0.5} 
and we consider homogeneous potentials in this paper. Moreover, $\Omega$ is said to have {\it an isolated singularity} (at the origin) if
$\Bbbk^{[n+1]}/(\Omega_{t_1},\cdots,\Omega_{t_{n+1}})$ 
is finite-dimensional. In this case, we simply say $\Omega$ is 
an i.s. potential. Note that having an isolated singularity is 
independent of the choices of generators 
$(t_1,\cdots,t_{n+1})$. For each $m\geq 2$, the potential
$\Omega:=\sum_{i=1}^{n+1} t_i^{m}$ has an isolated singularity 
(at the origin). In general, a generic homogeneous element of 
degree $\geq 2$ in $\Bbbk^{[n+1]}$ has isolated singularity
\cite[p.355]{MM}. So the results in this paper represent the most common scenario for i.s. potentials.

We will consider two special cases. Case 1: $\Omega$ is an 
i.s. potential with $\deg \Omega=n+1$. When $n=2$, 
$A_{\Omega}$ is called the {\it elliptic Poisson 
algebra} which has been studied by several authors 
\cite{HTWZ1, TWZ, MTU2, Pi1}. Case 2:  $\Omega$ is an 
i.s. potential with $\deg \Omega\geq n+3$. We would like to know more about the degree of the $n+2$ case. 

The following is a result concerning the algebras
in Case 1. 

\begin{theorem}
\label{xxthm0.6}
Let $\Omega\in \Bbbk^{[n+1]}$ be an i.s. potential of degree 
$n+1$.
\begin{enumerate}
\item[(1)]
If $N$ is $Q(P_{\Omega})$. Then 
$\bfd(N)=0$ and $\bfw(N)=1$.  
\item[(2)]
If $N$ is $Q(P_{\Omega-1})$, then $\bfd(N)=\bfw(N)=1$.
\end{enumerate}
\end{theorem}

For Case 2, we have a series of results. Recall that 
$A_{\Omega}$, $P_{\Omega}$, and $P_{\Omega-\xi}$ are defined 
in Construction \ref{xxcon0.5}. If $A$ is a Nambu-Poisson 
algebra, we use $\Aut_{Poi}(A)$ to denote the group of 
Nambu-Poisson algebra automorphisms of $A$.

\begin{theorem}[Proposition \ref{xxpro5.2}]
\label{xxthm0.7}
Let $\Omega\in \Bbbk^{[n+1]}$ be an i.s. potential of degree 
$\geq n+3$.
\begin{enumerate}
\item[(1)]
$$\Aut_{Poi}(Q(P_{\Omega}))=\Aut_{Poi}(P_{\Omega}).$$
If further $\Omega=\sum_{s=1}^{n+1} t_s^{d}$ with $d\geq n+3$,
then $\Aut_{Poi}(P_{\Omega})$ is explicitly given in Lemma 
{\rm{\ref{xxlem6.5}(1)}}.
\item[(2)]
Suppose $\xi\in \Bbbk^{\times}$. Then 
$$\Aut_{Poi}(Q(P_{\Omega-\xi}))=\Aut_{Poi}(P_{\Omega-\xi}).$$
If further $\Omega=\sum_{s=1}^{n+1} t_s^{d}$ with $d\geq n+3$,
then $\Aut_{Poi}(P_{\Omega-\xi})$ is explicitly given in Lemma 
{\rm{\ref{xxlem6.5}(2,3)}}.
\end{enumerate}
\end{theorem}

The automorphism group of Poisson algebras is computed in 
\cite{HTWZ1, HTWZ2, HTWZ3} using Poisson valuations. Motivated 
by the Dixmier conjecture and Poisson conjecture, see 
Section \ref{xxsec5}, we consider the Dixmier property 
[Definition \ref{xxdef5.3}].

\begin{theorem}[Theorem \ref{xxthm5.4}]
\label{xxthm0.8}
Let $\Omega\in \Bbbk^{[n+1]}$ be an i.s. potential of degree 
$\geq n+3$. Then both $P_{\Omega-\xi}$ and 
$Q(P_{\Omega-\xi})$, for any $\xi\in \Bbbk$, satisfy the 
Dixmier property.
\end{theorem}

Other properties, such as rigidity of grading/filtration,
are motivated by some work on noncommutative algebras
\cite[Corollary 0.3]{BZ2}, and the corresponding results 
in the Poisson case \cite[Theorems 0.10 and 0.11]{HTWZ3}. 
The following two results are in this direction. Some undefined 
terms will be defined in Section \ref{xxsec2}.

We say a ${\mathbb Z}$-graded algebra $A=\bigoplus_{i\in 
{\mathbb Z}} A_i$ is {\it connected graded} if $A_i=0$ 
for all $i>0$ and $A_0=\Bbbk$. So, $A$ lives in nonpositive 
degrees. Note that $P_{\Omega}$ is connected graded if we 
set the degree of $t_i$ to $-1$. This is the opposite of 
the usual definition of a connected graded algebra.

\begin{theorem}[Theorem \ref{xxthm5.5}]
\label{xxthm0.9}
Let $\Omega\in \Bbbk^{[n+1]}$ be an i.s. potential of degree 
$n+1+d_0$ where $d_0\geq 2$. Then $P_{\Omega}$ has a unique 
connected grading such that it is  Poisson $d_0$-graded
[Definition \ref{xxdef2.1}]. 
\end{theorem}

\begin{theorem}[Theorem \ref{xxthm5.6}]
\label{xxthm0.10}
Let $\Omega\in \Bbbk^{[n+1]}$ be an i.s. potential of degree 
$n+1+d_0$ where $d_0\geq 2$. Then $P_{\Omega-\xi}$ has a unique 
filtration ${\mathbb F}$ such that $\gr_{\mathbb F} 
(P_{\Omega-\xi})$ is a connected $d_0$-graded Nambu-Poisson domain. 
\end{theorem}

All the above results help compute 
related automorphism groups. For example, we have 

\begin{theorem}[Theorem \ref{xxthm5.7}]
\label{xxthm0.11} 
Let $\Omega\in \Bbbk^{[n+1]}$ be an i.s. potential of degree 
$\geq n+3$. Then, every Nambu-Poisson automorphism of 
$A_{\Omega}$ is graded and $\Aut_{Poi}(Q(A_{\Omega}))
=\Aut_{Poi}(A_{\Omega})=\Aut_{Poi}(P_{\Omega})$. Further,
$\Aut_{Poi}(A_{\Omega})$ is a finite subgroup of
$GL_{n+1}(\Bbbk)$.
\end{theorem}

\begin{remark}
\label{xxrem0.12}
This paper extends the ideas presented in \cite{HTWZ3} in multiple ways. Some new developments in the current paper are the following. 
\begin{enumerate}
\item[(1)]
In \cite{HTWZ3} we only consider ${\mathbb Z}$-valued 
valuations, and in this paper, we consider $\OO$-valued 
valuations for any ordered abelian group $\OO$.
\item[(2)]
In \cite{HTWZ3} we focus on Poisson algebras 
(= $2$-Lie Poisson algebras) and in this paper, we study
$n$-Lie Nambu-Poisson algebras for $n\geq 3$.
\item[(3)]
In \cite{HTWZ3} we introduce $\alpha$-, $\beta$-, and 
$\gamma$-type invariants. In this paper, we also define 
$\delta$-type invariants and use them to solve the 
isomorphism and embedding problems, see Theorems 
\ref{xxthm4.14} and \ref{xxthm4.16}.
\item[(4)]
We explore the possibility of defining $\varsigma$-type 
invariants to study the moduli of faithful $w$-valuations,
see Subsection \ref{xxsec4.5} for one example.
\item[(5)]
A realization lemma is given to connect the Nambu-Poisson
automorphisms and the usually algebraic automorphisms of 
$\Bbbk^{[n+1]}$ with fixed quasi-axis [Section \ref{xxsec6}]. 
This should be useful in the study of the algebraic 
automorphisms of $\Bbbk^{[n+1]}$.
\item[(6)]
We prove the non-existence of solutions for a large class of partial differential equations involving the Jacobian determinant in the field $\Bbbk(t_1,\cdots,t_{n+1})$. Solving PDEs involving the Jacobian determinant can help construct volume-preserving diffeomorphisms with given boundary data. This solution is crucial in establishing prescribed periodic orbits and ergodic mappings. The concept was first explored by Alpern in \cite{Al1976}, and later investigated by Anosov-Katok in \cite{AK1970}. Dacogogna demonstrated in \cite{Da1981} how this method can be applied to the minimization problem in the calculus of variations and nonlinear elasticity. See Section 
\ref{xxsec7}. 
\end{enumerate}
Also, note that the valuation method applies to other algebraic
systems. For example, let ${\mathcal P}$ be an operad with a 
morphism ${\mathcal Com}\to {\mathcal P}$ where ${\mathcal Com}$
is the operad encoding the unital commutative algebras, then 
the valuation method would be useful for studying different 
questions about ${\mathcal P}$-algebras. Some examples are 
given in \cite{GZ} related to multi-Poisson fields. 
\end{remark}

The paper is organized as follows. Section 1 recalls some 
basic definitions of $n$-Lie and Nambu-Poisson algebras and 
valuations. Section 2 gives the definitions of filtration and 
a degree function closely related to the notion of a 
valuation. In Section 3, we work out some examples. In Section 
4, we define more invariants such as the weighted version 
of depth and width. The proofs of the main theorems are given in Section 5. In Section 6, we study the Nambu-Poisson automorphism groups and in Section 7,  we investigate some partial differential equations involving the Jacobian determinant. %\red{Add section 6 and section 7.}

\section{Preliminaries}
\label{xxsec1}

In this section, we will review some basic concepts  
related to $n$-Lie algebras, Nambu-Poisson algebras, and 
valuations. First, we recall the definition of an $n$-Lie 
algebra.

\begin{definition}\cite{Fi}
\label{xxdef1.1}
Let $V$ be a vector space equipped with an $n$-ary multilinear 
operation $\pi:= \{-,\cdots,-\}: V^{\otimes n}\to V$. We say 
$V$ is an {\it $n$-Lie algebra} if the following hold:
\begin{enumerate}
\item[(1)]
$\pi$ is skew-symmetric, namely, $\{v_1,\cdots,v_n\}=0$
whenever $v_i=v_j$ for some $i\neq j$. Since we assume 
${\text{char}}\;\Bbbk=0$ (or over any field of characteristic
$\neq 2$), this is equivalent to
$$\{v_{\sigma(1)},\cdots,v_{\sigma(n)}\}=
{\rm{sgn}}(\sigma) \{v_1,\cdots,v_n\},$$
for all $\sigma\in S_n$, and 
\item[(2)]
the following Jacobi identity holds
{\small 
$$\{\{u_1,\cdots,u_n\},v_1,\cdots,v_{n-1}\}
=\sum_{s=1}^{n}
\{u_1,\cdots,u_{s-1},
\{u_{s},v_1,\cdots,v_{n-1}\},
u_{s+1},\cdots,u_{n}\},$$}
for all $u_{s},v_{t}\in V$.
\end{enumerate}
\end{definition}

These algebras, for general $n\geq 2$, were first introduced by 
Filippov \cite{Fi} in 1985. Before that, in 1973, Nambu \cite{Nam} 
studied a dynamical system which was defined as a Hamiltonian 
system concerning a ternary (or $3$-Lie) Poisson bracket. 
In 1995, Takhtajan \cite{Ta} reconsidered the subject, proposed 
a general, algebraic definition of a Nambu-Poisson bracket of 
order $n$, and gave the basic characteristic properties of this 
operation. In 2021, Umirbaev-Zhelyabin \cite{UZ} re-introduced 
an equivalent notion, called $n$-Lie Poisson algebra, and 
proved several very nice results about Nambu-Poisson 
algebras.

\begin{definition} \cite{Nam, Ta, UZ}
\label{xxdef1.2}
A commutative algebra $P$ is called a {\it Nambu-Poisson
algebra} or an {\it $n$-Lie Poisson algebra} if it is 
equipped with an $n$-ary multilinear operation 
$\pi:= \{-,\cdots,-\}: P^{\otimes n}\to P$ such that
\begin{enumerate}
\item[(1)]
$(P,\pi)$ is an $n$-Lie algebra, and 
\item[(2)]
$\pi$ satisfies the following Leibniz rule for every 
$1\leq i\leq n$
$$\begin{aligned}
\{v_1, \cdots, v_{i-1}, a b, v_{i+1},\cdots,v_n\}
=&a \{v_1, \cdots,v_{i-1}, b, v_{i+1},\cdots,v_n\}\\
& \quad + b \{v_1, \cdots, v_{i-1}, a, v_{i+1},\cdots,v_n\},
\end{aligned}
$$
for all $v_{s},a,b\in P$.
\end{enumerate}
\end{definition}

Note that a $2$-Lie Poisson algebra in the classical sense agrees with a Poisson 
algebra. We refer to the book 
\cite{LPV} for basic definitions of Poisson algebras.
Many well-known Poisson fields have an $n$-Lie version; see the
following example.

\begin{example}
\label{xxex1.3}
Let $\Omega\in \Bbbk^{[n+1]}(:=\Bbbk[t_1,\cdots,t_{n+1}])$ and 
$A_{\Omega}$ be the $n$-Lie Poisson algebra defined in 
Construction \ref{xxcon0.5}.
\begin{enumerate}
\item[(1)]
Let $q\in \Bbbk^{\times}$ and $\Omega=q t_1\cdots t_{n+1}$. 
Let $N_q$ be the subfield of $Q(A_{\Omega})$ generated by 
$t_1,\cdots,t_n$ for this particular $\Omega$. It is clear 
that $N_q$ is generated by $t_1,\cdots,t_n$ as a field with 
its $n$-Lie Poisson algebra structure determined by 
$\{t_1,\cdots,t_n\}=qt_1\cdots t_n$. We call $N_q$ the 
{\it $q$-skew Nambu-Poisson field} or {\it $q$-skew $n$-Lie 
Poisson field}.
\item[(2)]
In general, let $t_1,\cdots,t_n$ be nonzero elements in an 
$n$-Lie Poisson field $N$. Let $(a_{ij})_{n\times n}$ be 
an invertible $n\times n$-matrix over ${\mathbb Z}$. Let 
$y_i=\prod_{j=1}^{n} t_j^{a_{ij}}\in N$ for $1\leq i\leq n$. 
An easy computation shows that
\begin{equation}
\label{E1.3.1}\tag{E1.3.1}
\{y_1,\cdots, y_n\} y_1^{-1}\cdots y_{n}^{-1}
=\det[(a_{ij})_{n\times n}] 
\{t_1,\cdots,t_n\} t_1^{-1}\cdots t_{n}^{-1}.
\end{equation}
Using this fact, one sees that if $N_q=\Bbbk(t_1,\cdots, t_n)$ 
is $q$-skew $n$-Lie Poisson field, then 
$Q:=\Bbbk(y_1,\cdots,y_n)$ is a Nambu-Poisson subfield of $N_q$.
It is easy to see that $Q\cong N_{q'}$ where $q'=q
\det[(a_{ij})_{n\times n}]$.
\item[(3)]
Let $\Omega=t_{n+1}$. Let $N_{Weyl}$ be the subfield of 
$Q(A_{\Omega})$ generated by $t_1,\cdots,t_n$ for this 
particular $\Omega$. It is clear that $N_{Weyl}$ is generated 
by $t_1,\cdots,t_n$ as a field with its $n$-Lie Poisson algebra 
structure determined by $\{t_1,\cdots,t_n\}=1$.
We call $N_{Weyl}$ the {\it Weyl Nambu-Poisson field} or 
{\it $n$-Lie Weyl Poisson field}.
\end{enumerate}
\end{example}

\begin{definition}
\label{xxdef1.4}
Let $\OO$ be an abelian group. We say $\OO$ is an ordered group 
if it is equipped with a total order $<$ such that if $a<b$ 
then $a+c<b+c$ for all $c\in \OO$.
\end{definition}

For example, ${\mathbb Z}$ is an ordered group. There are several 
ways to make ${\mathbb Z}^{m}$ an ordered group when $m\geq 2$. 
One is the following: $(x_1,\cdots, x_m)<(y_1,\cdots,y_m)$ if 
and only if either (i) $\sum_i x_i <\sum_i y_i$ or (ii) 
$\sum_i x_i =\sum_i y_i$ and there is an $s$ between $1$ 
and $m$ such that $x_i=y_i$ for all $i<s$ and 
$x_{s}< y_{s}$.

\begin{convention}
\label{xxcon1.5} Throughout, we assume the following.
\begin{enumerate}
\item[(1)]
$\OO$ is a nontrivial abelian ordered group. 
\item[(2)]
$w$ is an element in $\OO$.
\item[(3)]
For every element $d\in \OO$, $d<\infty$ and 
$d+\infty=\infty=\infty+\infty$.
\end{enumerate}
\end{convention}

Next, we recall the definition of a valuation 
\cite[Definition 1.1]{HTWZ3}.

\begin{definition} 
\label{xxdef1.6}
Let $N$ be a commutative algebra over $\Bbbk$. For parts 
{\rm{(3)-(6)}}, let $N$ be a Nambu-Poisson algebra.
\begin{enumerate}
\item[(1)]
An {\it $\OO$-valued valuation} {\rm{(}}or simply a
{\it valuation}{\rm{)}} on $N$ is a map
$$\nu: N \to \OO \cup\{\infty\}$$
which satisfies the following properties: for all $a, b\in N$,
\begin{enumerate}
\item[(a)]
$\nu(a)=\infty$ if and only if $a=0$,
\item[(b)]
$\nu(a)=0$ for all $a\in \Bbbk^{\times}:=\Bbbk\setminus \{0\}$,
\item[(c)]
$\nu(ab)=\nu(a)+\nu(b)$ (assuming $n+\infty=\infty$ when $n\in
\OO\cup\{\infty\}$),
\item[(d)]
$\nu(a+b)\geq \min\{\nu(a),\nu(b)\}$, with equality if $\nu(a)
\neq \nu(b)$.
\end{enumerate}
\item[(2)]
A valuation $\nu$ is called {\it trivial} if $\nu(a)=0$ for 
all $a\in N\setminus\{0\}$. Otherwise, it is called {\it nontrivial}.
\item[(3)]
Let $w\in \OO$. A valuation $\nu$ is called an 
{\it $(\OO,w)$-valuation} 
{\rm{(}}or simply a {\it $w$-valuation}{\rm{)}} if
\begin{enumerate}
\item[(e)]
$\nu(\{b_1,\cdots,b_n\})\geq -w +\sum_{s=1}^n\nu(b_{s})$ 
for all $b_i\in N$.
\end{enumerate}
Note that a $({\mathbb Z},0)$-valuation is just a 
{\it $0$-valuation} given in Definition \ref{xxdef0.1}. 
%When $\nu$ is a $w$-valuation, we say $w$ is the {\it weight} of $\nu$ and write $w(\nu)=w$.
\item[(4)]
A $w$-valuation $\nu$ is called {\it classical}
if $\nu(\{b_1,\cdots, b_n\})>\sum_{s=1}^n \nu(b_{s})-w$ 
for all $b_i\in N$.
\item[(5)]
Let ${\mathcal V}_{w}(N)$ be the set of nontrivial 
$w$-valuations on $N$. 
\item[(6)]
Let $\nu_1$ (resp. $\nu_2$) be a $(\OO,w_1)$-valuation
(resp. $(\OO,w_2)$-valuation) on $N$. We say $\nu_1$ and 
$\nu_2$ are {\it equivalent} if there are positive integers
$a,b$ such that $a\nu_1=b\nu_2$. In this case, we write 
$\nu_1\sim \nu_2$. Otherwise we say $\nu_1$ and $\nu_2$ 
are {\it non-equivalent}.
\end{enumerate}
\end{definition}

By convention, a $w$-valuation is defined on a Nambu-Poisson 
algebra and a valuation can be defined on a commutative algebra 
(without considering the Nambu-Poisson structure). Similar 
to Definition \ref{xxdef1.6}(5), we define
$${\mathcal V}_{cw}(N):={\text{the set of nontrivial Poisson
classical $w$-valuations on $N$}}$$
and
$${\mathcal V}_{ncw}(N):={\text{the set of nontrivial Poisson
non-classical $w$-valuations on $N$}}.$$
It is clear that ${\mathcal V}_{w}(N)={\mathcal V}_{cw}(N)
\sqcup {\mathcal V}_{ncw}(N)$.

The following lemma is similar to \cite[Lemma 1.3]{HTWZ3}.
Let $A$ be an $\OO$-graded domain. Define the {\it support} 
of $A$ by
\begin{equation}
\label{E1.6.1}\tag{E1.6.1}
\Supp (A):=\{i\in \OO\mid A_i\neq 0\}.
\end{equation}
The following lemma is easy, and its proof is omitted.

\begin{lemma}
\label{xxlem1.7}
Suppose $\Bbbk$ is algebraically closed. Let $A$ be an 
$\OO$-graded domain. 
\begin{enumerate}
\item[(1)]
If $A$ is an $\OO$-graded field, then $\Supp(A)$ is a 
subgroup of $\OO$. 
\item[(2)] 
Suppose $\GKdim \Bbbk \Supp(A)=\GKdim A<\infty$. Then
$\dim_{\Bbbk} A_i$ is either 0 or 1.
\end{enumerate}
\end{lemma}

\begin{lemma}
\label{xxlem1.8}
Let $N\subseteq Q$ be two $n$-Lie Poisson fields.
\begin{enumerate}
\item[(1)]
Suppose $Q$ is algebraic over $N$.
There is a natural map induced by restriction
$${\mathcal V}_{w}(Q)\to {\mathcal V}_{w}(N).$$
As a consequence, if ${\mathcal V}_{w}(Q)\neq \emptyset$,
then ${\mathcal V}_{w}(N)\neq \emptyset$. 
\item[(2)]
For every $w'<w$,
$${\mathcal V}_{w'}(N)\subseteq {\mathcal V}_{cw}(N)
\subseteq {\mathcal V}_{w}(N).$$
As a consequence, if ${\mathcal V}_{w}(N)={\mathcal V}_{ncw}(N)$,
then ${\mathcal V}_{w'}(N)=\emptyset$ for all $w'<w$.
\item[(3)]
Let $d$ be a positive integer. Then the assignment
$\nu\to d \nu$ defines an injective map ${\mathcal V}_{w}(N)
\to {\mathcal V}_{dw}(N)$.
\item[(4)]
Part {\rm{(3)}} holds for ${\mathcal V}_{cw}(N)$ and 
${\mathcal V}_{ncw}(N)$.
\end{enumerate}
\end{lemma}

\begin{proof} 
(1) Let $\nu\in {\mathcal V}_{w}(Q)$. Let $\phi:=\nu\mid_{N}$.
It is clear that $\phi$ is a $w$-valuation on $N$.
It remains to show that $\phi$ is nontrivial. Suppose to
the contrary that $\phi$ is trivial. So $\phi(f)=0$
for all $0\neq f\in N$. We claim that if $0\neq q\in Q$,
then $\nu(q)=0$. If not, we can assume $\nu(q)<0$ (after
replacing $q$ by $q^{-1}$ if necessarily).
Since $Q$ is algebraic over $N$,
there are $a_i\in N$ with $a_0\neq 0$ such that
$q^{n+1} =a_0+a_1 q+\cdots+ a_{n} q^{n}$. Then
$$(n+1)\nu(q)=\nu(q^{n+1})=
\nu(\sum_{i=0}^{n} a_i q^i)
\geq \min_{i=0}^{n}\{ \nu(a_i q^{i})\}
\geq n\nu(q),$$ 
which contradicts the assumption 
$\nu(q)<0$. Therefore, we proved the claim.
Since $\nu$ is nontrivial, we obtain a contradiction.
The consequences are clear.

(2) Let $\nu$ be a $w'$-valuation. Then
$$\nu(\{a_1,\cdots,a_n\})\geq \sum_s \nu(a_s)-w'>\sum_s \nu(a_s)-w.$$
By definition, $\nu$ is a classical $w$-valuation 
on $N$. The main assertion follows. If 
${\mathcal V}_{w}(N)={\mathcal V}_{ncw}(N)$,
then ${\mathcal V}_{cw}(N)=\emptyset$. The consequence 
follows from the main assertion.

(3,4) Clear.
\end{proof}

A morphism $\phi: \OO\to \OO'$ between two ordered abelian
groups is called {\it order-preserving} if $w\leq w'$ in $\OO$
implies that $\phi(w)\leq \phi(w')$ in $\OO'$. Note that we do 
not assume that $w<w'$ in $\OO$ implies that $\phi(w)< 
\phi(w')$ in $\OO'$.

\begin{lemma}
\label{xxlem1.9}
Let $N$ be an $n$-Lie Poisson algebra. If $\phi: \OO\to \OO'$ 
is an order-preserving morphism between two ordered abelian 
groups, then, for every $(\OO,w)$-valuation $\nu$ on $N$,
$\phi\circ \nu$ is an $(\OO',\phi(w))$-valuation. As a 
consequence, if $\phi$ is an isomorphism, then $\phi$ induces a map
${\mathcal V}_{w}(N)\to {\mathcal V}_{\phi(w)}(N)$.
\end{lemma}

\begin{proof} This is an easy verification.
\end{proof}

The following lemma is an $n$-Lie Poisson version of 
\cite[Lemma 1.6]{HTWZ3}, that will be used in later 
sections.

\begin{lemma}
\label{xxlem1.10}
Let $N$ be an $n$-Lie Poisson field and $x_1,\cdots,x_n\in N$.
Let $X$ be the subfield generated by $x_1,\cdots,x_n$.
\begin{enumerate}
\item[(1)]
If $\GKdim X<n$, then $\{x_1,\cdots, x_n\}=0$.
\item[(2)]
Suppose $N$ has GK-dimension $n$ with nontrivial
$n$-Lie Poisson bracket. If $\{x_1,\cdots,x_n\}=0$, 
then $\GKdim X<n$.
\end{enumerate}
\end{lemma}

\begin{proof}
(1) Let $S:=\langle x_1,\cdots,x_n\rangle$ denote the set 
consisting of elements $x_1,\cdots,x_n$. Since $\GKdim X<n$, 
elements in $S$ are not algebraically independent. Pick a 
largest subset of $S$, say $S':=\langle x_1,\cdots,x_{j}\rangle$
for some $j<n$, whose elements are algebraically independent. 
Then $\GKdim X=j$. Let $X'$ be the subfield generated by $S'$. 
Then $\GKdim X=\GKdim X'$ and every element in $X$ is integral 
over $X'$.

For the subset $\bfx:=\langle x_1,\cdots,x_{n-1}\rangle \subset N$, 
let 
$$C_{\bfx}(N):=\{f\in N\mid \{x_1,\cdots, x_{n-1},f\}=0\}.$$ 
Since ${\text{char}}\; \Bbbk=0$, $C_{\bfx}(N)$ is a subfield 
of $N$ that is integrally closed in $N$ (namely, $C_{\bfx}(N)$ 
is equal to its integral closure in $N$). To see this, let 
$d(-)$ be the derivation $\{x_1,\cdots, x_{n-1},-\}$ and let 
$y\in N$ be integral over $C_{\bfx}(N)$. It suffices to show 
that $d(y)=0$. Let $h(t):=t^n+a_1 t^{n-1}+\cdots +a_n$ with 
coefficients $a_i$ in $C_{\bfx}(N)$ be the minimal polynomial 
of $y$. Then $0=d(0)=d(h(y))=h'(y) d(y)$. Since $\Bbbk$ is of 
characteristic zero, $h'(y)\neq 0$. This forces that $d(y)=0$ 
as required.

By definition, $x_i\in C_{\bfx}(N)$ for $i=1,\cdots, n-1$, and 
whence $X'\subseteq C_{\bfx}(N)$. Since $x_n$ is integral over 
$X'$, $x_n$ is integral over $C_{\bfx}(N)$. Hence $x_n\in 
C_{\bfx}(N)$ by the previous paragraph. This means that 
$\{x_1,\cdots, x_{n-1},x_n\}=0$.

(2) Suppose to the contrary that $\GKdim X=n$. In this case,
$X$ is a rational function field and $N$ is integral over $X$. 
By assumption, $\{x_1,\cdots,x_n\}=0$, which implies that 
$\{f_1,\cdots,f_n\}=0$ for all $f_i\in X$. 

We claim that, for every $i$, 
$$\{g_1,\cdots,g_i,f_{i+1},\cdots,f_n\}=0$$ 
for all $g_{s} \in N$ and $f_{t}\in X$. By the last paragraph, 
the claim holds for $i=0$. Suppose now the claim holds for 
$i\geq 0$. We want to prove it for $i+1$. Fix $g_1,\cdots,g_i,f_{i+2},
\cdots,f_n$ and let
$$C(N)=\{f\in N\mid 
\{g_1,\cdots,g_i,f,f_{i+1},\cdots,f_n\}=0\}.$$
By the induction hypothesis, $X\subseteq C(N)$.
By the proof of part (1), $C(N)$ is integrally closed
in $N$. Hence $C(N)=N$. This proves the claim 
for $i+1$. It follows by induction and the case
when $i=n$ that $N$ has trivial $n$-Lie Poisson structure.
This yields a contradiction.
\end{proof}

An ideal $I$ of an $n$-Lie Poisson algebra $P$ is called an 
{\it $n$-Lie Poisson ideal}, or {\it Nambu-Poisson ideal} or 
just {\it Poisson ideal} if 
$$\{I, P,\cdots, P\}\subseteq I.$$ 
A Poisson ideal $I$ is called {\it Poisson prime} if, for any 
two Poisson ideals $J_1, J_2$, $J_1 J_2\subseteq I$ implies 
that $J_1\subseteq I$ or $J_2\subseteq I$. The following 
result was proved by Goodearl \cite{Go1} in a more general 
setting of a commutative differential $\kk$-algebra.

\begin{lemma} \cite[Lemma 1.1(d)]{Go1}
\label{xxlem1.11}
If $P$ is a noetherian Nambu-Poisson algebra, then an ideal
$I$ of $P$ is a Poisson prime if and only if it is both
Poisson and prime. 
\end{lemma}

\section{Generalization of material from \cite{HTWZ3}}
\label{xxsec2}

Most of this section is a generalization of 
corresponding material in \cite{HTWZ3}. It is necessary
for this paper. However, to save space, some proofs are 
omitted, provided these proofs are either easy or similar 
to ones given in \cite{HTWZ3}.

Fix an integer $n\geq 3$ (or even $n=2$) 
and a Nambu-Poisson algebra stands for an $n$-Lie Poisson 
algebra. First, we recall the notion of a Nambu-Poisson $w$-graded 
algebra and a slightly weaker version of the valuation 
called a {\it filtration} of an algebra. As always, $w$ is 
an element in $\OO$ unless otherwise stated. 
By abusing the notation, let $0$ denote the 
group unit of the abelian group $\OO$.

\begin{definition}
\label{xxdef2.1}
A Nambu-Poisson algebra $A$ is called {\it $w$-graded} if
\begin{enumerate}
\item[(1)]
$A=\oplus_{i\in \OO} A_i$ is an $\OO$-graded
commutative algebra, and
\item[(2)]
$\{A_{i_1}, A_{i_2},\cdots, A_{i_n}\}\subseteq 
A_{(\sum_s i_s)-w}$ for all $i_1,\cdots,i_n \in \OO$.
\end{enumerate}
\end{definition}

Let $A_{\Omega}$ and $P_{\Omega}$ be the Nambu-Poisson
algebras defined in Construction \ref{xxcon0.5}. Suppose 
that $\deg t_i=1$ for all $i$ and that $\Omega$ is homogeneous. 
Then these are $w$-graded where $w=-(\deg \Omega-n-1)$. 

\begin{definition}
\label{xxdef2.2}
Let $A$ be an algebra. Let
${\mathbb F}:=\{F_{i}\mid i\in \OO\}$
be a chain of $\Bbbk$-subspaces of $A$. 
\begin{enumerate}
\item[(1)] \cite[pp. 173-174]{Zh2}
We say ${\mathbb F}$ is an {\it $\OO$-filtration} {\rm{(}}or
simply a {\it filtration}{\rm{)}} of $A$ if it satisfies
\begin{enumerate}
\item[(a)]
$1 \in F_0\setminus F_{>0}$ where
$F_{>0}:=\bigcup_{j>0} F_j$ and $F_{i}\supseteq F_{i'}$ 
for all $i<i'$ in $\OO$, 
\item[(b)]
$F_i F_j\subseteq F_{i+j}$ for all $i,j\in \OO$,
\item[(c)]
$\bigcup_{i\in \OO} (F_i\setminus F_{>i})=A\setminus \{0\}$,
where $F_{>i}:=\bigcup_{j>i} F_j$.
\end{enumerate}
\item[(2)]
Suppose that $A$ is a Nambu-Poisson algebra and that 
${\mathbb F}$ is a filtration of $A$. If further
\begin{enumerate}
\item[(d)]
$\{F_{i_1},\cdots, F_{i_n}\}\subseteq F_{(\sum_s i_s)-w}$ 
for all $i_1,\cdots,i_n\in \OO$,
\end{enumerate}
then ${\mathbb F}$ is called a {\it $w$-filtration} of $A$.
\end{enumerate}
\end{definition}

The {\it associated graded ring} of a filtration ${\mathbb F}$ 
of $A$ is defined to be
$$\gr_{\mathbb F} A:= \bigoplus_{i\in\OO} F_i/F_{>i}$$
which is an $\OO$-graded algebra. By 
\cite[Theorem 4.3(1)]{Zh2}, 
\begin{equation}
\label{E2.2.1}\tag{E2.2.1}
\GKdim \gr_{\mathbb F} A\leq \GKdim A.
\end{equation} 
For any element $a\in F_i$, let ${\overline{a}}$ denote the 
element $a+F_{>i}$ 
in the $i$th degree component $(\gr_{\mathbb F} A)_i:=F_i/F_{>i}$. 
Next, we add the Nambu-Poisson structure to the picture.

\begin{lemma}
\label{xxlem2.3}
Suppose ${\mathbb F}$ is a $w$-filtration of a Nambu-Poisson 
algebra $A$. Then $\gr_{\mathbb F} A$ is a Nambu-Poisson 
$w$-graded algebra.
\end{lemma}

\begin{proof} 
It is well-known that $\gr_{\mathbb F} A$ is a graded algebra
\cite[Section 4]{Zh2}. Its addition and multiplication are 
defined as follows:
${\overline{f}}+{\overline{g}}=\overline{f+g}$ when
${\overline{f}},{\overline{g}}\in (\gr_{\mathbb F} A)_i$
and ${\overline{f}}{\overline{g}}=\overline{fg}$ when
${\overline{f}}\in (\gr_{\mathbb F} A)_i$ and
${\overline{g}}\in (\gr_{\mathbb F} A)_j$.

Next, we define the Nambu-Poisson $w$-graded structure as 
follows. Let $\{\overline{f}_s\}_{s=1}^n$ be elements in
$(\gr_{\mathbb F} A)_{i_s}$ for $s=1,\cdots,n$ respectively
where $f_s\in F_{i_s}$. By definition,
$\{f_1,\cdots,f_n\}\in F_{\sum_s i_s-w}$. Let 
$\overline{\{f_1,\cdots,f_n\}}$ be the
class in $(\gr_{\mathbb F} A)_{\sum_s i_s-w}$. We define 
the Nambu-Poisson bracket on $\gr_{\mathbb F} A$ by 
$$\{{\overline{f}_1},\cdots,{\overline{f}_n}\}
:= \overline{\{f_1,\cdots,f_n\}}\in (\gr_{\mathbb F} A)_{\sum_s i_s-w}$$
for all ${\overline{f}_s}\in (\gr_{\mathbb F} A)_{i_s}$. 
It is easy to see that $\{-,\cdots,-\}$ is well-defined 
(or is independent of the choices of preimages $f_s$).

We claim that $\{-,\cdots,-\}$ is a Nambu-Poisson bracket. 
Let $a_s\in F_{i_s}$ for $s=1,\cdots,n-1$, $b\in F_j, c\in F_k$. 
Then in $(\gr_{\mathbb F} A)_{(\sum_{s<n}i_s)+j+k-w}$, we have
$$\begin{aligned}
\{\overline{a_1},\cdots,\overline{a_{n-1}}, 
\overline{b}\overline{c}\}
&=\{\overline{a_1},\cdots,\overline{a_{n-1}}, \overline{bc}\}\\
&=\overline{\{a_1,\cdots,a_{n-1},bc\}}\\
&=\overline{\{a_1,\cdots,a_{n-1},b\}c+\{a_1,\cdots,a_{n-1},c\}b}\\
&=\overline{\{a_1,\cdots,a_{n-1},b\}c}+\overline{\{a_1,\cdots,a_{n-1},c\}b}\\
&=\overline{\{a_1,\cdots,a_{n-1},b\}}\overline{c}
  +\overline{\{a_1,\cdots,a_{n-1},c\}}\overline{b}\\
&=\{\overline{a_1},\cdots,\overline{a_{n-1}},\overline{b}\}\overline{c}
  +\{\overline{a_1},\cdots,\overline{a_{n-1}},\overline{c}\}\overline{b}.
\end{aligned}
$$
So $\{\overline{a_1},\cdots,\overline{a_{n-1}}, -\}$ is a 
derivation. Similarly, one checks that the multilinear 
operation $\{-,\cdots,-\}$ is skew-symmetric and satisfies the 
Jacobi identity [Definition \ref{xxdef1.1}(2)]. Therefore 
$\gr_{\mathbb F} A$ is a Nambu-Poisson $w$-graded algebra.
\end{proof}

Given a filtration ${\mathbb F}$ we define the notion of a 
{\it degree} function, denoted by $\deg$, on elements in $A$ by
\begin{equation}
\label{E2.3.1}\tag{E2.3.1}
\deg a:=i, {\text{ if $a\in F_i\setminus F_{>i}$}} \quad
{\text{and}} \quad \deg(0)=+\infty.
\end{equation}
 
Below is an $\OO$-version of \cite[Lemma 2.4]{HTWZ3}. 
Its proof is omitted.

\begin{lemma}
\label{xxlem2.4}
Let ${\mathbb F}$ be a filtration of an algebra $A$
such that $\gr_{\mathbb F} A$ is a domain.
\begin{enumerate}
\item[(1)]
Then $\deg$ satisfies the following conditions for $a,b\in A$,
\begin{enumerate}
\item[(a)]
$\deg(a)\in \OO$ for any
nonzero element $a\in A$, and $\deg a=\infty$ if and only if
$a=0$,
\item[(b)]
$\deg(c)=0$ for all $c\in \Bbbk^{\times}$,
\item[(c)]
$\deg(ab)=\deg(a)+\deg(b)$,
\item[(d)]
$\deg(a+b)\geq \min\{\deg(a),\deg(b)\}$, with
equality if $\deg(a)\neq \deg(b)$.
\end{enumerate}
Namely, $\deg$ is a valuation on $A$.
\item[(2)]
If ${\mathbb F}$ is a $w$-filtration of $A$,
then $\deg$ is a $w$-valuation in the sense of
Definition {\rm{\ref{xxdef1.6}(3)}}.
\item[(3)]
$F_0(A)$ is  integrally closed in $A$.
\end{enumerate}
\end{lemma}

Conversely, if we are given a valuation (or a degree function)  
$$\nu(=\deg): A\to \OO\cup\{\infty\}$$ 
satisfying Lemma \ref{xxlem2.4}(1a,1b,1c,1d), then we can 
define a filtration ${\mathbb F}^{\nu}:=\{F_i^{\nu}\mid i\in \OO\}$
of $A$ by
\begin{equation}
\label{E2.4.1}\tag{E2.4.1}
F_i^{\nu}:=\{ a\in A\mid \nu(a)(=\deg(a))\geq i\}.
\end{equation}
If no confusion occurs, we will delete $^{\nu}$ from 
${\mathbb F}^{\nu}$ and $F_i^{\nu}$.

Here is an $\OO$/Nambu-Poisson version of \cite[Lemma 2.5]{HTWZ3}.
Its proof is omitted.

\begin{lemma}
\label{xxlem2.5}
Let $A$ be a Nambu-Poisson algebra with a valuation $\nu$. 
Then ${\mathbb F}:=\{F_i\}$ defined as in \eqref{E2.4.1} is a
filtration of $A$ such that $\gr_{\mathbb F} A$ is a domain.
If $\nu$ is a $w$-valuation, then $\gr_{\mathbb F} A$
is a Nambu-Poisson $w$-graded domain.
\end{lemma}

A filtration ${\mathbb F}$ is called {\it good} if 
$\gr_{\mathbb F} A$ is a domain. It is clear that there is a 
one-to-one correspondence between the set of valuations on $A$ 
and the set of good filtrations of $A$. In the Nambu-Poisson 
case we have the following version of \cite[Lemma 2.6]{HTWZ3}.
Its proof is omitted.

\begin{lemma}
\label{xxlem2.6}
Let $A$ be a Nambu-Poisson algebra. Then there is a 
one-to-one correspondence between the set of good 
$w$-filtrations of $A$ and the set of $w$-valuations 
on $A$.
\end{lemma}

Suppose $\gr_{\mathbb F} A$ is a domain. The corresponding 
degree function $\deg(=:\nu)$ (see \eqref{E2.3.1}) satisfies all 
conditions in Lemma \ref{xxlem2.4}(1). Let $Q$ be the fraction 
field of $A$. Now we define a {\it degree} on $Q$ by
\begin{equation}
\label{E2.6.1}\tag{E2.6.1}
\deg_Q(ab^{-1}):=\deg a-\deg b \in \OO\cup\{\infty\}
\end{equation}
for $a\in A$ and $b\in A\setminus\{0\}$. It is easy to
check that $\deg_Q$ satisfies all conditions in
Lemma \ref{xxlem2.4}(1). So we can also use $\deg_Q$ to define
a filtration ${\mathbb F}(Q):=\{F_i(Q)\}_{i\in \OO}$ 
of $Q$ by
\begin{equation}
\label{E2.6.2}\tag{E2.6.2}
F_i(Q):=\{a b^{-1} \in Q\mid \deg_Q (ab^{-1})\geq i\}
\quad {\text{for }} i\in \OO.
\end{equation}

When $A$ is an $\OO$-graded domain, the $\OO$-graded fraction 
ring of $A$ is denoted by $Q_{gr}(A)$.

\begin{lemma}
\label{xxlem2.7}
Retain the above notation.
\begin{enumerate}
\item[(1)]
Suppose ${\mathbb F}$ is a good $w$-filtration. 
Then ${\mathbb F}(Q)$ is a good $w$-filtration of
$Q$ and $\deg_Q$ is a $w$-valuation on $Q$. 
\item[(2)]
The associated graded ring $\gr_{{\mathbb F}(Q)} Q$ is
canonically isomorphic to the graded fraction ring
$Q_{gr}(\gr_{\mathbb F} A)$.
\end{enumerate}
\end{lemma}

\begin{proof}
See the proof of \cite[Lemma 2.7]{HTWZ3}.
\end{proof}

If ${\mathbb F}$ is the filtration determined by a valuation
$\nu$, we also use $\gr_{\nu}(A)$ for $\gr_{\mathbb F} (A)$.

\begin{lemma}
\label{xxlem2.8}
Let ${\mathbb F}$ and ${\mathbb G}$ be filtrations of $A$. 
Suppose ${\mathbb G}$ is a subfiltration of ${\mathbb F}$, 
namely, $G_i(A)\subseteq F_i(A)$ for all $i\in \OO$.
\begin{enumerate}
\item[(1)]
There is a natural graded algebra homomorphism 
$\phi:\gr_{\mathbb G} A\to \gr_{\mathbb F} A$.
\item[(2)]
${\mathbb G}={\mathbb F}$ if and only if $\phi$ is injective.
\end{enumerate}
\end{lemma}

\begin{proof} 
(1) See the proof of \cite[Lemma 2.8]{HTWZ3}. 

(2) One direction is clear. For the other direction, suppose
that ${\mathbb G}\neq {\mathbb F}$. Then there is an $f\in
F_j(A)\setminus G_j(A)$ for some $j$. Since ${\mathbb G}$ is 
a filtration, there is an $i\in \OO$ such that 
$f\in G_i(A)\setminus G_{>i}(A)$. Since $f\not\in G_j(A)$ and 
$f\in G_{i}(A)$, we have $i<j$. Consider the algebra map $\phi$ 
defined in the proof of part (1):
$$\phi: (\gr_{\mathbb G}(A))_i\to (\gr_{\mathbb F}(A))_i$$
which sends 
$$\widetilde{f}:=f+ G_{>i}(A)\mapsto \overline{f}:=f+F_{>i}(A).$$
Since $i<j$, $f\in F_{j}(A)\subseteq F_{>i}(A)$. 
This means that $\phi(\widetilde{f})=\overline{f}=0$. By the choice
of $i$, $\widetilde{f}\neq 0$. Therefore, $\phi$ is not injective 
as desired.
\end{proof}

Now, we introduce a general method to construct a filtration. 
Suppose that an algebra $A$ is generated by $\{x_k\}_k$ where
$k$ is in a fixed index set, and that we are given a degree 
assignment on the generators:
\begin{equation}
\label{E2.8.2}\tag{E2.8.2}
\deg(1)=0 \quad{\text{and}} \quad 
\deg(x_k)=d_k\in \OO \quad
{\text{for all}} \ k.
\end{equation}
Then we can define a set of $\Bbbk$-subspaces
${\mathbb F}^{ind}:=\{F_i^{ind}\}_{i\in \OO}$ of $A$ by
\begin{equation}
\label{E2.8.3}\tag{E2.8.3}
F_i^{ind}:={\text{the span of all monomials 
$\prod_{k} x_k^{n_k}$ with $\sum_k n_k d_k\geq i$ }}
\end{equation}
where only finitely many $n_k\neq 0$. It is clear that 
${\mathbb F}^{ind}$ is uniquely determined by the degree 
assignment given in \eqref{E2.8.2}. It is easy to see that 
${\mathbb F}^{ind}$ satisfies Definition \ref{xxdef2.2}(1b) 
and part of Definition \ref{xxdef2.2}(1a). It is not clear if 
${\mathbb F}^{ind}$ satisfies
\begin{equation}
\label{E2.8.4}\tag{E2.8.4}
1\in F^{ind}_0(A)\setminus F^{ind}_{>0}(A),
\end{equation}
and Definition \ref{xxdef2.2}(d)
\begin{equation}
\label{E2.8.5}\tag{E2.8.5}
A\setminus \{0\}=\bigcup_{i\in \OO} 
F^{ind}_i(A)\setminus F^{ind}_{>i}(A).
\end{equation}
However, in this paper, we usually consider the case when 
${\mathbb F}^{ind}$ is a subfiltration of another filtration 
${\mathbb F}$. Hence \eqref{E2.8.4} holds automatically and 
$$\bigcap_{i\in \OO} F^{ind}_i(A)\subseteq \bigcap_{i 
\in \OO} F_i(A)=\{0\},$$
which is a part of \eqref{E2.8.5}. In this case, \eqref{E2.8.5}
holds since $\{x_k\}_k$ is a generating set. When we use 
${\mathbb F}^{ind}$ or say that ${\mathbb F}^{ind}$  
is an {\it induced} filtration determined by the degree 
assignment given in \eqref{E2.8.2}, we always implicitly assume 
both \eqref{E2.8.4} and \eqref{E2.8.5}.

\begin{lemma}
\label{xxlem2.9}
Suppose that $A$ is generated by a set $\{x_k\}_k$ and that we 
are given a degree assignment $\deg$ as in \eqref{E2.8.2}. 
Write $\nu=\deg$. Let ${\mathbb F}$ be a filtration of $A$.
\begin{enumerate}
\item[(1)]
Let ${\mathbb F}^{ind}$ be the induced filtration defined by 
\eqref{E2.8.3}. Then $\gr_{{\mathbb F}^{ind}} A$ is generated 
by $\{\widetilde{x_k}\}_k$ where 
$\widetilde{x_k}:=x_k+F^{ind}_{>\nu(x_k)}(A)$.
\item[(2)]
Suppose that $\deg$ is the degree function corresponding to 
${\mathbb F}$ via \eqref{E2.3.1} and that ${\mathbb F}^{ind}$ 
be defined by using \eqref{E2.8.2}-\eqref{E2.8.3} based on 
the given generating set $\{x_k\}_k$. Then, there is a natural 
$\OO$-graded algebra homomorphism 
\begin{equation}
\label{E2.9.1}\tag{E2.9.1}
\phi^{ind}: \gr_{{\mathbb F}^{ind}} A\to \gr_{\mathbb F} A
\end{equation}
determined by $\phi^{ind}(\widetilde{x_k})=\overline{x_k}$
for all $k$, where $\overline{x_k}=x_k+F_{>\nu(x_k)}(A)$.
\item[(3)]
Suppose that $A$ is a Nambu-Poisson algebra and that 
${\mathbb F}^{ind}$ is the induced filtration defined as 
\eqref{E2.8.3}. Then ${\mathbb F}^{ind}$ is a $w$-
filtration if and only if 
\begin{equation}
\label{E2.9.2}\tag{E2.9.2}
\deg(\{x_{k_1},\cdots, x_{k_n}\})\geq \sum_{s=1}^{n} 
\deg(x_{k_s})-w
\end{equation} 
for all $k_s$.
\end{enumerate}
\end{lemma}

\begin{proof}
See the proof of \cite[Lemma 2.9]{HTWZ3}.
\end{proof}

The following result is a crucial lemma, which will be used 
multiple times to calculate $w$-valuation. We choose
to provide a proof although it is similar to
the proof of \cite[Lemma 2.11]{HTWZ3}. 

\begin{lemma}
\label{xxlem2.10}
Let $A$ be a domain of finite GK-dimension, say $d$, and 
generated by $\{x_k\}_{k\in S}$ for an index set $S$. Let 
$\nu$ be a valuation on $A$ and ${\mathbb F}$ be the 
filtration associated to $\nu$. Let ${\mathbb F}^{ind}$ be 
the induced filtration determined by \eqref{E2.8.3} and 
the degree assignment $\deg(x_k):=\nu(x_k)\in \OO$ for 
all $k\in S$. Let $I$ be the image of $\phi^{ind}$, or 
the subalgebra of $\gr_{\mathbb F}(A)$ generated by 
$\{{\overline {x_k}}\}_k$ where $\overline{x_k}$ is 
$x_k+F_{>\nu(x_k)}(A)$ for all $k\in S$. 
\begin{enumerate}
\item[(1)]
Suppose
\begin{enumerate}
\item[(1a)]
$\GKdim I=d$,
\item[(1b)]
$\gr_{{\mathbb F}^{ind}}(A)$ is a domain.
\end{enumerate}
Then ${\mathbb F}$ agrees with the filtration ${\mathbb F}^{ind}$.
As a consequence, there are natural isomorphisms of $\OO$-graded
algebras $\gr_{{\mathbb F}^{ind}}(A)\cong I=\gr_{\mathbb F}(A)$. 
\item[(2)]
Suppose
\begin{enumerate}
\item[(2a)]
all $\{\nu(x_k)\}_k$ are zero,
\item[(2b)]
$\GKdim I=d$.
\end{enumerate}
Then, $\nu$ is trivial.
\end{enumerate}
\end{lemma}

\begin{proof}
(1) Let ${\mathbb F}^{p}$ be the induced filtration 
${\mathbb F}^{ind}$ determined by \eqref{E2.8.3} and the 
degree assignment $\deg(x_k):=\nu(x_k)$ for all $k\in S$. 
By definition, ${\mathbb F}^{p}$ is a subfiltration of
${\mathbb F}$. By Lemma \ref{xxlem2.9}(2), there is an 
algebra homomorphism $\phi^{ind}$. By Lemma 
\ref{xxlem2.8}(2) it suffices to show that $\phi^{ind}$ is 
injective. By definition, $I$ is the image 
of the map $\phi^{ind}$, so we have 
$$\phi^{ind}: \gr_{{\mathbb F}^p} A\to I\subseteq \gr_{\mathbb F} A.$$
Suppose to the contrary that $\phi^{ind}$ is not injective. 
Since $\gr_{{\mathbb F}^p} A$ is a domain, by \eqref{E2.2.1}, we 
obtain that
$$\GKdim I < \GKdim \gr_{{\mathbb F}^p} A \leq \GKdim A=d$$
which contradicts (1a). Therefore, $\phi^{ind}$ is 
injective as required. The consequence is clear.

(2) Since all $\nu(x_k)=0$, $F_0(A)=A$ and $I$ is a subalgebra
of $(\gr_{\mathbb F} A)_0=F_0(A)/F_1(A)$. Since $\GKdim I=d=
\GKdim A=\GKdim F_0(A)$, we have 
$$\GKdim F_0(A)/F_{>0}(A)=\GKdim (\gr_{\mathbb F} A)_0
\geq \GKdim I= \GKdim F_0(A).$$
This implies that $F_{>0}(A)=0$, and thus the assertion
follows.
\end{proof}

\begin{lemma}
\label{xxlem2.11}
Let $A$ be a Nambu-Poisson domain generated by $\{x_k\}_{k\in S}$ 
for an index set $S$. Let $\nu$ be a $w$-valuation on 
$A$ and ${\mathbb F}$ be the filtration associated to $\nu$. 
\begin{enumerate}
\item[(1)]
Let ${\mathbb F}^{ind}$ be the induced filtration of $A$ 
determined by $\deg(x_k):=\nu(x_k)\in \OO$ for all $k$. 
If ${\mathbb F}^{ind}$ is a $w$-filtration, then 
$\phi^{ind}: \gr_{{\mathbb F}^{ind}} A\to \gr_{{\mathbb F}} A$ 
is a Nambu-Poisson algebra homomorphism.
\item[(2)]
Suppose that $A$ is a Nambu-Poisson $w$-graded algebra 
generated by homogeneous elements $\{x_k\}_{k\in S}$ with 
the degree assignment $\deg(x_k):=\nu(x_k)\in \OO$ for all 
$k\in S$. Then $\phi^{ind}: A\to \gr_{{\mathbb F}} A$ is an 
$n$-Lie Poisson algebra homomorphism.
\end{enumerate}
\end{lemma}

\begin{proof} See the proof of \cite[Lemma 2.12]{HTWZ3}.
\end{proof}

Let $P$ be a Nambu-Poisson $w$-graded domain. We fix this 
internal $\OO$-graded structure on $P$ which is called an 
{\it Adams grading}. We define the {\it Adams$^{Id}$ filtration} 
of $P$, denoted by ${\mathbb F}^{Id}$, determined by
\begin{equation}
\label{E2.11.1}\tag{E2.11.1}
F^{Id}_i(P):=\oplus_{n\geq i} P_n 
\quad {\text{for all $i\in \OO$}}.
\end{equation}
It is clear that ${\mathbb F}^{Id}$ is a Poisson
$w$-filtration such that the associated graded ring 
$\gr_{{\mathbb F}^{Id}} P$ is canonically isomorphic to $P$. 

Let $\Aut_{Ab}(\OO)$ be the automorphism group of 
the abelian group $\OO$ without considering the order.
Let $\zeta$ be an element in $\Aut_{Ab}(\OO)$. 
The {\it Adams$^{\zeta}$ filtration} of 
$P$, denoted by ${\mathbb F}^{\zeta}$, is defined by
$$F^{\zeta}_i(P):=\oplus_{n\geq i} P_{\zeta^{-1}(n)}
\quad {\text{for all $i\in \OO$}}.$$
It is clear that ${\mathbb F}^{\zeta}$ is a Poisson
$\zeta(w)$-filtration such that the associated graded
ring $\gr_{{\mathbb F}^{\zeta}} P$ is isomorphic to 
$P$ with grading shifted ($i\leftrightarrow \zeta^{-1}(i)$). 
Let $\nu^{Id}$ (resp. $\nu^{\zeta}$) denote the 
$w$-valuation (resp. $\zeta(w)$-valuation) 
on $P$ associated to ${\mathbb F}^{Id}$ (resp.
${\mathbb F}^{\zeta}$). We call both $\nu^{Id}$
and $\nu^{\zeta}$ the {\it Adams valuations} on
$P$. The following lemma is easy, and its proof is 
omitted.

\begin{lemma}
\label{xxlem2.12}
Let $P$ be an $\OO$-graded domain with $P_i\neq 0$ for all
$i\in \OO$.
\begin{enumerate}
\item[(1)]
$P$ has at least $|\Aut_{Ab}(\OO)|$ many filtrations, namely 
$\{{\mathbb F}^{\zeta}\mid \zeta\in \Aut_{Ab}(\OO)\}$. In these 
cases, the associated graded ring is isomorphic to $P$. If $P$ 
is a Nambu-Poisson $w$-graded domain, then the corresponding 
$\nu^{Id}$ {\rm{(}}resp. $\nu^{\zeta}${\rm{)}} is a $w$-valuation 
{\rm{(}}resp. $\zeta(w)$-valuation{\rm{)}}.
\item[(2)]
If $P$ is a Nambu-Poisson $0$-graded algebra, then the 
Nambu-Poisson field $Q(P)$ has at least $|\Aut_{Ab}(\OO)|$ 
distinct $0$-valuations 
$\{\nu^{\zeta}\mid \zeta\in \Aut_{Ab}(\OO)\}$.
\item[(3)]
Suppose an algebra $A$ has a filtration ${\mathbb F}^{c}$ such 
that $\gr_{{\mathbb F}^c} A$ is canonically isomorphic to $P$ 
as $\OO$-graded algebras. Then $A$ has a valuation associated 
to ${\mathbb F}^{c}$, denoted by $\nu^{c}$. If  further $A$ is 
a Nambu-Poisson algebra and ${\mathbb F}^{c}$ is a 
$w$-filtration, then $\nu^{c}$ is a $w$-valuation. 
\end{enumerate}
\end{lemma}

This lemma leads to the following definition.

\begin{definition}
\label{xxdef2.13}
Let $N$ be a Nambu-Poisson field and $\nu$ a $w$-valuation.
\begin{enumerate}
\item[(1)]
We say $\nu$ is {\it quasi-Adams} if $Q(\gr_{\nu} N)\cong N$.
\item[(2)]
$N$ is called {\it $w$-quasi-Adams} if (a) $N$ admits a 
nontrivial $w$-valuation and (b) every
$w$-valuation is quasi-Adams.
\item[(3)]
We say $\nu$ is {\it Weyl} if $Q(\gr_{\nu} N)\cong N_{Weyl}$
where $N_{Weyl}$ is the Weyl Nambu-Poisson field defined in 
Example \ref{xxex1.3}(2).
\end{enumerate}
\end{definition}

Next, we define some $\gamma$-type invariants, which is 
defined by using filtrations associated with $\nu$.

\begin{definition}
\label{xxdef2.14}
Let $N$ be a Nambu-Poisson field and $w,v$ be in $\OO$. 
\begin{enumerate}
\item[(1)]
The {\it ${^w\Gamma_{v}}$-cap} of $N$ is defined to be
$${^w\Gamma_{v}}(N):=
\bigcap_{\nu\in {\mathcal V}_{w}(N)} F^{\nu}_v(N).$$
where $\nu$ runs over all $w$-valuations of $N$.
\item[(2)]
We say $N$ is {\it ${^w\Gamma_{0}}$-normal} if $C:={^w\Gamma_{0}}(N)$ is 
a Nambu-Poisson subalgebra of $N$ such that $N=Q(C)$ and 
that $C$ is an affine normal domain.
\end{enumerate} 
\end{definition}

It is clear that ${^w\Gamma_{v}}(N)=\{a\in N\mid \nu(a)\geq v,
\forall \; {\text{$w$-valuations on $N$}}\}$. 
The following lemma is straightforward, and its proof 
is omitted.

\begin{lemma}
\label{xxlem2.15}
Let $f: N\to Q$ be a morphism between two Nambu-Poisson fields.
\begin{enumerate}
\item[(1)]
$f$ maps ${^w\Gamma_{v}}(N)$ to ${^w\Gamma_{v}}(Q)$ injectively.
\item[(2)]
${^w\Gamma_{v}}(N)\supseteq\, ^{w'}\Gamma_{v}(N)\supseteq \Bbbk$ for 
all $w<w'$ and all $v$, and ${^w\Gamma_{v}}(N)\supseteq \Bbbk$ 
if $v\leq 0$. 
\end{enumerate}
\end{lemma}

The following is a Nambu-Poisson version of 
\cite[Theorem 4.3]{HTWZ3} and its proof is omitted.

\begin{theorem}
\label{xxthm2.16}
Let $A$ be a Nambu-Poisson domain of GK-dimension at least 
2 and $N$ be the Nambu-Poisson fraction field $Q(A)$. 
Let $w>0$ be an element in $\OO$.
\begin{enumerate}
\item[(1)]
Let $\fp$ be a prime ideal of $A$ of height one such that
$A_{\fp}$ is regular. Then there is a unique nontrivial 
$w$-valuation, denoted by $\nu^{w,\fp}$, such that 
$$\nu^{w,\fp}(x)= i w, 
{\text{  if $x\in \fp^i A_{\fp}\setminus \fp^{i+1} A_{\fp}$}}.$$
\item[(2)]{\rm{(}}Controlling theorem{\rm{)}}
If $A$ is a noetherian normal domain, then 
${^w\Gamma_{0}}(N)\subseteq A$.
\end{enumerate}
\end{theorem}

\begin{remark}
\label{xxrem2.17} We recall some $\alpha$-type invariants 
introduced in \cite{HTWZ3}. This remark works for 
Nambu-algebras/fields $N$, but only for the case when 
$\OO={\mathbb Z}$.
\begin{enumerate}
\item[(1)]
For every $w\in {\mathbb Z}$, ${\mathcal V}_w(N)/\sim$ denote
the set of $w$-valuations $\nu$ such that $\nu$ is not 
equivalent to $a\nu$ for any positive rational number $a<1$. 
We define $\alpha_w(N)$ to be the cardinality of 
${\mathcal V}_w(N)/\sim$.
\item[(2)]
The above theorem makes $\alpha_1(N)$ infinite.
However when $w\neq 1$, $\alpha_{w}(N)$ could be finite.
For example, by Theorem \ref{xxthm3.9}, $\alpha_0(N)=0$ for 
some Nambu-Poisson fields.
\item[(3)]
One can also define $w(N)=\inf\{ w\mid \alpha_w(N)\neq 0\}$. 
Then one can check that $w(N)$ can only be $-\infty$, $0$, or $1$.
See \cite[Remark 2.15]{HTWZ3}.
\end{enumerate}
\end{remark}

Finally, we introduce a concept that will be used in later 
sections.

\begin{definition}
\label{xxdef2.18}
Let $N$ be a Nambu-Poisson field and $w\in \OO$. 
\begin{enumerate}
\item[(1)]
A $w$-valuation $\nu$ on $N$ is called a 
{\it faithful $w$-valuation} if the following hold.
\begin{enumerate}
\item[(a)]
The image of $\nu$ is $\OO \cup\{\infty\}$.
\item[(b)]
$\GKdim \gr_{\nu} N=\GKdim N$.
\item[(c)]
$\nu$ is nonclassical.
\end{enumerate}
Let ${\mathcal V}_{fw}(N)$ be the set of all faithful $w$-valuations
on $N$. 
\item[(2)]
A $w$-valuation $\nu$ on $N$ is called a {\it quasi-faithful $w$-valuation} 
if (b) and (c) in part (1) and the following hold.
\begin{enumerate}
\item[(a')]
$\OO/\OO'$ is torsion over ${\mathbb Z}$ where 
$\OO'$ is the image of $\nu$ on 
$N\setminus \{0\}$.
\end{enumerate}
Let ${\mathcal V}_{qfw}(N)$ be the set of all quasi-faithful $w$-valuations
on $N$. 
\end{enumerate}
\end{definition}

Note that faithful $0$-valuation is exactly the {\it faithful valuation} 
defined in Definition \ref{xxdef0.2}.

To conclude this section, we prove a lemma that will be used
in later sections.

\begin{lemma}
\label{xxlem2.19}
Suppose $\OO={\mathbb Z}^n$. Let $f: K\to Q$ be an embedding 
of Nambu-Poisson fields of GK-dimension $n$. Then 
$f$ induces a map from ${\mathcal V}_{qfw}(Q)\to
{\mathcal V}_{qfw}(K)$ via restriction.
\end{lemma}

\begin{proof} Let $\nu$ be a quasi-faithful $w$-valuation of 
$Q$. By restriction, $\nu$ becomes a $w$-valuation of $K$
(we still denoted it by $\nu$). Let $\OO'$ be the support of 
$\gr_{\nu} Q$. By definition, $\OO/\OO'$ is torsion. Since 
$\OO={\mathbb Z}^n$, we have $\OO'\cong {\mathbb Z}^n$. 

Let $S$ be the support of $\gr_{\nu} K$. We claim that 
${\rm{rk}}_{\mathbb Z}S=n$, which implies that $\OO/S$ is 
torsion. Suppose to the contrary that ${\rm{rk}}_{\mathbb Z}S
<n$. Then ${\rm{rk}}_{\mathbb Z} \OO'/S>0$. This means
there is an element $v$ in $\OO'\setminus S$ such that
$d v\not\in S$ for all nonzero integer $d$. Since 
$\OO'$ is the support of $\gr_{\nu} Q$, there is an 
element $f\in Q$ such that $\nu(f)=v$. Since $Q$ is algebraic
over $K$, there are $a_0,\cdots, a_h\in K$ with $a_0=1$ and 
$a_h\neq 0$ such that $a_0+a_1 f+ a_2 f^2 + \cdots + a_h f^h=0$. 
When $a_i\neq 0$, then $\nu(a_i f^i)=\nu(a_i)+i \nu(f)
=\nu(a_i)+i v$. So $\{\nu(a_i f^i)\}_{i=0}^{h}$ are pair-wise
distinct (for those $a_i\neq 0$). Thus
$$\infty =\nu(a_0+a_1 f+ a_2 f^2 + \cdots + a_h f^h)
=\min_{a_i\neq 0}\{\nu(a_i f^i)\}\neq \infty,$$
which yields a contradiction. Therefore ${\rm{rk}}_{\mathbb Z}S
=n$ and consequently, $\OO/S$ is torsion. Therefore 
$(a')$ in Definition \ref{xxdef2.18} holds.

For (b) in Definition \ref{xxdef2.18}, we note that
$\GKdim \gr_{\nu} K\geq \GKdim \Bbbk S=n$. Then 
$$\GKdim \gr_{\nu}K \leq \GKdim K=n\leq \GKdim \gr_{\nu} K$$
which implies (b) in Definition \ref{xxdef2.18}. 

Since $\GKdim \gr_{\nu}K=n=\GKdim \gr_{\nu} Q$ and 
$\gr_{\nu}Q$ has nonzero $n$-Lie Poisson bracket,  
by Lemma \ref{xxlem1.10}, $\gr_{\nu} K$ has nonzero 
$n$-Lie Poisson bracket. Therefore $\nu$ is a 
quasi-faithful $w$-valuation of $K$.
\end{proof}

\section{Computations related to i.s. potentials}
\label{xxsec3}

In this section, we work out some valuations for 
Nambu-Poisson fields related to the Nambu-Poisson 
algebras given in Construction \ref{xxcon0.5}. We start
with an easy lemma.

\begin{lemma} 
\label{xxlem3.1}
Let $\nu$ be a valuation on an algebra $B$ and $Y\subseteq B$ 
be a finite-dimensional subspace of $B$. Then, there is a 
$\Bbbk$-linear basis $\{x_s\}_{s=1}^{m}$ of $Y$ such that 
$\nu(x_1)=\nu(x_s)$ for all $s=1,\cdots,m$.
\end{lemma}

\begin{proof}
Let $\{y_s\}_{s=1}^{m}$ be a $\Bbbk$-linear basis of $Y$. 
Let $y\in Y$ be a nonzero element such that
$\nu(y)=\min\{\nu(y_s)\mid s=1,\cdots, m\}$. By valuation 
axioms, $\nu(y)=\min \{\nu(f)\mid f\in Y\}$. For each $s$, let 
$x_s=y_s+\lambda y$ for $\lambda\in \Bbbk$. For a generic 
$\lambda$, $\nu(x_s)=\nu(y)$ for all $s$ and $\{x_1\}_{s=1}^{m}$
are linearly independent. The assertion follows.
\end{proof}

We recall some notations introduced in Construction 
\ref{xxcon0.5}. Let $\Bbbk^{[n+1]}$ be the polynomial ring 
$\Bbbk[t_1,\cdots,t_{n+1}]$ that has an internal grading with 
$\deg(t_i)=1$, called the {\it standard Adams grading}. Let 
$\Omega\in \Bbbk^{[n+1]}$ be a homogeneous potential with 
respect to the Adams grading. Define
\begin{equation}
\label{E3.1.1}\tag{E3.1.1}
d_0=\deg \Omega -(n+1).
\end{equation}
A special case is when $d_0=0$, or $\deg \Omega=n+1$. In this 
case both $A_{\Omega}$ and $P_{\Omega}$ are Nambu-Poisson 
$0$-graded algebras. The Nambu-Poisson structure on $A_{\Omega}$
(and on $P_{\Omega}$ and $P_{\Omega-\xi}$) is determined by
\begin{equation}
\label{E3.1.2}\tag{E3.1.2}
\{t_1,\cdots, \widehat{t_i},\cdots,t_{n+1}\}
=(-1)^{n+1-i}\frac{\partial \Omega}{\partial t_i}=(-1)^{n+1-i}\Omega_{t_i}
\end{equation}
for every $i$. The following lemma is clear.

\begin{lemma}
\label{xxlem3.2}
Let $B$ be $A_{\Omega}$, or $P_{\Omega}$, or $P_{\Omega-\xi}$
and let $Y=V={\rm Span}_\kk(t_1,\ldots,t_{n+1})$. With the notation introduced in Lemma \ref{xxlem3.1}
let $\{x_1,\cdots,x_{n+1}\}$ be a basis of $V$ satisfying the 
conclusion of Lemma \ref{xxlem3.1}.
\begin{enumerate}
\item[(1)]
The $\Bbbk$-span of $\{x_1,\cdots, \widehat{x_i}, \cdots, 
x_{n+1}\}$ for $1\leq i\leq n+1$ is equal to the $\Bbbk$-span 
of $\{t_1,\cdots, \widehat{t_i},\cdots,t_{n+1}\}$ for 
$1\leq i\leq n+1$.
\item[(2)]
$A_{sing}:
=A_{\Omega}/(\Omega_{t_i}, 1\leq i\leq n+1)$ is equal to
$A_{\Omega}/(\{x_1,\cdots,\widehat{x_i},\cdots,x_{n+1}\},
1\leq i\leq n+1)$. 
\end{enumerate}
\end{lemma}

We will use the following very nice result of 
Umirbaev-Zhelyabin \cite{UZ}.

\begin{lemma} 
\label{xxlem3.3}
Suppose $\Omega$ is an i.s. potential. Let $\xi\in \Bbbk^{\times}$.
\begin{enumerate}
\item[(1)]
\cite[Theorem 1 and Corollary 2]{UZ}
$P_{\Omega-\xi}$ is a simple Nambu-Poisson algebra. 

\item[(2)]
$\gldim P_{\Omega-\xi}=n$.
\item[(3)]
Let $P$ be $A_{\Omega}$ or $P_{\Omega}$. Let $\phi:
P\to B$ be a morphism of Nambu-Poisson domains.
Suppose that {\rm{(a)}} $\phi$ is surjective and {\rm{(b)}} 
$\GKdim B<n$. Then $t_i\in \ker \phi$ for all
$i$ and $B=\Bbbk$.
\end{enumerate}
\end{lemma}

\begin{proof}
(2)  This follows easily from the 
Jacobian criterion for hypersurfaces. 

(3) Since $\GKdim B<n$, by Lemma \ref{xxlem1.10}(1), 
the Nambu-Poisson structure on $B$ is zero. 
Since $\phi(\Omega_{t_i})
=(-1)^{n+1-i}\{\phi(t_1),\cdots,\widehat{\phi(t_i)},\cdots,\phi(t_{n+1})\}=0$,
$\Omega_{t_i}\in \ker \phi$ for all $i$. Since $B$ is a domain and a factor ring of a finite-dimensional algebra
$A_{sing}$, $t_i\in \ker \phi$ 
for all $i$. Thus $B=\Bbbk$.
\end{proof}

\begin{notation}
\label{xxnot3.4}
For the rest of this section, we assume
\begin{enumerate}
\item[(1)]
$\Omega$ is an i.s. potential of degree $d_0+n+1$ for
some $d_0\geq 0$.
\item[(2)]
Let $V$ denote the $\Bbbk$-vector space
$\sum_{s=1}^{n+1} \Bbbk t_s$, which may be considered as a 
subspace of $A_{\Omega}$, or $P_{\Omega}$, or $P_{\Omega-\xi}$.
When $\nu$ is a valuation, then let $\{x_s\}_{s=1}^{n+1}$ be 
as in Lemma \ref{xxlem3.1} when we take $Y=V$.
\item[(3)]
Let $\rho=\min\{\nu(t_s)\mid s=1,\cdots, n+1\}\in \OO$. Then 
$\rho=\nu(x_i)$ for each $i$ as in Lemma \ref{xxlem3.1}. 
\item[(4)] For any $\rho\in \OO$, we use $\mathbb F^{Id\rho}$ to denote the induced filtration $\mathbb F^{ind}$ on $P_\Omega$ such that $\deg(t_i)=\rho$ for all $i$. In particular, if $\rho<0$,  we use $\mathbb F^{c\rho}$ to denote the induced filtration $\mathbb F^{ind}$ on $P_{\Omega-\xi}$ for $\xi\in \kk^\times$ such that $\deg(t_i)=\rho$ for all $i$.
\end{enumerate}
\end{notation}

\begin{lemma}
\label{xxlem3.5}
Let $B$ be either $A_{\Omega}$, $P_{\Omega}$, or $P_{\Omega-\xi}$ 
for some $\xi\in \Bbbk^{\times}$. Let $\nu$ be a $w$-valuation of 
the Nambu-Poisson field $N:=Q(B)$. Let $d_0\geq 0$ and $\rho$ be
as in Notation \ref{xxnot3.4}.
\begin{enumerate}
\item[(1)]
Suppose there is no $v\in \OO$ such that $v> 0$ and $w\geq d_0 v$. 
Then $\rho\geq 0$. As a consequence, ${^w\Gamma_{0}}(N)\supseteq B$. 
\item[(2)]
Suppose that $B=P_{\Omega-\xi}$ with $\xi\in \Bbbk^{\times}$ and 
that $d_0=0$. Then there is no $w$-valuation for every $w<0$.
\item[(3)]
Suppose that $B=P_{\Omega-\xi}$ with $\xi\in \Bbbk^{\times}$ and 
that $d_0>0$. Then, there is no nontrivial $0$-valuation.
\item[(4)]
Let $B$ be $P_{\Omega}$ or $P_{\Omega-\xi}$.
If $w=d_0 v>0$ for some $v\in \OO$ and there is no
$v'$ strictly between $v$ and $0$, then either $\rho\geq 0$ or 
$\rho<0$. Moreover if $\rho<0$, then $\rho=-v$ and ${\mathbb F}={\mathbb F}^{Id\rho}$ if $B=P_{\Omega}$ and
${\mathbb F}={\mathbb F}^{c\rho}$ if $B=P_{\Omega-\xi}$. 
\item[(5)]
Let $B$ be $P_{\Omega}$ with $d_0>0$. Then there is no
faithful $0$-valuation on $N$.
\end{enumerate}
\end{lemma}

\begin{proof} (1) 
By Lemma \ref{xxlem3.1} and its proof, there are linearly 
independent elements $\{x_s\}_{s=1}^{n+1}\subset V$ such that 
$\nu(x_s)=\rho$ for all $s$. Let $\Omega_{i}=\{x_1, \cdots,
\widehat{x_i},\cdots, x_{n+1}\}$ for each $1\leq i\leq n+1$. 
Then the Adams degree of $\Omega_{i}$ is $d_0+n$. By the choice 
of $\rho$ and valuation axioms, we have $\nu(\Omega_{i})\geq 
(d_0+n)\rho$ for $1\leq i\leq n+1$. 

Let $I$ be the subalgebra of $\gr_{\nu} N$ generated by 
$\overline{x_s}$ for all $s$.

Case 1: Suppose $\nu(\Omega_{i})>(d_0+n)\rho$ for all $i$.
Then $\overline{\Omega_{i}}=0$ in $I$, or equivalently,
$\Omega_{i}(\overline{x_1}, \cdots,\overline{x_{n+1}})=0$
in $I$ for all $i$. This means that $I$ is a factor ring of $A_{sing}$
by Lemma \ref{xxlem3.2}(2). Since $I$ is a domain, it 
forces that $I=\Bbbk$. Consequently, $\rho=0$. 

Case 2: Suppose $\nu(\Omega_{i})=(d_0+n)\rho$ for some $i$. 
Then 
$$(d_0+n) \rho=\nu(\Omega_{i})=\nu(\{x_1,\cdots,\widehat{x_i},
\cdots,x_{n+1}\}) \geq \sum_{s\neq i}\nu(x_s)-w=n\rho- w$$
which implies that $d_0 \rho \geq -w$, or equivalently, $d_0(-\rho)
\leq w$. By the hypothesis in part (1), $-\rho\leq 0$, or
$\rho\geq 0$ as required.

Since $B$ is generated by $V$, the consequence follows from the 
definition.

(2) Suppose to the contrary that $\nu$ is a $w$-valuation on 
$N$ with $w<0$. Since $d_0=0$, by the proof of part (1), only 
Case 1 can happen. As a consequence, $I=\Bbbk$ and $\rho=0$. 
In this case $F_0^{\nu}(B)=B$ and $F_0^{\nu}(B)/F_{>0}^{\nu}(B)
=\Bbbk$, which contradicts Lemma \ref{xxlem3.3}(1).

(3) Let $\nu$ be a $0$-valuation on $K$. By part (1), ${^0\Gamma_{0}}(N)\supseteq B$. Consequently, $F^{\nu}_0(B)=B$. Since
$B$ is Nambu-Poisson simple [Lemma \ref{xxlem3.3}(1)], the 
Nambu-Poisson ideal $F^{\nu}_1(B)$ is zero. Thus, $\nu$ is trivial.

(4) Retain the notation introduced in the proof of part (1).

Case 1: Suppose $\nu(\Omega_i)>(d_0+n)\rho$ for all $i$.
The exact proof shows that $I=\Bbbk$ and that $\rho=0$. 

Case 2: Suppose $\nu(\Omega_i)=(d_0+n)\rho$ for some $i$. Then 
$$(d_0+n) \rho=\nu(\Omega_i)=
\nu(\{x_1,\cdots, \widehat{x_i},\cdots,x_{n+1}\})
\geq \sum_{s\neq i}\nu(x_s)-w=n \rho - w$$
which implies that $d_0 \rho \geq -w$. 
If $\rho\geq 0$, we are done. Otherwise, we have $\rho<0$
and $w\geq d_0 (-\rho)$. By the hypothesis on $w$, $v=-\rho$. 

It is easy to check that ${\mathbb F}^{ind}$ on $B$ by using 
$\deg(x_s)=\rho$ for all $s$ is indeed a filtration of $B$. 
We have further
${\mathbb F}^{ind}={\mathbb F}^{Id\rho}$ if $B=P_{\Omega}$ and
${\mathbb F}={\mathbb F}^{c\rho}$ if $B=P_{\Omega-\xi}$ according to Remark \ref{xxnot3.4}(4). 

By the choice of $i$, $\overline{\Omega_{i}}\neq 0$, and consequently,
$\{\overline{x_1},\cdots, \widehat{\overline{x_i}},\cdots,
\overline{x_{n+1}}\}=\overline{\Omega_{i}}\neq 0$. This shows that 
$\{\overline{x_s}\}_{s\neq i}$ are algebraically independent 
by Lemma \ref{xxlem1.10}(1).
So $\GKdim I\geq n$. By Lemma \ref{xxlem2.10}(1),
${\mathbb F}={\mathbb F}^{ind}$. 

(5) Suppose to the contrary that $\nu$ is a faithful 
$0$-valuation on $B$. By part (1), $\rho\geq 0$. 

Case 1: Suppose $\rho>0$. Let $\{x_s\}_{s=1}^{n+1}$ be as in the 
proof of part (1). By definition, $\deg \overline{x_i}=\nu(x_i)
=\rho>0$ for all $i$. Thus $I$ is not $\Bbbk$ and whence 
$\GKdim I\geq 1$. So, Case 1 of the proof of part (1) is impossible. 
This implies that only Case 2 can happen. Then 
$\{\overline{x_1},\cdots, \widehat{\overline{x_i}},
\cdots, \overline{x_{n+1}}\}=\overline{\Omega_i}\neq 0$
for some $i$. So $\{\overline{x_s}\}_{s\neq i}$ are
algebraically independent by Lemma \ref{xxlem1.10}(1). 
Then $\GKdim I=n$. By Lemma \ref{xxlem2.10}(1),
${\mathbb F}={\mathbb F}^{ind}$ since $B$ is graded to 
start with.

Since $\rho>0$, for every $i$,
$$\nu(\{x_1,\cdots,\widehat{x_i},\cdots, x_{n+1}\})
=\nu(\Omega_i)\geq (d_0+n) \rho=\sum_{s\neq i}\nu(x_s)-(-d_0\rho).$$
By Lemma \ref{xxlem2.9}(3), ${\mathbb F}^{ind}$ is a
classical $0$-valuation. So $\nu$ is not a faithful 
$0$-valuation, yielding a contradiction. 

Case 2: $\rho=0$. Then $F^\nu_0(B)=B$. Since $\nu$ is a faithful $0$-valuation, $B/F^\nu_{>0}B$ is a n-Lie Poisson domain of GK dimension $<n$. By Lemma \ref{xxlem1.10}(1), we get $\overline{\Omega_i}=\{\overline{x_1},\cdots \widehat{\overline{x_i}},\cdots,\overline{x_{n+1}}\}=0$ for all $i$. This implies that $B/F^\nu_{>0}B$ is a factor ring of $A_{sing}$ by Lemma \ref{xxlem3.2}(2) and hence $B/F^\nu_{>0}B=\kk$. Then all $i$, 
$h_i:=\overline{x_i}\in \Bbbk^{\times}$. 
Let $y_i=x_i -c_i x_1$ where $c_1=0$ and 
$c_i=\frac{h_i}{h_1}$ for all $i\geq 2$. Then 
$\nu(y_i)>0$ for all $i\geq 2$, and 
$$\begin{aligned}
\nu(\Omega_i)&=\nu(\{x_1,\cdots \widehat{x_i},\cdots,x_{n+1}\})\\
&=\nu(\{x_1, y_2+c_2 x_1, \cdots, \widehat{y_i+c_i x_1},\cdots, 
y_{n+1}+c_{n+1} x_{1}\})\\
&=\nu(\{x_1, y_2, \cdots, \widehat{y_i}, \cdots, y_{n+1}\})>0
\end{aligned}$$
for each $i=2,\cdots,n+1$. A similar argument can be made to show
that $\nu(\Omega_1)>0$. This means that each $\Omega_i$ does not 
contain the term $y_1^{d_0+n}$ when written as a polynomial
in $y_1,y_2,\cdots,y_{n+1}$, or $\Omega_i\in (y_2,\cdots,y_{n+1})$ 
-- the ideal generated by $y_2,\cdots, y_{n+1}$. So 
$\GKdim A/(\Omega_1,\cdots, \Omega_{n+1}) \geq 1$, which contradicts to the 
hypothesis of $\Omega$ being an i.s. potential by 
Lemma \ref{xxlem3.2}(2).
\end{proof}

\begin{lemma}
\label{xxlem3.6}
Let $B$ be $P_{\Omega}$. Suppose that $\nu$ is a $w$-valuation 
on $B$ and that ${\mathbb F}$ is the associated filtration of 
$B$. Let $I$ be the subalgebra of $\gr_{\nu} B$ generated by 
$\{{\overline{x_i}}\}_{s=1}^{n+1}$.
\begin{enumerate}
\item[(1)]
Suppose $d_0=0$ and $w=0$. If $F_0(B)=B$, then either 
${\mathbb F}$ is a trivial filtration or $F_0/F_1\cong \Bbbk$ 
and $\rho>0$. As a consequence, if $\rho=0$, then $\nu$ 
is trivial.
\item[(2)]
Suppose $\rho>0$. If $\GKdim I=n$, then ${\mathbb F}$ is 
the filtration ${\mathbb F}^{Id\rho}$.
\item[(3)]
Suppose $\rho<0$. If $\GKdim I=n$, then ${\mathbb F}$ is 
the filtration ${\mathbb F}^{Id\rho}$. As a consequence, 
$F_0(B)=\Bbbk$.
\end{enumerate}
\end{lemma}

\begin{proof}
(1) If ${\mathbb F}$ is nontrivial, then $F_{>0}(B)\neq 0$. 
Since ${\mathbb F}$ is a $0$-filtration, $F_{>0}(B)$ is a
Nambu-Poisson ideal of $B$. And $B/F_{>0}(B)$ is a Nambu-Poisson 
algebra of GK-dimension $<n$. Then the Nambu-Poisson bracket on
$B/F_{>0}(B)$ is trivial by Lemma \ref{xxlem1.10}(1). We apply 
\eqref{E3.1.2} to conclude that $\Omega_{t_i}\in F_{>0}(B)$ for 
all $1\le i\le n+1$. Because $\Omega$ has an isolated 
singularity, $B/F_{>0}(B)$ is finite-dimensional (and hence 
isomorphic to $\Bbbk$) and $F_{>0}(B)$ is generated by $t_i$ 
for $i=1,\cdots,n+1$. The main assertion follows.

Now we assume that $\rho=0$. By the valuation axiom, $F_0(B)=B$. 
The consequence follows from the main assertion.

(2) Let ${\mathbb F}^{ind}$ be the induced filtration with 
$\deg(x_s)=\rho>0$ for $s=1,\cdots,n+1$. Then 
${\mathbb F}^{ind}$ is equal to ${\mathbb F}^{Id\rho}$ and 
$\gr_{{\mathbb F}^{ind}} B\cong B$ since $B$ is graded to start 
with. The assertion follows from Lemma \ref{xxlem2.10}(1).

(3) Let ${\mathbb F}^{ind}$ be the induced filtration with 
$\deg(x_s)=\rho<0$ for $s=1,\cdots,n+1$. Then 
${\mathbb F}^{ind}$ is equal to ${\mathbb F}^{Id\rho}$ and 
$\gr_{{\mathbb F}^{ind}} B\cong B$ since $B$ is graded to start with. 
The assertion follows from Lemma \ref{xxlem2.10}(1). The 
consequence is clear.
\end{proof}

Now, it is ready to work out all faithful $0$-valuations on some 
Poisson fields. We say two filtrations are equivalent if their
associated valuations are equivalent.

\begin{theorem}
\label{xxthm3.7} 
Let $B=P_{\Omega}$ and let $N$ be the Nambu-Poisson fraction 
field $Q(B)$. Assume that $d_0=0$ and $\OO={\mathbb Z}$. 
\begin{enumerate}
\item[(1)]
$({\mathcal V}_{0}(N)/\sim) =\{\nu^{Id},\nu^{-Id}\}=
{\mathcal V}_{f0}(N)$. As a consequence, $\alpha_0(N)=2$.
\item[(2)]
Let $\nu$ be a nontrivial $0$-valuation on $N$. Then
$Q(\gr_{\nu}(N))\cong N$. As a consequence, every nontrivial 
$0$-valuation on $N$ is equivalent to a faithful $0$-valuation.
\item[(3)]
$\bfd(N)=0$ and $\bfw(N)=1$ and $N$ is quasi-Adams.
\end{enumerate}
\end{theorem}

\begin{proof}
(1) Let $\nu$ is a nontrivial $0$-valuation on $N$. 
By Lemma \ref{xxlem3.6}(1), $\rho\neq 0$. So 
we have two cases to consider.

Case 1: $\rho>0$. Let $I$ be the subalgebra of 
$\gr_\mathbb FP$ generated by ${\overline x}_s$ for 
$s=1,\cdots,n+1$. Since $B$ is Nambu-Poisson $0$-graded 
with Adams grading $|t_i|=1$ for all $i$, ${\mathbb F}^{ind}$ 
is equivalent to ${\mathbb F}^{Id}$. Then $\phi^{ind}: B\to 
I\subseteq \gr_{\mathbb F} B$ is a Nambu-Poisson algebra 
homomorphism [Lemma \ref{xxlem2.11}(2)]. Since $\GKdim I\geq 1$, 
by Lemma \ref{xxlem3.3}(3), $\GKdim I\geq n$. By Lemma 
\ref{xxlem3.6}(2) and its proof, ${\mathbb F}^{ind}={\mathbb F}$. 
Consequently, ${\mathbb F}$ is equivalent to ${\mathbb F}^{Id}$ 
or $\nu$ is equivalent to $\nu^{Id}$.

Case 2: $\rho<0$. Let $I$ be the subalgebra of $\gr_{\mathbb F}B$ 
generated by ${\overline x}_i$ for $i=1,\cdots,n+1$. Since $B$ 
is Nambu-Poisson $0$-graded with Adams grading $|t_i|=1$ for 
all $i$, ${\mathbb F}^{ind}$ is equivalent to ${\mathbb F}^{-Id}$. 
Then $\phi^{ind}: B\to I\subseteq \gr_{\mathbb F} B$ is a 
Nambu-Poisson algebra homomorphism [Lemma \ref{xxlem2.11}(2)]. 
Since $\GKdim I\geq 1$, by Lemma \ref{xxlem3.3}(3), $\GKdim I
\geq n$. By Lemma \ref{xxlem3.6}(3) and its proof, 
${\mathbb F}^{ind}={\mathbb F}$. Consequently, ${\mathbb F}$ is 
equivalent to ${\mathbb F}^{-Id}$ or $\nu$ is equivalent to 
$\nu^{-Id}$.

Therefore, up to equivalence $\sim$, there are only two nontrivial 
$0$-valuations, namely $\nu^{Id}$ and $\nu^{-Id}$. It is clear that they are faithful. The assertions 
follow.

(2) Clear from the proof of part (1) and the fact $\gr_{\nu^{\pm Id}}
B\cong B$, see discussion before Lemma \ref{xxlem2.12}.

(3) These are consequences of parts (1) and (2).
\end{proof}

Next we study the Nambu-Poisson algebra $P_{\Omega-\xi}$ 
for $\xi\neq 0$ when $d_0=0$.

\begin{lemma}
\label{xxlem3.8}
Suppose that $\OO={\mathbb Z}$ and that $w=d_0\geq 0$. Then 
$B:=P_{\Omega-\xi}$ has at least one faithful $w$-valuation 
$\nu^{c}$ as given in Lemma {\rm{\ref{xxlem2.12}(3)}}. 
\end{lemma}

\begin{proof}
Note that $B$ has a filtration determined by \eqref{E2.8.3} and 
the degree assignment $\deg(t_i)=-1=-|t_i|$ for $i=1,\cdots,n+1$. 
Or equivalently, we define a filtration 
${\mathbb F}^{c}:=\{F_i^{c}\}$ on $P$ by 
$$F_{-i}^{c}(B)=\{\sum_j a_j f_j\in B\mid a_j\in \Bbbk, f_j 
{\text{ are monomials of Adams degree $\leq i$}}\}$$ 
for all $i$. Then, the associated graded ring is 
$\gr_{{\mathbb F}^{c}} B \cong P_{\Omega}$. Since 
$|\Omega|=n+1+d_0$, \eqref{E3.1.2} imply that the 
filtration is a $w$-filtration by Lemma \ref{xxlem2.9}(3). 
The assertion follows.
\end{proof}

Now, we are ready to work out another example.

\begin{theorem}
\label{xxthm3.9}
Let $B$ be $P_{\Omega-1}$ and $N$ be the Nambu-Poisson fraction 
field $Q(B)$. Assume that $\OO={\mathbb Z}$ and that $d_0=0$. 
\begin{enumerate}
\item[(1)]
There is a unique nontrivial $0$-valuation up to equivalence,
or ${\mathcal V}_{0}(N)/\sim$ is a singleton. As a consequence,
$\alpha_{0}(N)=1$.
\item[(2)]
Let $\nu$ be a nontrivial $0$-valuation on $N$. Then, $\nu$ is
determined by the grading of $B$ as given in Lemma 
\ref{xxlem3.8}. As a consequence,
$Q(\gr_{\nu}(N))\cong Q(P_{\Omega})$ or $N\to_{\nu}
Q(P_{\Omega})$.
\item[(3)]
$\bfd(N)=1$ and $\bfw(N)=1$.
\end{enumerate}
\end{theorem}

\begin{proof}
(1,2) By Lemma \ref{xxlem3.8} and its proof, there is a faithful 
$0$-valuation, denoted by $\nu^{c}$, such that $N\to_{\nu^{c}} 
Q(P_{\Omega})$. It remains to show every nontrivial $0$-valuation 
is equivalent to $\nu^{c}$ that is given in Lemma \ref{xxlem3.8}. 
The rest of the proof is similar to the proof of Theorem 
\ref{xxthm3.7}.

Let $\nu$ be a nontrivial $0$-valuation on $N$, and we are using
the notations introduced in Notation \ref{xxnot3.4}. Let 
$\mathbb F$ be the $0$-filtration on $B$ associated to $\nu$. 
We have three cases to consider.

Case 1: $\rho=0$. Then $F_0(B)=B$. Since $B$ is Nambu-Poisson simple 
[Lemma \ref{xxlem3.3}(1)], $F_1(B)=0$. Thus, $\nu$ is trivial, 
yielding a contradiction.

Case 2: $\rho>0$. This contradicts the fact that $\Omega-1=0$ in $B$. 

Case 3: $\rho<0$. Let $I$ be the subalgebra $\gr_{\nu}(N)$
generated by $\{\overline{x}_s\}_{s=1}^{n+1}$. Let 
${\mathbb F}^{ind}$ be the induced filtration determined
by $\deg(x_s)=\rho<0$ for all $s=1,\cdots, n+1$. 
Up to an equivalence, we may assume that $\rho=-1$. Then 
${\mathbb F}^{ind}$ is a $0$-filtration with 
$\gr_{{\mathbb F}^{ind}} P\cong P_{\Omega}$.
By Lemma \ref{xxlem2.11}(1), 
$$\phi^{ind}: \gr_{{\mathbb F}^{ind}} B\to I\subseteq 
\gr_{\mathbb F} B$$
is a Nambu-Poisson algebra morphism. Note that 
$\gr_{{\mathbb F}^{ind}} B$ is $P_{\Omega}$. 
Since $\GKdim I\geq 1$, by Lemma \ref{xxlem3.3}(3), 
$\GKdim I\geq n$. By Lemma \ref{xxlem2.10}(1),
${\mathbb F}={\mathbb F}^{ind}$. So, we obtain a unique valuation 
in this case, which is given in Lemma \ref{xxlem3.8}. 
The assertion is proved.

(3) Let $\nu$ be the unique Poisson valuation on $N$ given in 
part (2). By parts (1,2), $Q(P_{\Omega-1})\to_{\nu}Q(P_{\Omega})$ 
is the only possible arrow. Hence $\bfw(N)=1$. By Theorem 
\ref{xxthm3.7}, there is no arrow $Q(P_{\Omega})\to_{\nu} K$ 
with $K\not\cong Q(P_{\Omega})$. Thus $\bfd(P)=1$. 
\end{proof}

\begin{theorem}
\label{xxthm3.10}
Suppose that $\OO={\mathbb Z}$ and that $d_0\geq 2$. Let $B$ 
be either $A_{\Omega}$, or $P_{\Omega}$, or $P_{\Omega-\xi}$ 
and let $N=Q(B)$. Then, for every $w$ between $1$ and $d_0-1$, 
${^w\Gamma_{0}}(N)=B$.
\end{theorem}

\begin{proof} By Lemma \ref{xxlem3.5}(1), 
${^w\Gamma_{0}}(N)\supseteq B$. By Theorem \ref{xxthm2.16}(2),
${^1\Gamma_{0}}(N)\subseteq B$. Hence
$$B\supseteq\, {^1\Gamma_{0}}(N)\supseteq\, ^2\Gamma_0(N)\supseteq\, \cdots 
\supseteq\, ^{d_0-1}\Gamma_{0}(N) \supseteq B.$$
The assertion follows.
\end{proof}

\begin{theorem}
\label{xxthm3.11}
Suppose that $\OO={\mathbb Z}$ and that $d_0\geq 1$. Let 
$B$ be $P_{\Omega-\xi}$ where $\xi\neq 0$ and let $N=Q(B)$.
Then 
$${^w\Gamma_{0}}(N)=\begin{cases} N & w\leq 0,\\
B & 1\leq w\leq d_0-1,\\
\Bbbk & w\geq d_0.
\end{cases}$$
\end{theorem}

\begin{proof}
It follows from Lemma \ref{xxlem3.5}(3) that 
${^w\Gamma_{0}}(N)=N$ if $w\leq 0$. If $1\leq w\leq d_0-1$ 
(if $d_0=1$, then there is no such $w$), the assertion 
follows from Theorem \ref{xxthm3.10}. Now let $w=d_0$. 
By Theorem \ref{xxthm2.16}(2), ${^w\Gamma_{0}}(N)\subseteq B$. 
By Lemma \ref{xxlem2.12}(3), there is an $d_0$-valuation 
$\nu:=\nu^{c}$ such that $\nu(t_i)=-1$ for all $i$ 
(see Lemma \ref{xxlem3.5}(4)). Consequently, $F^{\nu}_0(P)
=\Bbbk$. Then 
$${^{d_0}\Gamma_{0}}(N)\subseteq {^{d_0}\Gamma_{0}}(N)\cap 
F^{\nu}_0(N)\subseteq B\cap F^{\nu}_0(N)=F^{\nu}_0(B)=\Bbbk.$$
By Lemma \ref{xxlem2.15}(2), ${^w\Gamma_{0}}(N)=\Bbbk$ 
for all $w\ge d_0$. 
\end{proof}

\section{Secondary invariants}
\label{xxsec4}

In this section, we recall or introduce several secondary 
invariants related to the $w$-valuations. These invariants 
are helpful in understanding Nambu-Poisson fields from 
different aspects. 

\subsection{$\alpha$-type invariants}
\label{xxsec4.1}
Any invariant defined directly by using $w$-valuations is called
{\it of $\alpha$-type}. For example, $\alpha_w(N)$ in Remark 
\ref{xxrem2.17}(1) and $w(N)$ in Remark \ref{xxrem2.17}(3) are 
of $\alpha$-type. For the rest of this section, we discuss 
secondary invariants of other types.

\subsection{$\beta$-type invariants}
\label{xxsec4.2}

\begin{definition}
\label{xxdef4.1}
Let $K$ and $Q$ be two Nambu-Poisson fields and $w$ an element 
in $\OO$. We say $K$ {\it $w$-controls} $Q$ if there is a 
$w$-valuation $\nu$ on $K$ such that 
$$Q(\gr_{\nu}(K))\cong Q.$$ 
We write $K \to_{\nu} Q$ in this case.
\end{definition}

Any invariants that are involved with arrows $\to_{\nu}$ 
are called {\it of $\beta$-type}. 

\begin{definition}
\label{xxdef4.2}
Let $N$ be a Nambu-Poisson field and $w$ an element in $\OO$.
\begin{enumerate}
\item[(1)]
The {\it $w$-depth} of $N$ is defined to be
$$\bfd_w(N):=\sup\{n \mid 
\{N:=K_0\to_{\nu_1} K_1\to_{\nu_2} K_2\to_{\nu_3} 
\cdots \to_{\nu_{n}}K_n\} \}$$
where $K_i\not\cong K_j$, for all $0\leq i\neq j\leq n$ and 
each $\nu_i$ is a faithful $w$-valuation. 
\item[(2)]
The {\it $w$-width} of $K$ is defined to be
$$\bfw_w(N):=\# \left[
\{ Q\mid N\to_{\nu} Q {\text{ where $\nu$ is a 
faithful $w$-valuation}}\}/\cong \right]$$
where $\{\quad\} /\cong$ means the set of isomorphism classes.
\end{enumerate}
\end{definition}

When $w=0$, we will omit the subscript $w$ in the above 
definition as in previous sections. For example 
$\bfd_0(N)=\bfd(N)$ (resp. $\bfw_0(N)=\bfw(N)$) as in 
Definition \ref{xxdef4.2}(1) (resp. Definition 
\ref{xxdef4.2}(2)). We have already provided some results
concerning depth and width. We will prove one more result
about $\beta$-invariants.

Let $m\geq n$ and $\bfq:=\{q_{i_1,\cdots,i_n}
\mid 1\leq i_1< \cdots <i_n\leq m\}$ be a subset of $\Bbbk$. 
Let $T(\bfq)$ denote Nambu-Poisson torus 
$\Bbbk[x_1^{\pm 1},\cdots,x_m^{\pm 1}]$ with Nambu-Poisson 
bracket determined by 
$\{x_{i_1},\cdots, x_{i_n}\}=q_{i_1 \cdots i_n} x_{i_1}
\cdots x_{i_n}$ for $q_{i_1\cdots i_n}\in \bfq$ 
for all $1\leq i_1<\cdots< i_n\leq m$. 

\begin{theorem}
\label{xxthm4.3}
Let $\OO={\mathbb Z}$. Let $T$ be the Nambu-Poisson torus 
$T(\bfq)$ defined as above and $N$ be the Nambu-Poisson 
fraction field of $T$. Assume that $T$ is Nambu-Poisson simple.
\begin{enumerate}
\item[(1)]
There is a one-to-one correspondence between the set 
${\mathcal V}_{f0}(N)$ and the set 
$\{(v_1,\cdots,v_m)\in {\mathbb Z}^m \mid \gcd(v_1,\cdots,v_m)=1\}$. 
Further, for every $0$-valuation $\nu$, 
$Q(\gr_{\nu}(N))\cong N$. As a consequence, $\bfd_0(N)=0$, 
$\bfw_0(N)=1$, and $N$ is quasi-Adams.
\item[(2)]
There is no $w$-valuation for all $w<0$.
\end{enumerate}
\end{theorem}

\begin{proof}
(1) First, we prove that there is a one-to-one correspondence 
between the set of $0$-valuations and the set ${\mathbb Z}^m$. 
Let $\nu$ be a $0$-valuation on $N$. Let $v_i=\nu(x_i)$ for 
$i=1,\cdots,m$. Then $(v_1,\cdots,v_m) \in {\mathbb Z}^m$. 
Conversely, given $(v_1,\cdots,v_m)\in {\mathbb Z}^m$ we 
claim that there is a unique $0$-valuation $\nu$ such that 
$\nu(x_i)=v_i$ for all $i$. Let ${\mathbb F}^{ind}$ be 
the induced filtration determined by \eqref{E2.8.3} with the 
degree assignment $\deg(x_i)=v_i$ and $\deg(x_i^{-1})=-v_i$ 
for all $i=1,\cdots,m$. By the definition of $T$, it 
is easy to see that $T$ is in fact ${\mathbb Z}^m$-graded 
with degree assignment $\deg(x_i^{\pm 1})=\pm v_i$ for all $i$. 
By Lemma \ref{xxlem2.12}(1), ${\mathbb F}^{ind}$ agrees with 
the Adams$^{Id}$ filtration ${\mathbb F}^{Id}$. As a 
consequence, ${\mathbb F}^{ind}$ is a good filtration which 
provides a $0$-valuation $\nu$ such that $\nu(x_i)=v_i$ for 
all $i$. It remains to show that such a $\nu$ is unique. 

Now let $\nu$ be any $0$-valuation on $N$ such that 
$\nu(x_i)=v_i$ for all $i$ and let ${\mathbb F}$ be the 
associated filtration. By the valuation axioms [Definition
\ref{xxdef0.1}], $\nu(x_i^{-1})=-v_i$ for all $i$. By the 
argument in the previous paragraph and Lemma \ref{xxlem2.9}(2), 
we have a sequence of algebra morphisms
$$T\cong \gr_{{\mathbb F}^{Id}} T=
\gr_{{\mathbb F}^{ind}} T\to I\subseteq \gr_{\mathbb F} T.$$
It is clear that $T$ is Nambu-Poisson $0$-graded. By Lemma 
\ref{xxlem2.11}(1), 
$\phi^{ind}: \gr_{{\mathbb F}^{ind}} T\to \gr_{\mathbb F} T$ 
is a Nambu-Poisson $0$-graded algebra morphism. Since $T$ 
is Nambu-Poisson simple, $\phi^{ind}$ is injective. Thus 
$\GKdim I\geq \GKdim T=m$. By Lemma \ref{xxlem2.10}(1),
${\mathbb F}$ agrees with ${\mathbb F}^{ind}(={\mathbb F}^{Id})$. 
Therefore $\nu$ agrees with $\nu^{Id}$. This proves the 
uniqueness. 

Note that ${\nu}$ is a faithful $0$-valuation if and only if
$\gcd\{\nu(x_i),\cdots,\nu(x_m)\}=1$. The main assertion follows.
By the above proof, $\gr_{\nu} T\cong T$ for any Poisson
$0$-valuation $\nu$. By Lemma \ref{xxlem2.7}(2),
$Q(\gr_{\nu} N)\cong N$. Or equivalently, $N\to_{\nu} N$
for all $0$-valuations $\nu$. By definition, we have 
$\bfd_0(N)=0$, $\bfw_0(N)=1$, and $N$ is quasi-Adams.

(2) By the proof of part (1), for every $0$-valuation 
$\nu$, we have $Q(\gr_{\nu}(N))\cong N$. So $\nu$ is 
nonclassical. The assertion follows from Lemma 
\ref{xxlem1.8}(2).
\end{proof}

\subsection{$\gamma$-type invariants}
\label{xxsec4.3}

Recall that any invariant defined by using filtration
${\mathbb F}^{\nu}$ is called of $\gamma$-type. For example,
${^w\Gamma_{v}}(N)$ defined in Definition \ref{xxdef2.14} 
are $\gamma$-type invariants. In this subsection, we will use 
some $\gamma$-type invariants to distinguish some Nambu-Poisson 
fields. For part of this section, we assume

\begin{hypothesis}
\label{xxhyp4.4}
Let $\OO={\mathbb Z}^n$ for some $n\geq 2$.
\end{hypothesis} 

We also assume that $N$ is an $n$-Lie Poisson field of 
GK-dimension $n$ (=transcendence degree $n$). The same idea is expected to apply to $n$-Lie Poisson fields of 
higher GK-dimension. The following lemma is easy, and its 
proof is omitted.

\begin{lemma}
\label{xxlem4.5}
Let $A$ is a ${\mathbb Z}^n$-graded field such that 
$\dim_{\Bbbk} A_w=1$ for all $w\in {\mathbb Z}^n$. Then
$A\cong \Bbbk[x_1^{\pm 1},\cdots, x_{n}^{\pm 1}]$ where 
$\deg(x_i)=(0,\cdots, 0,1,0,\cdots,0)\in {\mathbb Z}^n$ 
where $1$ is in the $i$th position. 
\end{lemma}

Let us consider a class of Nambu-Poisson algebras. Let 
$q\in \Bbbk^{\times}$ and 
$\kappa:=(\kappa_1,\cdots,\kappa_n)\in \OO$. Let 
$T(q, \kappa_1, \cdots,\kappa_n)$ or $T(q, \kappa)$ be the 
Nambu-Poisson torus $\Bbbk[x_1^{\pm 1},\cdots, x_n^{\pm 1}]$ 
with the Nambu-Poisson structure determined by 
\begin{equation}
\label{E4.5.1}\tag{E4.5.1}
\{x_1,\cdots,x_{n}\} x_1^{-1} \cdots x_{n}^{-1}=q x_1^{\kappa_1}
\cdots x_n^{\kappa_n}.
\end{equation}
If each $\kappa_i\geq -1$, then $T(q, \kappa)$ is a Nambu-Poisson 
subalgebra of the localization $(A_{\Omega})_{S^{-1}}$ where 
$A_{\Omega}$ is given as in Construction \ref{xxcon0.5} with 
$S=\{x_i\}_{i=1}^n$ and where
\begin{equation}
\label{E4.5.2}\tag{E4.5.2}
\Omega:= q x_1^{1+\kappa_1}
\cdots x_n^{1+\kappa_n} x_{n+1}.
\end{equation}

Note that $T(q, \kappa)$ plays an important role in the study of 
$\OO$-valuations due to the following lemma and Definition 
\ref{xxdef4.11}.

\begin{lemma}
\label{xxlem4.6}
Suppose $\Bbbk$ is algebraically closed and assume Hypothesis 
\ref{xxhyp4.4}. Let $P$ be an $\OO$-graded field such that 
$\OO/\OO'$ is torsion where $\OO'$ is the support of $P$. 
Suppose that $P$ is Nambu-Poisson $\kappa$-graded for some 
$\kappa\in \OO$ of GK-dimension $n$. Then $P\cong T(q,\kappa)$.
\end{lemma}

\begin{proof} Under the hypothesis, $\OO'$ is an abelian group 
of rank equal to the rank of $\OO$. So $\GKdim \Bbbk \OO'=
\GKdim \Bbbk \OO=n$. By Lemma \ref{xxlem1.7}(2), $\dim_{\Bbbk}
P_i$ is either 1 or 0. 

Since $\OO'$ is free over ${\mathbb Z}$ of rank $n$, $\OO'
\cong {\mathbb Z}^n$. By Lemma \ref{xxlem4.5},
$P\cong \Bbbk[x_1^{\pm 1}, \cdots, x_n^{\pm 1}]$ as an 
$\OO'$-graded (whence $\OO$-graded) algebra. Since $P$ is 
Poisson $\kappa$-graded,
$$\{x_1,\cdots,x_n\}=qx_1^{1+\kappa_1}\cdots x_n^{1+\kappa_n}$$
for some $q\in \Bbbk$ and $(\kappa_1,\cdots,\kappa_n)=\kappa$.
This equation is equivalent to \eqref{E4.5.1}.
Since $P$ is generated by $x_1^{\pm 1},\cdots, 
x_n^{\pm 1}$. The Nambu-Poisson structure of
$P$ is completely determined by \eqref{E4.5.1}.
Therefore, the assertion follows.
\end{proof}

Let $N_{q, \kappa}$ be the Nambu-Poisson fraction field of 
$T(q, \kappa)$. By definition, $N_q\cong N_{q, 0,\cdots,0}$ 
and $N_{Weyl}\cong N_{1, (-1,\cdots,-1)}$. If $\kappa\in \OO$, 
$\gcd(\kappa)$ stands for $\gcd(\kappa_1,\cdots,\kappa_n)$. The 
following lemma is a generalization of a result in \cite{GZ}. 

\begin{lemma}
\label{xxlem4.7}
Retain the above definition. Assume that $q, q'\neq 0$ and 
that $\kappa\neq 0$.
\begin{enumerate}
\item[(1)]
$T(q, \kappa)\cong T(1, \kappa)$.
\item[(2)]
$T(q, \kappa)$ is isomorphic to 
$T(q, \kappa_0,0,\cdots,0)$ where $\kappa_0:
=\gcd(\kappa)$.
In particular, $T(q, \kappa_0,0,\cdots,0)\cong 
T(q', \kappa_0,d_2\kappa_0,\cdots,d_n \kappa_0)$
for all $d_i\in {\mathbb Z}$. 
\item[(3)] 
$T(q, \kappa)$ is isomorphic to 
$T(q', \kappa')$ if and only if 
$\gcd(\kappa')=\gcd(\kappa)$.
\item[(4)]
$T(q, \kappa_0,0,\cdots,0)$ is isomorphic to 
$T(q', \kappa'_0,0,\cdots,0)$
if and only if $\kappa_0=\pm \kappa'_0$.
\item[(5)]
$T(q, \kappa)$ is Nambu-Poisson simple.
\item[(6)]
There is a Poisson algebra embedding $T(1, \kappa)\to 
T(1, d\kappa)$ for all $d\neq 0$.
\end{enumerate}
\end{lemma}

\begin{proof}
(1) Without loss of generality, we may assume $\kappa_1\neq 0$.
The assertion follows from the isomorphism determined by
$\phi: x_1\to c x_1, x_i\to x_i$ for all $i\geq 2$ where 
$qc^{\kappa_1}=1$.

(2) Let $a_1:=(a_{11},\cdots,a_{1n})$ be the vector
$\kappa/\kappa_0$ in $\OO$. Then $\gcd(a_1)=1$. By 
induction on $n$, one can find vectors $a_2,\cdots,a_n$ in 
$\OO$ such that $\sum_{i=1}^n a_i {\mathbb Z}=\OO$. Let 
$\phi$ be the isomorphism 
$T(q',\kappa_0,0,\cdots,0)\to T(q, \kappa)$ (for some $q'\neq 0$) 
sending 
$$x_i\to x_1^{a_{i1}}\cdots x_{n}^{a_{in}}, \quad 
{\text{for all $i=1,\cdots,n$.}}$$
One can easily see that $\phi$ induces a Nambu-Poisson algebra 
isomorphism. We may assume $q'=q$ by part (1).

(3) Suppose we have an isomorphism $\phi: T:=T(q,\kappa)
\to T(q',\kappa')=:T'$. Since $\OO:={\mathbb Z}^n$ is an 
ordered group, only invertible elements in the group algebra 
$\Bbbk \OO$ are group elements in $\OO$ up to a nonzero scalar. 
This means that only invertible elements in $T$ are 
$x_1^{m_1}\cdots x_n^{m_n}$ for some $m_1,\cdots,m_n\in 
{\mathbb Z}$. Then, up to scalars, $\phi$ is
determined by a matrix in $\begin{pmatrix} a_{11}& \cdots & a_{1n}
\\ \cdots & \cdots &\cdots \\
a_{n1}& \cdots & a_{nn}\end{pmatrix}
\in GL_n({\mathbb Z})$, namely,
$$\phi: x_1\to x^{a_{11}} \cdots x_{n}^{a_{1n}}, 
\cdots, x_{n}\to x_1^{a_{n1}} \cdots x_{n}^{a_{nn}}.$$
Applying $\phi$ to \eqref{E4.5.1}, we obtain that each 
$\kappa'_i$ (for $i=1,\cdots,n$) is an integral combination 
of $\kappa_1,\cdots,\kappa_n$. Thus $\gcd(\kappa)$ divides 
$\gcd(\kappa')$. Since $\phi$ has an inverse, 
$\gcd(\kappa')$ divides $\gcd(\kappa)$. Therefore
$\gcd(\kappa)=\gcd(\kappa')$. This is one 
direction. The other direction follows by using part (2) 
twice and part (1).

(4) This is a special case of part (3).

(5) If $T$ is not Nambu-Poisson simple, then $T$ has a 
nonzero Nambu-Poisson prime ideal $I$ by Lemma \ref{xxlem1.11}. 
Then $T/I$ has GK-dimension at most $n-1$. By Lemma 
\ref{xxlem1.10}(1), the Nambu-Poisson structure on
$T/I$ is zero. By \eqref{E4.5.1}, $x_1^{\kappa_1}\cdots x_{n}^{\kappa_n}
=0$. But, by definition, $x_1^{\kappa_1}\cdots x_{n}^{\kappa_n}$ is 
an invertible element, a contradiction. Therefore 
$T$ is Poisson simple.

(6) 
Let $x_1,\cdots,x_n$ be the set of generators for 
$T(1,n\kappa_0,0,\cdots,0)$.
Let $u=x_1^n$. Then, by \eqref{E1.3.1}, we have 
$$\begin{aligned}
\{u,x_2,\cdots,x_n\} u^{-1} x_2^{-1}\cdots x_{n}^{-1}
&=n\{x_1,x_2,\cdots,x_n\} x_1^{-1}\cdots x_{n}^{-1}\\
&=n x_1^{n \kappa_0} =n u^{\kappa_0},
\end{aligned}
$$
which implies that $u$ and $x_2,\cdots,x_n$ generate 
$T(n,\kappa_0,0,\cdots,0)$. The assertion follows.
\end{proof}

By Lemma \ref{xxlem4.7}(1,2), if $\kappa\neq 0$, then 
$N_{q,\kappa}\cong N_{1, (\gcd(\kappa),0,\cdots,0)}$. 
For every $k\geq 1$, we let $N(k)$ be 
$N_{1, (k,0,\cdots,0)}$. It is clear that $N_{Weyl}
\cong N(1)$ and that every $N_{q,\kappa}$ with $\kappa
\neq 0$ is isomorphic to $N(k)$ for $k=\gcd(\kappa)$. 
Next, we show that $N(k)$ are non-isomorphic for distinct
$k$ (also a generalization of a result of \cite{GZ} in 
the Poisson setting).

\begin{lemma}
\label{xxlem4.8}
Let $k\geq 1$ and $N:=N(k)$ be defined as above.
Let $\OO={\mathbb Z}$.  
\begin{enumerate}
\item[(1)]
${^w\Gamma_{0}}(N)=\Bbbk[x_1]$ for $1\leq w\leq k-1$.
\item[(2)]
${^w\Gamma_{0}}(N)=\Bbbk$ for all $w\geq k$. 
\end{enumerate}
\end{lemma}

\begin{proof}
(1) By definition, $A=\Bbbk[x_1,\cdots,x_n]$ is a Nambu-Poisson 
subalgebra of $N$ with 
\begin{equation}
\label{E4.8.1}\tag{E4.8.1}
\{x_1,\cdots,x_n\}=(x_1\cdots x_n) x_1^k
\end{equation}
such that $Q(A)=N$. By Theorem \ref{xxthm2.16}(2), 
${^1\Gamma_{0}}(N)\subseteq A$. Note that $B:=\Bbbk[x_1,x_{2}^{-1},\cdots,
x_n^{-1}]$ is also a Nambu-Poisson subalgebra of $N$ (with 
\begin{equation}
\notag
\{x_1,x_2^{-1}\cdots,x_n^{-1}\}=(-1)^{n-1}(x_1x_2^{-1}\cdots x_n^{-1}) x_1^k
\end{equation}
such that $Q(B)=N$. By Theorem \ref{xxthm2.16}(2), 
${^1\Gamma_{0}}(N)\subseteq B$. Then 
${^1\Gamma_{0}}(N)\subseteq A\cap B=\Bbbk[x_1]$. 
As a consequence, ${^w\Gamma_{0}}(N)\subseteq \Bbbk[x_1]$
for all $w\geq 1$. 

Let $\nu$ be a $w$-valuation on $N$. Then 
$$k\nu(x_1)=\nu(x_1^k)=\nu(\{x_1,\cdots,x_n\}x_1^{-1}\cdots x_n^{-1})
\geq 0-w.$$
Hence $\nu(x_1)\geq -\frac{w}{k}$.
If $1\leq w<k$, $\nu(x_1)\geq 0$ and $x_1\in {^w\Gamma_{0}}(N)$.
So ${^w\Gamma_{0}}(K)=\Bbbk[x_1]$. 

(2) If $w\geq k$, there is a
$w$-valuation defined by $\nu(x_1)=-1$ and $\nu(x_s)=0$
for all $s=2,\cdots,n$. So $x_1\not\in {^w\Gamma_{0}}(N)$. 
In fact, one sees that any polynomial $f(x_1)$ of positive 
degree is not in ${^w\Gamma_{0}}(N)$. So ${^w\Gamma_{0}}(N)=
\Bbbk$. 
\end{proof}

As an immediate consequence, we have

\begin{corollary}
\label{xxcor4.9}
The Nambu-Poisson fields $\{N(k)\}$ are non-isomorphic
for distinct $k\geq 1$. Moreover, there is no embedding from $N(k)$ to $N(k')$ for $1\le k'<k$.
\end{corollary}

\subsection{$\delta$-type invariants}
\label{xxsec4.4}

In this subsection, we will introduce some $\delta$-type 
invariants which are defined by the associated graded 
rings $\gr_{\nu}(N)$. For a large of this subsection we 
assume Hypothesis \ref{xxhyp4.4} where $N$ is an $n$-Lie 
Poisson field. We will set up some foundation for all 
$\delta$-type invariants.  

\begin{definition}
\label{xxdef4.10}
The {\it $\kappa$-invariant} of $T(q,\kappa)$
for $q\neq 0$ is defined to be
$$\kappa(T(q,\kappa):=\gcd(\kappa)$$ 
which is a nonnegative integer.
\end{definition}

By Lemma \ref{xxlem4.7}(3), $\kappa(T(q,\kappa))$ is a
Nambu-Poisson algebra invariant of $T(q,\kappa)$. Now we 
define some $\delta$-type invariants of Nambu-Poisson fields 
related to $\gr_{\nu} (N)$. 

\begin{definition}
\label{xxdef4.11} 
Suppose that $\OO={\mathbb Z}^n$. Let $N$ be an 
$n$-Lie Poisson field of GK-dimension $n$ and let 
$w\in \OO$. 
\begin{enumerate}
\item[(1)]
The {\it $\kappa_w$-invariant} of $N$ is defined to be
$$\kappa_{w}(N):=\{ \kappa(\gr_{\nu} N) \mid \gr_{\nu} N
\cong T(q,\kappa)\}$$
where $\nu$ runs over all $w \in{\mathcal V}_{qfw}(N)$.
If ${\mathcal V}_{qfw}=\emptyset$, we write 
$\kappa_w(N)=\emptyset$. In general, $\kappa_w(N)$ is a 
collection of nonnegative integers.
\item[(2)]
Suppose $0\in \kappa_0(N)$. The {\it $\varrho$-invariant} 
of $N$ is defined to be
$$\varrho(N):=\{q \mid 
\nu\in {\mathcal V}_{qf0}(N),\; 
\gr_{\nu} (N)\cong T(q,0)\}.$$
If there is no such $\nu$, we write $\varrho(N)=\emptyset$.
Generally, $\varrho(N)$ is a collection of nonzero 
scalars in $\Bbbk$.
\end{enumerate}
\end{definition}

The following computation can be used to solve 
isomorphism and embedding problems for the class of 
Nambu-Poisson fields $N_{q}$.

\begin{lemma}
\label{xxlem4.12}
Let $q,q'\neq 0$.
\begin{enumerate}
\item[(1)]
$T(q,0)\cong T(q',0)$ if and only if $q=\pm q'$.
\item[(2)]
$T(q,0)$ can be embedded to $T(q',0)$ if and only if
$q/q'$ is an integer.
\item[(3)]
Any Nambu-Poisson algebra morphism $f: T(q,0)\to T(q,0)$ 
is an isomorphism.
\end{enumerate}
\end{lemma}

\begin{proof} (1,2) The proof is similar to the proof of Lemma 
\ref{xxlem4.7}(3). The key point is that only invertible 
elements in $T(q,0)$ are of the form $x_1^{m_1}\cdots x_n^{m_n}$ 
for some $m_1,\cdots,m_n\in {\mathbb Z}$. Then, by \eqref{E1.3.1},
if there is an embedding $f: T(q,0)\to T(q',0)$ then $q/q'$ is 
an integer. We skip the rest of the details.

(3) Let $y_i=f(x_i)$. Then up to a nonzero scalar
$y_i=x_1^{a_{i1}}\cdots x_n^{a_{in}}$ for all $i$. 
If $f: T(q,0)\to T(q,0)$ is not an isomorphism, 
then $\det[(a_{ij})_{n\times n}]\neq \pm 1$. Then, the assertion
follows from part (1) and \eqref{E1.3.1}.
\end{proof}

Part (3) of the above lemma says that $T(q,0)$ has the Dixmier
property [Definition \ref{xxdef5.3}]. We fix the canonical 
embedding $T(q,0)\to N_q$ for the rest of this subsection.

\begin{lemma}
\label{xxlem4.13}
Suppose $\kk$ is algebraically closed and $\OO={\mathbb Z}^n$. Let $N$ be the Nambu-Poisson 
field $N_q$ where $0\neq q\in \Bbbk$.
\begin{enumerate}
\item[(1)]
Let $\nu$ be a $0$-valuation on $N$. Then $\nu$ is determined 
by the collection of elements {\rm{(}}or vectors{\rm{)}} 
$\nu(x_1),\cdots, \nu(x_n)$ in $\OO$. Further, the 
corresponding filtration ${\mathbb F}$ on $T(q,0)$ agrees 
with ${\mathbb F}^{ind}$ which is determined by 
$\deg(x_i)=\nu(x_i)$ for all $i$. If we assign the grading on 
$T(q,0)$ by $\deg(x_i)=\nu(x_i)$, then 
$$Q_{gr}(\gr_{\nu} T(q,0))\cong \gr_{\nu} N
\cong Q_{\gr}(T(q,0))$$ 
as Nambu-Poisson $0$-graded domains.
\item[(2)]
$\nu$ is quasi-faithful if and only if ${\rm{rk}}_{\mathbb Z} 
(\sum_{i=1}^n \nu(x_i) {\mathbb Z})=n$. In this case, we have
\begin{equation}
\label{E4.13.1}\tag{E4.13.1}
\gr_{\nu} T(q,0)\cong \gr_{\nu} N\cong T(q,0)
\end{equation}
as Nambu-Poisson $0$-graded domains.
\item[(3)]
$\varrho(N)=\{q, -q\}$.
\item[(4)]
Let $\nu$ and $\mu$ be two distinct quasi-faithful 
$0$-valuations. Then there are monomials $a,b$ in 
$T(q,0)$ such that $\nu(a)>\mu(a)$ and $\nu(b)<\mu(b)$.
\item[(5)]
For a quasi-faithful $0$-valuation $\nu$, there is a unique
quasi-faithful $0$-valuation $\nu^{-}$ such that 
$\nu(f)+\nu^{-}(f)\leq 0$ for all $f\in T(q,0)$.
Further, $\nu(f)+\nu^{-}(f)=0$ for $f\in T(q,0)$ if and only 
if $f$ is a monomial $x_1^{m_1}\cdots x_n^{m_n}$ for some 
$m_s\in {\mathbb Z}$.
\end{enumerate}
\end{lemma}

\begin{proof}
(1) Let $\nu$ be a $0$-valuation of $N$ with corresponding 
filtration ${\mathbb F}$. Let $v_i=\nu(x_i)\in \OO$. Then $-v_i=
\nu(x_i^{-1})$. Then $T(q,0)$ is an $\OO$-graded domain with 
grading determined by $\deg(x_i)=v_i$ and $\deg(x_i^{-1})=-v_i$ 
for all $i$. This grading produces a good induced filtration 
${\mathbb F}^{ind}$ via \eqref{E2.8.3}. Note that, by Lemma 
\ref{xxlem2.12}, ${\mathbb F}^{ind}$ also agrees with the 
Adams$^{Id}$ filtration ${\mathbb F}^{Id}$ given in 
\eqref{E2.11.1}. It is clear that ${\mathbb F}^{ind}$
is a $0$-filtration since $T(q,0)$ itself is a Nambu-Poisson 
$0$-graded domain. Further $\gr_{{\mathbb F}^{ind}} 
T(q,0)\cong T(q,0)$ (as $T(q,0)$ is a $\OO$-graded algebra
to start with). By Lemma \ref{xxlem2.11}(2), 
$\phi^{ind}: T(q,0)\to \gr_{\nu} T(q,0)$ is a Nambu-Poisson 
algebra homomorphism. Since $T(q,0)$ is Nambu-Poisson simple, 
$\phi^{ind}$ is injective and the image
of $\phi^{ind}$ has GK-dimension $n$. By Lemma \ref{xxlem2.10}(1),
${\mathbb F}^{ind}={\mathbb F}$. Hence
$Q_{gr}(\gr_{\nu} T(q,0))\cong Q_{gr}
(\gr_{{\mathbb F}^{ind}}T(q,0))$. By Lemma \ref{xxlem2.7}(1),
\begin{equation}
\label{E4.13.2}\tag{E4.13.2}
\gr_{\nu} (N)=Q_{gr}(\gr_{\nu} T(q,0))
\cong Q_{gr}(T(q,0)).
\end{equation}
Therefore, the assertion follows.

(2) The first assertion follows from part (1). The second
assertion follows from part (1) and the fact that $T(q,0)$
is a graded field (whence $Q_{gr}(T(q,0))=T(q,0)$) when the 
support ${\mathbb S}(T(q,0))$ has rank $n$ over ${\mathbb Z}$ by Lemma \ref{xxlem1.7}(2).

(3) It is clear that $q\in \varrho(N)$. Since $T(q,0)\cong T(-q,0)$
by the isomorphism $\phi: x_1,\to x_2, x_2\to x_1,x_i\to x_i$
for all $i\geq 3$. Then $-q\in \varrho(N)$. Since \eqref{E4.13.1}
holds for all quasi-faithful $0$-valuations $\nu$, we obtain
that $\varrho(N)=\{q,-q\}$ by Lemma \ref{xxlem4.12}(1) as required.

(4) Suppose to the contrary that $\mu(a)\leq \nu(a)$ for all
monomials $a\in T(q,0)=:T$. Let ${\mathbb G}$ (resp. ${\mathbb F}$) 
be the filtration corresponding to $\mu$ (resp. $\nu$). 
Since $\mu$ (resp. $\nu$) is a quasi-faithful $0$-valuation, part (2) and Lemma \ref{xxlem1.7}(2) implies that if $a$ and $b$ are distinct monomials we have $\mu(a)\neq \mu(b)$ (resp. $\nu(a)\neq \nu(b)$). Hence $\mu(a)\leq \nu(a)$ for all $a\in T$ and ${\mathbb G}$ is a 
subfiltration of ${\mathbb F}$. By the proof of Lemma \ref{xxlem2.8},
there is an $\OO$-graded algebra homomorphism 
$$\phi: \gr_{\mu}(T)\to \gr_{\nu}(T).$$
We claim that $\phi$ preserves the Nambu-Poisson structure.
For $f\in G_i(T)$ (resp. $f\in F_i(T)$), 
let $\overline{f}$ (resp. $\widetilde{f}$) be the element 
$f+G_{>i}(T)$ (resp. $f+ F_i(T)$). Then 
$$\begin{aligned}
\phi(\{\overline{f_1},\cdots,\overline{f_n}\})
&=\phi(\overline{\{f_1,\cdots,f_n\}})\\
&=\widetilde{\{f_1,\cdots,f_n\}}\\
&=\{\widetilde{f_1},\cdots, \widetilde{f_n}\}.
\end{aligned}
$$
Hence, $\phi$ is a Nambu-Poisson algebra homomorphism.
Since $\gr_{\mu}(T)$ is simple, $\phi$ is injective. 
By Lemma \ref{xxlem2.8}(2), ${\mathbb G}={\mathbb F}$,
yielding a contradiction. 

(5) Let $\nu$ a quasi-faithful $0$-valuation of $N$. By 
part (1), $\nu$ is determined by the collection of 
vectors $[\nu(x_1),\cdots, \nu(x_n)]$. 

Now we define a grading on $T(q,0)$ by $\deg(x_i)=
-\nu(x_i)$. Then $T(q,0)$ is a Nambu-Poisson 
$0$-graded domain. Let $\nu^{-}$ be the 
$0$-valuation $\nu^{Id}$ associated to this grading 
as given in Lemma \ref{xxlem2.12}(1). Then, it is routine 
to check that all properties in part (5) are satisfied
for this pair $(\nu,\nu^{-})$. It remains to show the 
uniqueness of $\nu^{-}$. Let $\mu$ be another 
$0$-valuation that is different from $\nu^{-}$.
By part (4), there is a monomial $a$ such that
$\mu(a)>\nu^{-}(a)$. Then 
$$\nu(a)+\mu(a)> \nu(a)+\nu^{-}(a)=0$$
which shows that $(\nu,\mu)$ does not satisfy
the condition specified.
\end{proof}

Our first application of $\delta$-type invariant is the
following.

\begin{theorem}
\label{xxthm4.14}
Let $q,q'\neq 0$. Then $N_q\cong N_{q'}$ if and only 
if $q=\pm q$. 
\end{theorem}

\begin{proof}
If $q'=\pm q$, by Lemma \ref{xxlem4.12}(1), $T(q,0)\cong T(q',0)$.
After localization, we obtain that $N_q\cong N_{q'}$. The 
converse follows from Lemma \ref{xxlem4.13}(3).
\end{proof}

We can also use $\delta$-type invariants to determine an 
embedding problem.

\begin{lemma}
\label{xxlem4.15} 
Suppose $\kk$ is algebraically closed. Let $N$ and $Q$ be two Nambu-Poisson fields such that there is an embedding $N\to Q$. Then, for every scalar 
$a\in \varrho(Q)$, there is an $b\in \varrho(N)$ such that 
$a/b$ is an integer.
\end{lemma}

\begin{proof}
In this proof, we take $\OO={\mathbb Z}^n$. By the definition 
of $\varrho(Q)$, there is a quasi-faithful $0$-valuation $\nu$ 
on $Q$ such that $\gr_{\nu}(Q)\cong T(a,0)$. Let $\mu$ be the 
restriction of $\nu$ on $N$. By Lemma \ref{xxlem2.19}, $\mu$ 
is a quasi-faithful $0$-valuation on $N$. The embedding $N\to Q$ 
induces an injective map $\gr_{\mu} (N)\to \gr_{\nu}(Q)$. 
By Lemma \ref{xxlem4.5} and the proof of Lemma \ref{xxlem4.6}, 
$\gr_{\mu}(N) \cong T(b,0)$ for some $b$. So $b\in \varrho(N)$ 
and $T(a,0)\to T(b,0)$ is an embedding. By Lemma \ref{xxlem4.12}, 
$a/b$ is an integer as required. 
\end{proof}

Here is a solution to the embedding problem for the class
$N_q$ for all $q\in \Bbbk^{\times}$. 

\begin{theorem}
\label{xxthm4.16}
Let $q,q'\neq 0$. Then $N_q\hookrightarrow N_{q'}$ if and only 
if $q/q'$ is an integer. 
\end{theorem}

\begin{proof} 
If $d:=q/q'$ is an integer, then $T(q,0) \hookrightarrow T(q',0)$ by sending $x_1\to x_1^{d}$ and $x_i\to x_i$ for all $i\geq 2$. So $N_q$ embeds into $N_{q'}$. 

Conversely, suppose that $N_q\hookrightarrow N_{q'}$. Without loss of generality, we can assume $\kk$ to be algebraically closed. By Lemma \ref{xxlem4.13}(3), $q' \in \varrho(N_{q'})$. By Lemma 
\ref{xxlem4.15} there is an $a\in \varrho(N_q)$ such that 
$a/q'$ is an integer. By Lemma \ref{xxlem4.13}(3), $a$ is either 
$q$ or $-q$. Hence, $q/q'$ is an integer. So, the assertion follows.
\end{proof}

\subsection{$\varsigma$-type invariants}
\label{xxsec4.5}
In this subsection, we propose a class of $\varsigma$-type 
invariants. By definition, these are invariants related to
the moduli of a certain family of valuations. In other words, 
these are geometric invariants. At this point, we don't
have a formal language to deal with the moduli problem
of valuations. Hopefully, the initial examples should indicate that these invariants are reasonable. 

\begin{definition}
\label{xxdef4.17}
Let $\OO={\mathbb Z}$ and let $N$ be a Nambu-Poisson
field. If the class valuations ${\mathcal V}_{fw}(N)$
is parametrized by a space {\rm{(}}or an algebraic 
variety{\rm{)}} $X$, then this space is denoted by 
$\varsigma_{fw}(N)$.
\end{definition}

For the rest of this subsection, we assume that 
$\OO={\mathbb Z}$. We have the following lemma about 
$1$-valuations.

\begin{lemma}
\label{xxlem4.18}
Let $B$ be a Nambu-Poisson algebra that is an affine normal 
domain. Let $\nu$ be a valuation on $N:=Q(B)$ such that 
{\rm{(a)}} $F_0(B)=B$ and {\rm{(b)}} $\GKdim F_0(B)/F_1(B)
=\GKdim B-1$. Then, $\nu$ is equivalent to a $1$-valuation.
\end{lemma}

\begin{proof} By definition $F_0(B)/F_1(B)$ is a domain.
Hence $F_1(B)$ is a prime ideal of $B(=F_0(B))$. Let 
$S=F_0(B)\setminus F_1(B)$. Then, every element in $S$ has 
valuation $0$. Then, the valuation of every element in 
$A:=B_{S^{-1}}$ is nonnegative. As a consequence, 
$F_0(A)=A$. Since $F_1(B)$ is prime of height 1, $A$ is 
local and regular of Krull dimension 1. This means that 
$A$ is a DVR with a generator $\omega$ of the maximal ideal 
of $A$.  By the proof of \cite[Lemma 2.10(3)]{HTWZ3}, we have 
that there is a positive integer $h$ such that
$$F_i(A)=\begin{cases} A & i\leq 0,\\
\omega^{\lceil i/h \rceil} A & i>0.\end{cases}
$$
Then $\nu$ is equivalent to another valuation, denoted 
by $\nu'$, with filtration defined by
$$F'_i(A)=\begin{cases} A & i\leq 0,\\
\omega^{i} A & i>0.\end{cases}
$$
By Theorem \ref{xxthm2.16}(1), the $\nu'$ is $1$-valuation. 
The assertion follows.
\end{proof}

In the next proposition, we work out some faithful 
$n$-valuations on some Nambu-Poisson fields. Let $\Reg B$ denote 
the regular locus of a commutative algebra $B$. If $B$ is affine, 
then $\Reg B$ is an open variety of $\Spec B$. If $B$ is a 
Nambu-Poisson algebra, let $B_{cl}$ be the factor ring 
$B/(\{B,\cdots, B\})$.

\begin{proposition}
\label{xxpro4.19}
Suppose $\Bbbk$ is algebraically closed.
Let $B$ be a Nambu-Poisson algebra that is an affine domain of 
GK-dimension $n$. Let $N$ be the Nambu-Poisson fraction field 
$Q(B)$. Let $\fm$ be a maximal ideal in $\Reg B$. 
\begin{enumerate}
\item[(1)] 
There is a unique faithful $n$-valuation, denoted by $\nu^{n,\fm}$, 
of $N$ such that $F_i^{\nu}(B)=\fm^i$ for all $i\geq 0$.
\item[(2)]
$\nu^{n,\fm}$ is classical if and only of $\fm\in \Spec B_{cl}$.
\item[(3)]
$\nu^{n,\fm}$ is either classical or Weyl.
\end{enumerate}
\end{proposition}

\begin{proof}
(1) Since $\fm\in \Reg B$, the localization $A:=B_{\fm}$ is a 
regular local ring with maximal ideal, denoted by $\fp$. Define 
a filtration 
${\mathbb F}^{n,\fm}:=\{F^{n,\fm}_i\mid i\in {\mathbb Z}\}$ on 
$A$ by
\begin{equation}
\label{E4.19.1}\tag{E4.19.1}
F^{n,\fm}_i(A)=\begin{cases} \fp^i & i\geq 0.\\
A& i<0.\end{cases}
\end{equation}
Since $A$ is regular local, $\gr_{{\mathbb F}^{n,\fp}} A$ 
is isomorphic to a polynomial $\Bbbk[t_1,\cdots,t_n]$
where $\Bbbk$ is the residue field in this case. Hence 
${\mathbb F}^{n,\fm}$ is a good filtration. The uniqueness
follows from the facts that $\gr_{\mathbb F}B\cong 
\gr_{\mathbb F} A\cong \Bbbk[t_1,\cdots,t_n]$ and that 
$B\cap \fp^i=\fm^i$.

Next, we show that the valuation associated with 
${\mathbb F}^{n,\fm}$ is an $n$-valuation. Let 
$S:=\{s_1,\cdots,s_n\}$ be a regular system of parameters.
Then $\fm=\sum_{\alpha=1}^n s_{\alpha} A$. As a consequence, 
for all $i\geq 0$, $F^{n,\fm}_i=\sum_{\alpha_1,\cdots,\alpha_i}
s_{\alpha_1}\cdots s_{\alpha_i} A$. Using the fact that
$\{s_{\alpha_1}, s_{\alpha_2},
\cdots, s_{\alpha_n}\}\in A$ and $\{A,\cdots,A\}\subseteq A$,
we obtain that $\{F^{n,\fm}_{i_1},\cdots, F^{n,\fm}_{i_n}\}
\subseteq F^{n,\fm}_{\sum_{s=1}^{n} i_{s}-n}$. Therefore
${\mathbb F}^{n,\fm}$ is a good $n$-filtration. We use $\nu^{n,\fm}$
to denote the associated $n$-valuation.

By the above argument, there is 
one-to-one correspondence between $\Reg B$ and the set 
of Nambu-Poisson $n$-valuations $\{\nu^{n,\fm}\mid 
\fm\in \Reg B\}$.

(2) Since $\gr_{\mathbb F} B$ is generated in degree 1, we
have 
$$\begin{aligned}
{\text{$\nu^{n,\fm}$ is classical}}
&\Leftrightarrow
{\text{ $\{F_1(B),\cdots,F_1(B)\}\subseteq F_{n-n+1}(B)$}}\\
&\Leftrightarrow
{\text{ $\{\fm,\cdots,\fm\}\subseteq \fm$}}\\
&\Leftrightarrow
{\text{ $\{\Bbbk+\fm,\cdots,\Bbbk+\fm\}\subseteq \fm$}}\\
&\Leftrightarrow
{\text{ $\{B,\cdots,B\}\subseteq \fm$}}\\
&\Leftrightarrow
{\text{ $\fm\in \Spec B_{cl}$.}}
\end{aligned}
$$

(3) Since $B$ has GK-dimension $n$, $\gr_{\mathbb F} 
B$ is isomorphic to $\Bbbk[t_1,\cdots,t_n]$ where $t_1,
\cdots, t_n$ are in degree 1. If $\nu^{n,\fm}$ is not 
classical, then $\{t_1,\cdots,t_n\}\neq 0$. Since 
$\gr_{\mathbb F} B$ is Nambu-Poisson $n$-graded, 
$\{t_1,\cdots,t_n\}\in \Bbbk^{\times}$.
So $\nu^{n,\fm}$ is Weyl.
\end{proof}

\begin{example}
\label{xxexa4.20} 
In this example, $n=2$, that means that we are working on 
ordinary Poisson algebras instead of Nambu-Poisson algebras.
Assume that $\Bbbk$ is algebraically closed. Let $d_0\geq 3$
and let $\Omega=t_1^{3+d_0}+t_2^{3+d_0}+t_3^{3+d_0}$ which 
is an i.s. potential of degree $\geq 6$ in $\Bbbk^{[3]}$. Let 
$B$ be $P_{\Omega-1}$ as defined in Construction \ref{xxcon0.5}. 
Let $N:=Q(B)$. We will prove that there is a natural one-to-one 
correspondence between faithful $2$-valuation $\nu$ and the 
maximal ideal $\fm\in \MaxSpec B$, namely, 
\begin{equation}
\label{E4.20.1}\tag{E4.20.1}
{\mathcal V}_{f2}(N)\longleftrightarrow \MaxSpec B=
\MaxSpec {^1\Gamma_{0}}(N).
\end{equation}
Under this setting, it is reasonable to conjecture that all
faithful $2$-valuations on $N$ are ``parametrized'' by the variety
$\MaxSpec B$ and that $\varsigma_{f2}(N)=\MaxSpec B$. 

Note that \eqref{E4.20.1} fails for $N_q$ or $N_{Weyl}$. 
\end{example}

\begin{proof}[Proof of the statement in Example \ref{xxexa4.20}]
Note that $B$ is regular [Lemma \ref{xxlem3.3}(2)].
Hence, every element in $\MaxSpec B$ is in $\Reg B$.
For every $\fm\in \MaxSpec B$, by Proposition 
\ref{xxpro4.19}(1), there is a unique faithful
$2$-valuation $\nu^{2,\fm}$ such that $F_i^{\nu}(B)
=\fm^i$ for all $i\geq 0$. This gives a map 
from $\MaxSpec B\to {\mathcal V}_{f2}(N)$. 

Conversely, let $\nu$ be a faithful $2$-valuation of $N$. 
It remains to show that $\nu=\nu^{2,\fm}$ for some 
$\fm\in \MaxSpec B$. By Theorem \ref{xxthm3.11}, 
${^2\Gamma_{0}}(N)=B$. This implies that $F^{\nu}_0(B)
=B$. Let $\fp=F^{\nu}_1(B)$ which must be a prime 
ideal of $B$ as $F^{\nu}_0(B)/F^{\nu}_1(B)$ is a domain. 
Since $\nu$ is nontrivial, $\fp\neq 0$. If $\fp$ has height
one, then Lemma \ref{xxlem4.18} implies that 
$\nu$ is equivalent to a $1$-valuation, contradicting
with $\nu$ being faithful. Therefore $\fp\in \MaxSpec B$.
It suffices to show that $\nu=\nu^{2,\fp}$. 

Let $(a,b,c)=(\nu(t_1),\nu(t_2),\nu(t_3))$ and let $I$ be the 
subalgebra of $\gr_{\nu}(B)$ generated by $\overline{t_1}$,
$\overline{t_2}$, and $\overline{t_3}$. By symmetry, we may 
assume that $a\leq b\leq c$. Since $F^{\nu}_0(B)=B$, $a\geq 0$. 
Since $\Omega-\xi=0$ in $B$, we get $a=0$. Moreover, $\nu$ being a 2-valuation implies that 
$$0=\nu((3+d_0) t_1^{2+d_0})
=\nu(\{t_2,t_3\})\geq \nu(t_2)+\nu(t_3)-2=b+c-2$$
which implies that $b+c\leq 2$. So we only need to consider the 
following four cases for $(a,b,c)$: $(0, 0, 0)$, $(0, 0, 1)$, 
$(0, 0, 2)$, $(0, 1, 1)$.

Case 1: $(a,b,c)=(0,0,0)$. We have already shown that $F_1^{\nu}(B)
\in \MaxSpec B$. In this case this is equivalent to $I=\Bbbk$. 
Hence $\overline{t_i}=c_i\in \Bbbk^{\times}$ for $i=1,2,3$. Let 
$y_i=t_i-c_i$ for $i=1,2,3$.  Since
$$0=\nu((3+d_0)t_3^{2+d_0})=\nu(\{t_1, t_2\})
=\nu(\{y_1,y_2\}) \geq \nu(y_1)+\nu(y_2)-2,$$
we obtain that $\nu(y_1)=\nu(y_2)=1$. Similarly,
$\nu(y_3)=1$. 

Let $J$ be the subalgebra of $\gr_{\nu} B$ generated by
$\overline{y_1}$, $\overline{y_2}$, and $\overline{y_3}$. 
If $\GKdim J\leq 1$, then by \cite[Lemma 1.3(3)]{HTWZ1}, 
$\overline{y_1}=d \overline{y_2}$ for some $d\in \Bbbk^{\times}$,
or equivalently, $y_1=dy_2 +z_1$ with $\nu(z_1)\geq 2$.
Then 
$$0=\nu((3+d_0)t_3^{2+d_0})=\nu(\{t_1,t_2\})=
\nu(\{y_1,y_2\})=\nu(\{z_1,y_2\})
\geq \nu(z_1)+\nu(y_2)-2>0,$$
yielding a contradiction. Therefore $\GKdim J=2$. 
By recycling the notation, let ${\mathbb F}^{ind}$ be 
the (new) induced filtration determined
by $\deg(y_1)=\deg(y_2)=\deg(y_3)=1$. Then $\nu$ 
is uniquely determined by $(\overline{t_1},\overline{t_2},
\overline{t_3})\in \Bbbk^3$ which must satisfy the condition
$(\overline{t_1})^{3+d_0}+(\overline{t_2})^{3+d_0}
+(\overline{t_3})^{3+d_0}-1=0$, see Proposition \ref{xxpro4.19}. 
Therefore, $\nu$ is uniquely determined by 
$$(\overline{t_1},\overline{t_2},
\overline{t_3})\in \MaxSpec B,$$
namely, $\nu^{2,\fp}$ where $\fp=(\overline{t_1},\overline{t_2},
\overline{t_3})$.

Case 2: $(a,b,c)=(0, 0, 1)$. Then, the induced filtration
${\mathbb F}^{ind}$ agrees with the one given in Lemma 
\ref{xxlem4.18} with $F_1^{ind}(B)=(t_3)$. If $\GKdim I=2$, 
then, by Lemma \ref{xxlem4.18}, ${\mathbb F}^{ind}$ is a 
filtration that is equivalent to $1$-valuation. So $\nu$ is 
not a faithful $2$-valuation, yielding a contradiction. It 
remains to consider that case when $\GKdim I\leq 1$. Then 
$I=\Bbbk[\overline{t_3}]$ by \cite[Lemma 1.3(3)]{HTWZ1}.
As a consequence $t_1=c_1+y_1$ and $t_2=c_2+y_2$
where $c_1,c_2\in \Bbbk^{\times}$ and $\nu(y_1),
\nu(y_2)\geq 1$. Then 
$$0=\nu((3+d_0)t_2^{2+d_0})=\nu(\{t_3,t_1\})=
\nu(\{t_3,y_1\})
\geq \nu(t_3)+\nu(y_1)-2,$$
which implies that $\nu(y_1)=1$. Similarly,
$\nu(y_2)=1$. Further $c_1^{3+n}+c_2^{3+n}-1=0$.
By re-cycling notations, let ${\mathbb F}^{ind}$ be 
the induced determined by $\deg(y_1)=\deg(y_2)=\deg(t_3)=1$. 
If $\GKdim I=2$, by Proposition \ref{xxpro4.19}, there is 
a unique $2$-valuation determined by 
$\fp:=(\overline{t_1},\overline{t_2},0)\in \Bbbk^3$ 
with $\fp\in \MaxSpec B$ (or with 
$\overline{t_1}^{3+d_0}+\overline{t_2}^{3+d_0}-1=0$).

If $\GKdim I\leq 1$, then $t_1=c_1+d_1 t_3+\gamma$ and 
$y=c_2+d_2 t_3+\delta$ where $c_1,d_1,c_2,d_2\in 
\Bbbk^{\times}$ and $\nu(\gamma), \nu(\delta)\geq 2$. Then 
$$0=\nu((3+d_0)t_2^{2+d_0})=\nu(\{t_3,t_1\})
=\nu(\{t_3,\gamma\})\geq \nu(t_3)+\nu(\gamma)-2>0,$$
yielding a contradiction. 

Case 3: $(a,b,c)=(0,0,2)$. Then, the induced 
filtration agrees with the one given in Lemma \ref{xxlem4.18} 
(with $\deg(t_3)=2$). If $\GKdim I=2$, then, by Lemmas 
\ref{xxlem2.10}(1) and \ref{xxlem4.18}, $\nu$ is not a 
faithful $2$-valuation, yielding 
a contradiction. It remains to consider that case when 
$\GKdim I\leq 1$. Then $I=\Bbbk[s]$ by 
\cite[Lemma 1.3(3)]{HTWZ1}. As a consequence, $t_1=c_1+y_1$ 
and $t_2=c_2+y_2$
where $c_1,c_2\in \Bbbk^{\times}$ and $\nu(y_1),
\nu(y_2)\geq 1$. Then 
$$0=\nu((3+d_0)t_2^{2+d_0})=\nu(\{t_3,t_1\})
=\nu(\{t_3,y_1\})
\geq \nu(t_3)+\nu(y_1)-2>0,$$
yielding a contradiction. So, in this case there is no 
faithful $2$-valuation.

Case 4: $(a,b,c)=(0, 1, 1)$. If $\GKdim I=2$, we
have an induced filtration ${\mathbb F}^{ind}$ 
that agrees with ${\mathbb F}$. In this case 
$$0=\Omega-\xi=t_1^{3+d_0}-\xi$$
in $F^{\nu}_0(B)/F^{\nu}_1(B)$. Then 
$\overline{t_2}=\overline{t_3}=0$ and 
$\overline{t_1}=  c_1\in \Bbbk^{\times}$
satisfying $c_1^{3+d_0}=\xi$. Then ${\mathbb F}^{ind}$
is corresponding to $\nu^{2, \fm}$ where 
$\fm=(c_1,0,0)\in \MaxSpec B$. Moreover, an easy 
computation shows that $\{\overline{t_2},\overline{t_3}\}
=(3+d_0)c_1^{2+d_0}$. So $Q(\gr_{{\mathbb F}^{ind}} P)\cong 
N_{Weyl}$. By Lemma \ref{xxlem2.10}(1), $\nu=\nu^{2,\fm}$.

If $\GKdim I\leq 1$, then $t_2=d t_3+\alpha$ 
where $\nu(\alpha)>1$. Then 
$$0=\nu((3+d_0)t_1^{2+d_0})=\nu(\{t_2,t_3\})=
\nu(\{\alpha,t_3\})=
\geq \nu(\alpha)+\nu(t_3)-2>0,$$
yielding a contradiction. This finishes the proof.
\end{proof}

\section{Applications}
\label{xxsec5}

For the rest of the paper, we give several applications of 
valuations. To save space, we won't provide the most general results but rather provide some essential ideas on using valuations to solve problems. Hence, each topic will be 
very brief. Throughout this section, we assume that 
$\OO={\mathbb Z}$ and that
\begin{equation}
\label{E5.0.1}\tag{E5.0.1}
{\text{$\Omega$ is an i.s. potential in $\Bbbk^{[n+1]}$ of
degree $n+1+d_0$ for $d_0\geq 0$.}}
\end{equation}

\subsection{Automorphisms of Nambu-Poisson fields}
\label{xxsec5.1}
In this section, we compute the automorphism group of some 
Nambu-Poisson fields.

\begin{lemma}
\label{xxlem5.1}
Let $\Omega$ be as in \eqref{E5.0.1} with $d_0\geq 1$ and 
let $B$ be $P_{\Omega-\xi}$ as defined in Construction 
\ref{xxcon0.5} where $\xi\in \Bbbk$. Let $N:=Q(B)$ and    
$\nu$ be a $d_0$-valuation on $N$ with associated filtration $\mathbb F$. Suppose $\nu(f)<0$ for 
some $f\in B$. 
\begin{enumerate}
\item[(1)]
If $\xi=0$, then $\nu=\nu^{-Id}$ as given in Lemma 
\ref{xxlem2.12}(1) which is defined by the the filtration 
${\mathbb F}^{-Id}$ that is determined by the degree 
assignment $\deg(t_s)=-1$ for all $s$. As a consequence, 
$\gr_{\mathbb F} B\cong B$.
\item[(2)]
If $\xi\neq 0$, then $\nu$ agrees with the valuation $\nu^{c}$ 
given in Lemma \ref{xxlem2.12}(3) which is defined by the 
filtration ${\mathbb F}^{c}$ that is determined by the degree 
assignment $\deg(t_s)=-1$ for all $s$. As a consequence, 
$\gr_{\mathbb F} B\cong P_{\Omega}$.
\end{enumerate}
In both cases, $Q(\gr_{\nu}N)\not\cong N_{Weyl}$.
As a consequence, 
\begin{enumerate}
\item[(3)]
$\nu(t_s)=-1$ for all $s=1,\cdots,n+1$ and $F^{\nu}_0(P)=\Bbbk$.
\end{enumerate}
\end{lemma}

\begin{proof}
By valuation axioms, one of $\nu(t_s)$ must be negative. By 
definition, $\rho<0$ (see Notation \ref{xxnot3.4}(3) for the
definition of $\rho$).

(1,2) Since $\rho<0$, the assertions follows from 
lemma \ref{xxlem3.5}(4).

(3) By the proof of parts (1,2), $F_0(B)=\Bbbk$.
\end{proof}

Let $B$ be a filtered algebra with $F_0(B)=\Bbbk$ and $\sigma$ 
be an algebra automorphism. We say $\sigma$ is {\it linear} 
(with respect to the filtration) if 
$\sigma(F_{-1}(B))\subseteq F_{-1}(B)$. 

\begin{proposition}
\label{xxpro5.2}
Let $\Omega$ be as in \eqref{E5.0.1} where $d_0\geq 2$. Let 
$P_{\Omega-\xi}$ be defined as in Construction \ref{xxcon0.5}. 
Consider the natural filtration on $P_{\Omega-\xi}$ determined 
by the degree assignment $\deg(t_s) = -1$ for all $s$. 
\begin{enumerate}
\item[(1)]
$\Aut_{Poi}(Q(P_{\Omega}))=\Aut_{Poi}(P_{\Omega})$ and 
$\sigma\in \Aut_{Poi}(P_{\Omega})$ is graded.
\item[(2)]
Every Nambu-Poisson algebra automorphism of $P_{\Omega}$ lifts 
to a Nambu-Poisson algebra automorphism of $A_{\Omega}$.
\item[(3)]
If $\xi\neq 0$, then $\Aut_{Poi}(Q(P_{\Omega-\xi}))=
\Aut_{Poi}(P_{\Omega-\xi})$ and $\sigma\in 
\Aut_{Poi}(P_{\Omega-\xi})$ is linear. Further,
$\sigma(t_s)\in V$ where $V=\sum_{s=1}^{n+1} \Bbbk t_s$.
\item[(4)]
If $\xi\neq 0$, then every Nambu-Poisson algebra automorphism of
$P_{\Omega-\xi}$ lifts to a Nambu-Poisson algebra automorphism of
$A_{\Omega}$.
\end{enumerate}
\end{proposition}

\begin{proof}
(1) It is clear that $\Aut_{Poi}(P_{\Omega})\subseteq 
\Aut_{Poi}(Q(P_{\Omega}))$ by localization. For the opposite 
inclusion, let $\sigma\in \Aut_{Poi}(Q(P_{\Omega}))$. Since 
$P_{\Omega}=^1\Gamma_0(Q(P_{\Omega}))$ [Theorem \ref{xxthm3.10}], 
$\sigma$ restricts to a Poisson algebra automorphism of 
$P_{\Omega}$. Therefore, the first assertion follows. Let $\nu$ 
be a faithful $d_0$-valuation on $Q(P_{\Omega})$ such that 
$F^{\nu}_0(P_{\Omega})\neq P_{\Omega}$. Then $\nu(f)<0$ for 
some $f\in P_{\Omega}$. By Lemma \ref{xxlem5.1}(1), 
$\nu=\nu^{-Id}$, and consequently, such an $\nu$ is unique. 
By the uniqueness of such $\nu$, $\nu\circ \sigma=\nu$. 
This implies that $\nu(\sigma(f))=-1$ whenever $\nu(f)=-1$, 
or equivalently, $\sigma$ is linear. For each $i$, write 
$\sigma(t_i)=y_i+a_i$ where $y_s\in V=\sum_{s=1}^{n+1} \Bbbk t_s$
and $a_s\in \Bbbk$. Applying $\sigma$ we have
$$\begin{aligned}
\Omega_{t_i}&(y_1+a_1,\cdots,y_{n+1}+a_{n+1})
=\sigma(\Omega_{t_i})\\
&=\sigma((-1)^{n+1-i}\{t_1,\cdots, \widehat{t_i},\cdots,t_{n+1}\})\\
&=(-1)^{n+1-i}\{y_1+a_1,\cdots, \widehat{y_i+a_i},\cdots,y_{n+1}+a_{n+1}\}\\
&=(-1)^{n+1-i}\{y_1,\cdots, \widehat{y_i},\cdots,y_{n+1}\}
\in W
\end{aligned}
$$
where $W=\sum_{s=1}^{n+1}\Bbbk \Omega_{t_s}$. So the constant 
term of $\Omega_{t_i}(y_1+a_1,\cdots,y_{n+1}+a_{n+1})$ is 0. 
This implies $\Omega_{t_i}(a_1,\cdots,a_{n+1})=0$ for each $i$.
Since $\Omega$ has an isolated singularity at the origin, 
$a_i=0$ for all $i$. Thus, $\sigma$ is graded. 

(2) By part (1), every Nambu-Poisson algebra automorphism 
$\sigma$ of $P_{\Omega}$ is graded and hence is determined by 
$\sigma(t_1),\ldots,\sigma(t_{n+1})\in (P_{\Omega})_{-1}$ such 
that the following identity holds 
\[
(-1)^{n+1-i}\{\sigma(t_1),\ldots,\widehat{\sigma(t_i)},\ldots,
\sigma(t_{n+1})\}=\sigma(\Omega_{t_i})
\]
for all $1\le i\le n+1$. Note the above identities hold in 
$(P_{\Omega})_{-(n+d_0)}=(A_\Omega)_{-(n+d_0)}$. Since
the only relation of $P_{\Omega}$ is in degree $n+1+d_0$, 
$\sigma$ can be lifted to a unique graded Nambu-Poisson 
algebra automorphism of $A_\Omega$. 

(4,3) The proofs are similar to the proofs of parts (1,2).
\end{proof}

\subsection{Dixmier property}
\label{xxsec5.2}
We refer to \cite[Section 6.1]{HTWZ3} for a brief discussion 
about Dixmier conjecture and Dixmier property of noncommutative algebras and 
Poisson algebras. 

\begin{definition}
\label{xxdef5.3}
Let $B$ be a Nambu-Poisson algebra. We say $B$ satisfies the 
{\it Dixmier property} if every injective Nambu-Poisson algebra 
morphism $f: B\to B$ is bijective.
\end{definition}

Our main result in this subsection is the following.

\begin{theorem}
\label{xxthm5.4}
Let $\Omega$ be as in \eqref{E5.0.1} where $d_0\geq 2$.
Let $B$ be either $P_{\Omega}$ or $P_{\Omega-\xi}$ defined 
as in Construction \ref{xxcon0.5}. Both $B$ and $Q(B)$ 
satisfy the Dixmier property.
\end{theorem}

\begin{proof} Let $N=Q(B)$ and let $f: N\to N$ be an 
(injective) Nambu-Poission algebra homomorphism. By Theorem 
\ref{xxthm3.10}, ${^1\Gamma_{0}}(N)=B$. By Lemma \ref{xxlem2.15}(1), 
$f$ restricts to a map $f\mid_{B}: B\to B$. Since $N=Q(B)$, 
it remains to show that $B$ satisfies the Dixmier property. 
Let $\nu$ be a $d_0$-valuation on $N$ such that $F_0(B)\neq B$. 
By Lemma \ref{xxlem5.1}, (a) $\nu=\nu^{-Id}$, (b) 
$\nu(t_i)=-1$ for all $i=1,\cdots, n+1$ and (c) $F_0(B)=\Bbbk$.
Let $\mu$ be the $d_0$-valuation $\nu\circ f\mid_{B}$ on $B$ 
with filtration ${\mathbb F}'$. Then $F'_0(B)=\Bbbk$. Lemma 
\ref{xxlem5.1} implies that $\mu=\nu$. Since $\mu=\nu\circ f$,  the 
$\nu$-values of $f(t_i)$ is $-1$ for all $i$, or equivalently, 
$f(t_1),\ldots,f(t_{n+1})$ are in $V$. Since $f$ is injective, 
we obtain that $f(V)=V$. Since $B$ is generated by $V$, $f$ is 
bijective as required. 
\end{proof}

\subsection{Rigidity of grading}
\label{xxsec5.3}

In this subsection, we prove the following.

\begin{theorem}
\label{xxthm5.5}
Let $\Omega$ be as in \eqref{E5.0.1} where $d_0\geq 2$. Let 
$P_{\Omega}$ be defined as in Construction \ref{xxcon0.5}.
Then $P_{\Omega}$ has a unique connected grading 
such that it is Nambu-Poisson $d_0$-graded.
\end{theorem}

\begin{proof} Let $B$ denote $P_{\Omega}$. By Construction 
\ref{xxcon0.5}, there is a grading (which is the opposite 
of the Adams grading) such that $B$ is Nambu-Poisson 
$d_0$-graded. 

Suppose $B$ has a new connected grading such that $B$ is 
Nambu-Poisson $d_0$-graded. By definition, this is a 
non-positive grading, and its opposite grading is a non-negative 
grading called a new Adams grading. Let $\mu$ be the 
$d_0$-valuation associated with the new Adams grading. In this 
case $\mu(f)<0$ for some $f\in B$. By Lemma \ref{xxlem5.1}(1,3), 
$\mu=\nu^{-Id}$ and $\mu(t_i)=-1$ for $s=1,\cdots,n+1$. Write 
$t_i=y_i+a_i$ where $y_i$ is homogeneous of new grading $-1$ 
and $a_i\in \Bbbk$ for all $i$. Since every linear combination 
of $t_1,\cdots, t_{n+1}$ has $\mu$-value $-1$, $y_1,\cdots,
y_{n+1}$ are linearly independent. Now
$$\begin{aligned}
0&=\Omega(t_1,\cdots, t_{n+1})
=\Omega(y_1+a_1,\cdots, y_{n+1}+a_{n+1})\\
&=\Omega(a_1,\cdots, a_{n+1})
+\sum_{s=1}^{n+1}
\Omega_{t_s}(a_1,\cdots,a_{n+1}) y_s+ldt
\end{aligned}
$$
where $ldt$ is a linear combination of terms of the new degree 
$\leq -2$. Then, the coefficients of linear terms 
$\Omega_{t_s}(a_1,\cdots,a_{n+1})=0$ for all $s$. Since 
$\Omega$ is an i.s. potential, $a_i=0$ 
for all $i$. Thus, $t_i=y_i$ is homogeneous of degree $-1$ in 
the new grading. Since $B$ is generated by $t_i$, the new 
grading agrees with the given grading. 
\end{proof}

\subsection{Rigidity of filtration}
\label{xxsec5.4}

In this subsection, we prove the following.

\begin{theorem}
\label{xxthm5.6}
Let $\Omega$ be as in \eqref{E5.0.1} where $d_0\geq 2$. Let 
$P_{\Omega-\xi}$ be defined as in Construction 
\ref{xxcon0.5} where $\xi\neq 0$. Then $P_{\Omega-\xi}$ 
has a unique filtration ${\mathbb F}$ such that 
$\gr_{\mathbb F} (P_{\Omega-\xi})$ is a connected 
graded Nambu-Poisson $d_0$-graded domain. 
\end{theorem}

\begin{proof} Let $B$ be $P_{\Omega-\xi}$. Let ${\mathbb F}$
be the original filtration by the construction and let 
$\nu^{c}$ be the associated $d_0$-valuation given in Lemma 
\ref{xxlem5.1}(2).

Suppose $B$ has a new filtration, say ${\mathbb F}'$, such that 
$\gr_{{\mathbb F}'} B$ is a Nambu-Poisson $d_0$-graded domain. 
Let $\nu'$ be the associated $d_0$-valuation on $B$. In this case 
$\nu'(f)<0$ for some $f\in B$. By Lemma \ref{xxlem5.1}(1,3), 
$\nu'=\nu^{c}$. This implies that ${\mathbb F}'={\mathbb F}$ 
as required.
\end{proof}

\subsection{Automorphism group of $n$-Lie Poisson polynomial rings}
\label{xxsec5.5}
Here, we give an example of computing the automorphism group 
of some $n$-Lie Poisson polynomial rings.

\begin{theorem}
\label{xxthm5.7} 
Let $\Omega$ be as in \eqref{E5.0.1} where $d_0\geq 2$. Let 
$A_{\Omega}$ and $P_{\Omega}$ be as in Construction 
\ref{xxcon0.5}. Then 
$$\Aut_{Poi}(Q(A_{\Omega}))=
\Aut_{Poi}(A_{\Omega})=\Aut_{Poi}(P_{\Omega}).$$
Further, every Nambu-Poisson algebra automorphism of $A_{\Omega}$ 
is Adams graded and $\Aut_{Poi}(A_{\Omega})$ is a finite subgroup 
of $GL_{n+1}(\Bbbk)$.
\end{theorem}

\begin{proof} By localization, $\Aut_{Poi}(A_{\Omega}) \subseteq 
\Aut_{Poi}(Q(A_{\Omega}))$. Since ${^1\Gamma_{0}} (Q(A_{\Omega}))
=A_{\Omega}$ [Theorem \ref{xxthm3.10}], $\Aut_{Poi}(Q(A_{\Omega}))
\subseteq \Aut_{Poi}(A_{\Omega})$. Hence, we proved the first
equation. 

Next, we show that every Nambu-Poisson automorphism of $A_{\Omega}$ 
is Adams graded and $Aut_{Poi}(A_{\Omega}) =\Aut_{Poi}(P_{\Omega})$.
Let $\phi$ be a Nambu-Poisson algebra automorphism of 
$A_{\Omega}$. Since the Poisson center of $A_{\Omega}$ is 
$\Bbbk[\Omega]$ \cite[Proposition 1]{UZ}, $\phi(\Omega)=
a\Omega+b$ for some $a\in \Bbbk^{\times}$ and $b\in \Bbbk$. 
Since $A_{\Omega}/(a\Omega+b)$ has finite global dimension 
if and only if $b\neq 0$ [Lemma \ref{xxlem3.3}(2)], we obtain
that $A_{\Omega}/(a\Omega+b)\cong A_{\Omega}/(\Omega)$ if and 
only if $b=0$. Thus $\phi(\Omega)=a\Omega$, and consequently, 
$\phi$ induces a  Nambu-Poisson algebra automorphism of 
$P_{\Omega}$. For each Nambu-Poisson algebra automorphism $\phi$ 
of $A_{\Omega}$, let $\phi'$ denote the induced automorphism of 
$P_{\Omega}$, which is Adams graded by Proposition 
\ref{xxpro5.2}(1). By Proposition \ref{xxpro5.2}(2), we can 
lift $\phi'$ to a unique, Adams graded, Nambu-Poisson algebra 
automorphism of $A_{\Omega}$, denoted by $\sigma$. Clearly 
$\sigma'=\phi'$. Let $\varphi= \phi\circ \sigma^{-1}$. Then 
$\varphi'=Id_{P_{\Omega}}$. Since $\sigma$ preserves the Adams 
grading, it suffices to show that $\varphi$ is the identity. Since 
$\varphi'$ is the identity, we have
\begin{align}
\notag%\label{E5.7.1}\tag{E5.7.1}
\varphi(t_i)&=t_i+\Omega f_i, \quad {\text{for all $i$}}
\end{align}
where $f_i\in A_{\Omega}$. Since $\Omega$ does not
have a linear term, a computation shows that 
$\varphi(\Omega)=\Omega+ \Omega \Phi(f_i)$
where $\Phi(f_i)\in (A_{\Omega})_{\geq 1}$.
Since $\varphi$ preserves the Poisson center $\Bbbk[\Omega]$ 
of $A_{\Omega}$ \cite[Proposition 1]{UZ}, 
$\varphi(\Omega)=\Omega+\Omega \Phi$ where 
$\Phi\in \Bbbk[\Omega]$. Since $\varphi$ is an 
algebra automorphism of $\Bbbk[\Omega]$, $\Phi=0$ and
$\varphi(\Omega)=\Omega$.

Let $(x_1,\cdots,x_{n+1})$ be a new $\Bbbk$-linear basis of
$V$ such that 
$$\Omega=x_{n+1}^{n+1+d_0}+{\text{linear combination of other
monomials in degree $n+1+d_0$}}.$$
Let ${\mathbb B}:=\{1, x_1,\cdots,x_{n+1}\} \cup \{b_s\}$ be 
a fixed $\Bbbk$-linear basis of $P_{\Omega}$ consisting of 
monomial elements $x_1^{s_1}\cdots x_{n+1}^{s_{n+1}}$
where $s_i\geq 0$ for all $i=1,\cdots,n$ and $0\leq s_{n+1}
\leq n+d_0$. We consider ${\mathbb B}$ as a fixed subset of 
monomial elements in $A_{\Omega}$ by a lifting. By the 
induction on the degree of elements, every element $f$ in 
$A_{\Omega}$ is of the form 
$$f=1 f^{1}(\Omega)+\sum_{i=1}^{n+1} x_i f^{x_i}(\Omega)
+\sum_{b_s} b_s f^{b_s}(\Omega)$$ 
where each $f^{\ast}(\Omega)$ is in $\Bbbk[\Omega]$. We now fix
one $k$ and write
$$\varphi(x_k)=1f^{1}(\Omega)+\sum_{i=1}^{n+1} x_i f^{x_i}(\Omega)
+\sum_{b_s} b_s f^{b_s}(\Omega).$$

For each $\xi\neq 0$, let $\pi_{\xi}$ be the quotient map from 
$A_{\Omega} \to A_{\Omega}/(\Omega-\xi)=:P_{\Omega-\xi}$. It is 
clear that the image of ${\mathbb B}$ is a $\Bbbk$-linear basis 
of $P_{\Omega-\xi}$. For simplicity, we continue to use 
$b_s$ etc for basis elements in $P_{\Omega-\xi}$.

Since $\varphi(\Omega)=\Omega$, we have $\varphi(\Omega-\xi)=\Omega-\xi$.
Let $\varphi'_{\xi}$ be the induced automorphism of 
$A_{\Omega-\xi}$ modulo $\Omega-\xi$. Then $\varphi'_{\xi}$ is a 
Nambu-Poisson algebra automorphism of $P_{\Omega-\xi}$ and
$$\varphi'_{\xi}(x_k)=1f^{1}(\xi)+\sum_{i=1}^{n+1} x_i f^{x_i}(\xi)
+\sum_{b_s} b_s f^{b_s}(\xi).$$
By Proposition \ref{xxpro5.2}(3), $\varphi'_{\xi}$ is linear and
$\varphi(V)\subseteq V$. Thus $f^1(\xi)=0$ and $f^{b_s}(\xi)=0$ 
for all $\xi\neq 0$. Hence $f^1(\Omega)=0$ and $f^{b_s}(\Omega)=0$, 
consequently, 
$\varphi(x_k)=\sum_{s=1}^{n+1} x_s f^{x_s}(\Omega)$.
Since $\varphi'(x_k)=x_k$, we have $f^{x_k}(\Omega)-1\in
\Omega\Bbbk[\Omega]$ and $f^{x_j}(\Omega)\in \Omega \Bbbk[\Omega]$
for all $j\neq k$. So we have $a_{ij}(\Omega)\in \Bbbk[\Omega]$ 
such that
$$\varphi(x_k)=x_k+\sum_{s=1}^{n+1} x_s \Omega a_{ks}(\Omega).$$
By the next lemma \ref{xxlem5.8}, $a_{ij}(\Omega)=0$ for all
$i,j$ whence $\varphi(x_k)=x_k$ for each $k=1,\cdots, n+1$ as 
required. By definition, $\varphi=\phi\circ \sigma^{-1}$. 
Hence $\phi=\sigma$. Therefore every Nambu-Poisson automorphism 
of $A_{\Omega}$ is Adams graded and $\Aut_{Poi}(A_{\Omega})
=\Aut_{Poi}(P_{\Omega})$.

Since $\phi$ is graded, there is a natural injective map 
$\phi\to \phi\mid_{V}\in GL(V)=GL_{n+1}(\Bbbk)$. Via this map 
we may consider $G:=\Aut_{Poi}(A_{\Omega})$ as a subgroup of 
$GL_{n+1}(\Bbbk)$. It remains to show that $G$ is finite. After 
replacing $\Bbbk$ with its algebraic closure, we may assume 
that $\Bbbk$ is algebraically closed. Let $B$ be the graded 
commutative algebra $P_{\Omega}$ which is $A/(\Omega)$ and let 
$X$ be the corresponding smooth hypersurface $\Proj B$. By 
\cite[p.347]{MM}, $\Aut(X)$ is finite. There is a natural 
group homomorphism $\pi$ from $\Aut_{gr.alg}(B)\to \Aut(X)$ 
whose kernel consists of automorphisms of the form 
$\eta_{\xi}: B\to B, b\mapsto \xi^{\deg b} b$ for some 
$\xi\in \Bbbk$. Thus, we have an exact sequence
\begin{equation}
\label{E5.7.1}\tag{E5.7.1}
1\to \ker \pi\cap G \to G \to \Aut(X).
\end{equation}
One can check that $\ker \pi\cap G=\{\eta_{\xi}
\mid \xi^{d_0}=1\}$. Thus $|\ker \pi\cap G|= d_0$. Now
\begin{equation}
\notag%\label{E5.7.2}\tag{E5.7.2}
|G|=|\ker \pi\cap G| |\Aut(X)|
\leq d_0|\Aut(X)|<\infty.
\end{equation}
\end{proof}

By \cite[Theorem 5]{MM}, $\Aut(X)$ is trivial for a generic 
$\Omega$. In this case, by \eqref{E5.7.1},
$\Aut_{Poi}(A_{\Omega})=\ker \pi\cap G=C_{d_0}$. 

\begin{lemma}
\label{xxlem5.8}
Suppose $\Omega$ is an i.s. potential of degree $n+1+d_0$
with $d_0\geq 2$. 
\begin{enumerate}
\item[(1)]
Let $f_1,\cdots,f_{n+1}$ be elements in $V:=\sum_{s=1}^{n+1} 
\Bbbk t_s$ such that $\sum_{s=1}^{n+1} f_s \Omega_{t_s}=0$. 
Then $f_i=0$ for all $i=1,\cdots,n+1$.
\item[(2)]
Let $(x_1,\cdots,x_{n+1})$ be a basis of $V$ and let 
$\varphi$ be a Nambu-Poisson automorphism of $A_{\Omega}$ 
such that $\varphi(\Omega)=\Omega$ and 
$$
\varphi(x_i)=x_i+ \sum_{s=1}^{n+1}x_s\Omega a_{is}(\Omega)
$$
for some polynomials $a_{ij}(t)\in \kk[t]$. Then 
$\varphi=Id$.
\end{enumerate}
\end{lemma}

\begin{proof}
(1) Let $A=\Bbbk[t_1,\cdots, t_{n+1}]$ with Adams degree 
$\deg(t_s)=1$ for all $s$. Let $M$ be the free module 
$\bigoplus_{s=1}^{n+1} A \Omega_{t_s}$ which is isomorphic 
to $A^{\oplus (n+1)}[-(n+d_0)]$ as graded $A$-modules. The 
Koszul complex associated with the sequence 
$(\Omega_{t_1},\cdots, \Omega_{t_{n+1}})$ is
\begin{equation}
\label{E5.8.1}\tag{E5.8.1}
0\to \Lambda^{n+1}(M)
\to \cdots \Lambda^{2}(M)\xrightarrow{\partial_2} M
\xrightarrow{\partial_1}
A\to A/(\Omega_s)
\to 0.
\end{equation}
Since each $\Omega_{t_i}$ is homogeneous, the above 
is a complex of graded $A$-modules. Since $\Omega$ is
an i.s. potential, $(\Omega_{t_1},\cdots, \Omega_{t_{n+1}})$
is a regular sequence. Consequently, complex \eqref{E5.8.1} is exact.

Suppose $\sum_{s=1}^{n+1} f_s \Omega_{t_s}=0$ for some
elements $f_i\in V$ that are not all zero. Then 
$m:=(f_1,\cdots,f_{n+1})$, considered as an element 
in $M$, is in the kernel of the map $\partial_1$ 
in the above exact Koszul complex. Then $m$ is in the
image of $\partial_2$. But $\deg m=n+d_0+1$ in $M$ and 
the lowest degree in $\Lambda^{2}(M)$ is $2(n+d_0)$
that is strictly larger than $n+d_0+1$. So $m$ cannot 
be in the image of $\partial_2$. Since \eqref{E5.8.1}
is exact, this forces that $m=0$ in $M$. Therefore
$f_s=0$ for all $s$.

(2) After a linear change of basis, we have
$$
\varphi(t_i)=t_i+ \sum_{s=1}^{n+1}x_s\Omega w_{is}(\Omega),
\quad i=1,\cdots, n+1
$$
for some polynomials $w_{is}(t)$. We need to prove that 
$w_{ij}(t)=0$ for all $i,j$. Suppose to the contrary that 
some $w_{ij}\neq 0$ for some $(i,j)$. Then we can 
write
$$\begin{aligned}
\varphi(t_i)&=t_i+ f_i \Omega^{m}+hdt
\quad i=1,\cdots, n+1
\end{aligned}
$$
where $f_1,\cdots,f_{n+1}\in V$ are not all zero, $m\geq 1$, 
and $hdt$ stands for the linear combination of higher Adams 
degree terms. Through Taylor expansion,
$$\begin{aligned}
\Omega&=\phi(\Omega)\\
&=\Omega(t_1+f_1 \Omega^m+ hdt, \cdots, t_{n+1}+f_{n+1} 
\Omega^m+hdt)\\
&=\Omega(t_1,\cdots,t_{n+1})+(\sum_{s=1}^{n+1} f_s
\Omega_{t_s}(t_1,\cdots,t_{n+1}))\Omega^m+hdt
\end{aligned}
$$
which implies that $\sum_{s=1}^{n+1} f_s\Omega_{t_s}=0$. 
By part (1), $f_i=0$ for all $i$, yielding a contradiction.
\end{proof}

\subsection{Embedding problem}
\label{xxsec5.6}
Recall that the $q$-skew Nambu-Poisson field $N_{q}$ 
is defined in Example \ref{xxex1.3}(1). Also, see related 
Nambu-Poisson fields in Lemma \ref{xxlem4.13}.

\begin{corollary}
\label{xxcor5.9}
Let $N_{q}$ be the $q$-skew Nambu-Poisson field and $K$ be
another Nambu-Poisson field. 
\begin{enumerate}
\item[(1)]
If ${\mathcal V}_{-1}(K)\neq \emptyset$, then there is no 
Poisson algebra morphism from $N_{q}$ to $K$.
\item[(2)]
If $K$ is either the Weyl Poisson field $N_{Weyl}$ or 
Nambu-Poisson field $Q(P_{\Omega})$ where $\deg 
\Omega\neq \sum_{s=1}^{n+1} \deg(t_s)$, then there is 
no Nambu-Poisson algebra morphism from $N_{q}$ to $K$.
\end{enumerate}
\end{corollary}

\begin{proof}
(1) Suppose there is an embedding $N_q\to K$. Let $\nu$ be an 
element in ${\mathcal V}_{-1}(K)$. Consider $\nu$ as a 
$0$-valuation on $K$, it is a classical valuation. This means 
the Nambu-Poisson structure on $\gr_{\nu} K$ is zero. 
Therefore, the restriction of $\nu$ on $N_q$ is also classical. 
This contradicts Lemma \ref{xxlem4.13}(1).

(2) If $K=N_{Weyl}$, it is easy to check that 
${\mathcal V}_{-1}(N_{Weyl})\neq \emptyset$. If $K=Q(P_{\Omega})$, 
then ${\mathcal V}_{-1}(K)\neq \emptyset$ by Lemmas 
\ref{xxlem1.8}(2) and \ref{xxlem2.12}(1). So, the assertion 
follows from part (1). 
\end{proof}

\begin{corollary}
\label{xxcor5.10}
Let $\Omega_0$ be an i.s. potential. Let $N$ be the 
Nambu-Poisson field $Q(P_{\Omega_0-1})$ and $K$ be another Nambu-Poisson field of the same GK-dimension as $N_q$. 
\begin{enumerate}
\item[(1)]
If ${\mathcal V}_{-1}(K)\neq \emptyset$, then there is no 
Poisson algebra morphism from $N$ to $K$.
\item[(2)]
If $K$ is either $N_{Weyl}$ or $Q(P_{\Omega})$ where 
$\deg \Omega\neq \sum_{i=1}^{n+1}\deg t_i$, then there is  no 
Nambu-Poisson algebra morphism from $N$ to $K$.
\end{enumerate}
\end{corollary}

\begin{proof}
(1) Suppose there is an embedding $N\to K$. Let $\nu$ be an 
element in ${\mathcal V}_{-1}(K)$ and let $\mu$ be the 
restriction of $\nu$ on $N$. Consider $\mu$ as a 
$0$-valuation, it is classical. By Lemma \ref{xxlem1.8}(1),
$\mu$ is nontrivial. By Theorem \ref{xxthm3.9}(2),
$Q(\gr_{\mu} N)\cong Q(P_{\Omega_0})$, and whence 
$\mu$ is nonclassical, yielding a contradiction.

(2) The proof is similar to the proof of Corollary 
\ref{xxcor5.9}(2).
\end{proof}

\begin{corollary}
\label{xxcor5.11}
Suppose $\Omega$ is an i.s. potential of degree 
$n+1+d_0$ for $d_0\geq 2$. Let $N$ be the Nambu-Poisson 
fraction field of $P_{\Omega-\xi}$ where $\xi\in 
\Bbbk^{\times}$ and $K$ be another Nambu-Poisson 
field of the same GK-dimension as $N_q$.
\begin{enumerate}
\item[(1)]
If ${\mathcal V}_{0}(K)\neq \emptyset$, then there is no 
Nambu-Poisson algebra morphism from $N$ to $K$.
\item[(2)]
If $K$ is either $N_{q}$, or $N_{Weyl}$, or $Q(P_{\Omega_0})$ 
where $\Omega_0$ is homogeneous of degree $n+1$, then there 
is no Nambu-Poisson algebra morphism from $N$ to $K$.
\end{enumerate}
\end{corollary}

\begin{proof}
(1) Suppose there is an embedding $N\to K$. Let $\nu$ be an 
element in ${\mathcal V}_{0}(K)$ and let $\mu$ be the 
restriction of $\nu$ on $N$. By Lemma \ref{xxlem1.8}(1),
$\mu$ is nontrivial. This contradicts Lemma \ref{xxlem3.5}(3).

(2) In all cases of $K$, one can easily see that 
${\mathcal V}_{0}(K)\neq \emptyset$. The assertion 
follows from part (1).
\end{proof}

Without using valuations, it isn't easy to show 
the part (2) of the above Corollaries. 

\section{A realization lemma}
\label{xxsec6}

In this section, we want to connect the Nambu-Poisson automorphisms
of $A_{\Omega}$ with the ordinary algebra automorphisms of 
$A_{\Omega}$. This helps us to understand certain automorphisms 
of $A_{\Omega}$ with fixed quasi-axis. 

Let $S$ be a subset of $A:=\Bbbk[t_1,\cdots,t_{n+1}]$. Let
$$\Aut(A\mid S)
:=\{\sigma\in \Aut_{alg}(A)\mid \sigma(f)\in S,
{\text{ for every $f\in S$}}\}.$$
For example, let $\Omega$ be any element in $A$
which may not be homogeneous. Then
\begin{equation}
\label{E6.0.1}\tag{E6.0.1}
\Aut(A\mid \{\Omega\})
=\{\sigma\in \Aut_{alg}(A)\mid \sigma(\Omega)=
\Omega\}
\end{equation}
which is called the {\it automorphism group of $A$ with fixed 
quasi-axis $\Omega$}. Similarly, we have
$$\Aut(A\mid \Bbbk^{\times}\Omega+\Bbbk)
=\{\sigma\in \Aut_{alg}(A)\mid \sigma(\Omega)
=a\Omega+b, {\text{ for some $a\in\Bbbk^{\times}$
and $b\in \Bbbk$}}\}.$$

Let $A$ be a Nambu-Poisson algebra with bracket 
$\pi:=\{-,\cdots,-\}: A^{\otimes n}\to A$ and let $e$ be 
in $\Bbbk^{\times}$. Then, we can define a new 
Nambu-Poisson structure on $A$ by
$$\{f_1,\cdots, f_n\}:= e\{f_1,\cdots,f_n\} \quad 
{\text{for all $f_i\in A$.}}$$
This new Nambu-Poisson algebra is denoted by $A(e)$. 

\begin{definition}
\label{xxdef6.1}
Let $A$ and $B$ be two Nambu-Poisson algebras. A $\Bbbk$-algebra
morphism $\phi: A\to B$ is called an {\it $e$-morphism}
if it is a Nambu-Poisson algebra morphism from $A\to B(e)$ for
some $e\in \Bbbk^{\times}$. If $e$ is not 
specified, we say $\phi$ is an {\it $\epsilon$-morphism}.

Similarly, one can define $e$-($\epsilon$-)versions of Nambu-Poisson
endomorphism, Nambu-Poisson isomorphism, Nambu-Poisson automorphism, etc.
\end{definition}

Let $A$ be a Nambu-Poisson algebra. Then we let
$$\Aut_{\Bbbk^{\times} Poi}(A)=\{\phi\in \Aut_{alg}(A)\mid 
{\text{ $\phi$ is an $\epsilon$-automorphism of $A$}}\}.$$
The next lemma is easy.

\begin{lemma}
\label{xxlem6.2}
Let $A$ be a Nambu-Poisson algebra and $e\in \Bbbk^{\times}$.
The following hold.
\begin{enumerate}
\item[(1)]
$\nu$ is a $w$-valuation on $A$ if and only if 
it is a $w$-valuation on $A(e)$
\item[(2)]
$\Aut_{Poi}(A)=\Aut_{Poi}(A(e))$.
\item[(3)]
$\Aut_{\Bbbk^{\times} Poi}(A)=\Aut_{\Bbbk^{\times} Poi}(A(e))$.
\end{enumerate}
\end{lemma}

By the above lemma, we can replace Nambu-Poisson morphisms 
by $\epsilon$-morphisms when $w$-valuations are 
concerned. Below is the main result -- a realization lemma.

\begin{theorem}[Realization Lemma]
\label{xxthm6.3} 
Let $\Omega\in A:=\Bbbk^{[n+1]}$
be an irreducible element of degree at least 2. Then 
$$\Aut(A\mid \Bbbk^{\times}\Omega+\Bbbk)
=\Aut_{\Bbbk^{\times} Poi}(A_{\Omega}).$$
\end{theorem}

\begin{proof} Let $\phi\in \Aut_{\Bbbk^{\times} Poi}(A_{\Omega})$.
By \cite[Proposition 1]{UZ}, the Poisson center of $A_{\Omega}$ 
is $\Bbbk[\Omega]$. It is easy to check that $\phi$ preserves the
Poisson center. So $\phi$ induces an automorphism of
$\Bbbk[\Omega]$. Thus $\phi(\Omega)=a\Omega+b$ for some 
$a\in \Bbbk^{\times}$ and $b\in \Bbbk$. This means that 
$\phi\in \Aut(A\mid \Bbbk^{\times}\Omega+\Bbbk)$. 

Conversely, let $\phi\in \Aut(A\mid \Bbbk^{\times}\Omega+\Bbbk)$.
So $\phi(\Omega)=a\Omega+b$ for some $a\in \Bbbk^{\times}$ and 
$b\in \Bbbk$. Let $\{-,\cdots,-\}$ be the Nambu-Poisson bracket 
defined in Construction \ref{xxcon0.5}. Let
$\frac{\partial (f_1,\cdots, f_{n+1})}
{\partial (t_1,\cdots, t_{n+1})}$ be the Jacobian determinant
defined as in \eqref{E0.5.1}. Now we compute
$$\begin{aligned}
\phi(\{f_1,\cdots,f_n\})
&=\phi\left( \frac{\partial (f_1,\cdots, f_n, \Omega)}
{\partial (t_1,\cdots, t_{n+1})}\right)\\
&=\frac{\partial (\phi(f_1),\cdots,\phi(f_n), \phi(\Omega))}
{\partial (\phi(t_1),\cdots, \phi(t_{n+1}))}\\
&=\frac{\partial (\phi(f_1),\cdots, \phi(f_n), \phi(\Omega))}
{\partial (t_1,\cdots, t_{n+1})}
\cdot
\left(\frac{\partial (\phi(t_1),\cdots,\phi(t_{n+1}))}
{\partial (t_1,\cdots, t_{n+1})}\right)^{-1}\\
&=\frac{\partial (\phi(f_1),\cdots, \phi(f_n), a\Omega))}
{\partial (t_1,\cdots, t_{n+1})} J^{-1}\\
&=J^{-1}a \{\phi(f_1),\cdots, \phi(f_n)\}
\end{aligned}
$$
where $J:=\frac{\partial (\phi(t_1),\cdots, \phi(t_{n+1}))}
{\partial (t_1,\cdots, t_{n+1})}\in \Bbbk$ as 
$\phi$ is an algebra automorphism of $A$.
Therefore $\phi\in \Aut_{\Bbbk^{\times} Poi}(A_{\Omega})$.
\end{proof}

The above theorem is helpful for computing 
$\Aut(A\mid \Bbbk^{\times}\Omega+\Bbbk)$ for certain $\Omega$.

\begin{proposition}
\label{xxpro6.4}
Let $\Omega$ be an i.s. potential of degree $n+1+d_0$ with 
$d_0 \geq 2$. Let $A_{\Omega}$ and $P_{\Omega}$ be as in 
Construction \ref{xxcon0.5}. 
\begin{enumerate}
\item[(1)]
$$\Aut_{\Bbbk^{\times} Poi}(Q(A_{\Omega}))=
\Aut_{\Bbbk^{\times} Poi}(A_{\Omega})=\Aut_{\Bbbk^{\times} Poi}
(P_{\Omega}).$$
Further, every Nambu-Poisson $e$-automorphism of $A_{\Omega}$ 
is Adams graded and $\Aut_{\Bbbk^{\times} Poi}(A_{\Omega})$ is a 
subgroup of $GL_{n+1}(\Bbbk)$.
\item[(2)]
Every algebra automorphism in $\Aut(A\mid \Bbbk^{\times}\Omega+\Bbbk)$
is graded. Consequently, $\Aut(A\mid \Bbbk^{\times}\Omega+\Bbbk)$ is
a subgroup of $GL_{n+1}(\Bbbk)$.
\end{enumerate}
\end{proposition}

\begin{proof}
(1) The proof of Theorem \ref{xxthm5.7} works in the 
$e$-automorphism setting.

(2) This follows from (1) and Theorem \ref{xxthm6.3}.
\end{proof}

For the rest of this section, we consider one potential 
$\Omega$ explicitly. Let
\begin{equation}
\label{E6.4.1}\tag{E6.4.1}
\Omega=\sum_{s=1}^{n+1} t_s^{n+1+d_0}
\end{equation}
where $d\geq 2$. Let 
\begin{equation}
\label{E6.4.2}\tag{E6.4.2}
G_0:=\{(a_1,\cdots,a_{n+1})\in (\Bbbk^{\times})^{n+1} \mid
\prod_{s=1}^{n+1} a_s =a_i^{n+d_0+1}, \forall \; i=1,\cdots,n+1\}
\end{equation}
and
\begin{equation}
\label{E6.4.3}\tag{E6.4.3}
G_1:=\{(a_1,\cdots,a_{n+1})\in (\Bbbk^{\times})^{n+1} \mid
\prod_{s=1}^{n+1} a_s =a_i^{n+d_0+1}=1, \forall \; i=1,\cdots,n+1\}.
\end{equation}
Note that both $G_0$ and $G_1$ are finite groups.

\begin{lemma}
\label{xxlem6.5}
Let $\Omega$ be as in \eqref{E6.4.1} where $d_0\geq 2$.
Let $P_{\Omega}$ and $P_{\Omega-\xi}$ be defined as in 
Construction \ref{xxcon0.5} where $\xi\neq 0$.
\begin{enumerate}
\item[(1)]
There is a short exact sequence of groups
$$\{1\}\to G_0\to \Aut_{Poi}(P_{\Omega})\to S_{n+1}\to \{1\},$$
and $\Aut_{Poi}(P_{\Omega})\cong S_{n+1}\ltimes G_0$ if and only 
if $d_0$ is odd.
\item[(2)]
If $n+d_0$ is even, there is a short exact sequence of groups
$$\{1\}\to G_1\to \Aut_{Poi}(P_{\Omega-\xi})\to A_{n+1}\to \{1\},$$
and $\Aut_{Poi}(P_{\Omega-\xi})\cong A_{n+1}\ltimes G_1$.
\item[(3)]
If $n+d_0$ is odd, there is a short exact sequence of groups
$$\{1\}\to G_1\to \Aut_{Poi}(P_{\Omega-\xi})\to S_{n+1}\to \{1\},$$
and $\Aut_{Poi}(P_{\Omega-\xi})\cong S_{n+1}\ltimes G_1$.
\end{enumerate}
\end{lemma}

\begin{proof} (1) Let $P=P_{\Omega}$. We divide the proof 
into several steps. By Proposition \ref{xxpro5.2}(2), every 
automorphism lifts to a graded automorphism of $A_{\Omega}$.
Let $V=\sum_{i=1}^{n+1} \Bbbk x_i$. For every element 
$f\in V$, let $W_f$ denote the $\Bbbk$-linear span of 
elements in $\{f,V,\cdots, V\}$ (considered as a subspace 
of $\Bbbk^{[n+1]}$). 

Step 1: If $f\in \bigcup_{i=1}^{n+1}\Bbbk x_i$, we claim that 
every element in $W_f$ does not have an isolated singularity at 
zero. Without loss of generality, we may assume that $f=x_1$ and 
$g_j= \sum_{i=1}^{n+1} c_{ji} x_i$ for $j=2,\cdots, n$. Then 
$$\begin{aligned}
\{f,g_2,\cdots, g_{n}\}&
=\sum c_{\ast,\ast} \{x_1,x_{i_2},\cdots,x_{i_n}\}\\
&\in \sum_{t=2}^{n+1} \Bbbk x_t^{n+d_0}.
\end{aligned}
$$ 
This implies that $W_f=\sum_{t=2}^{n+1} \Bbbk x_t^{n+d_0}$.
It is clear that every element in $W_f$ does not have 
isolated singularity at zero.

Step 2: If $f\not\in \bigcup_{i=1}^{n+1}\Bbbk x_i$, we claim 
that there is an element in $W_f$ such that it has an isolated 
singularity at zero. Write $f=x_1+\sum_{s=2}^{n+1} e_{s} x_s$
where $e_s\in \Bbbk$ and $e_2\neq 0$. Let $\{g_2,\cdots, g_n\}$ be
$\{x_2,\cdots, \widehat{x_i},\cdots, x_{n+1}\}$, then, for every
$i\geq 2$,
$$h_i:=\pm \{f,g_2,\cdots,g_n\}
=\pm (n+1+d) x_i^{n+d_0} + e_i (n+1+d_0) x_1^{n+d_0}.$$
By definition, $h_i\in W_f$. Let $w_i\in \Bbbk^{\times}$ 
such that $\sum_{i=2}^{n+1} e_i w_i\neq 0$, then 
$$\sum_{i=2}^{n+1} w_i h_i=
(n+1+d_0)(\sum_{i=2}^{n+1} \pm w_i x_i^{n+d_0})+(n+1+d_0)
(\sum_{i=2}^{n+1} e_i w_i) x_1^{n+d_0}$$
which has an isolated singularity at zero. Note that 
$\sum_{i=2}^{n+1} w_i h_i\in W_f$, so the claim is proved.

Step 3: Let $\phi$ be a Nambu-Poisson algebra automorphism 
of $P$. By Steps 1 and 2, $\phi$ preserves 
$\bigcup_{i=1}^{n+1} \Bbbk x_i$. Hence there are 
$\sigma\in S_{n+1}$ and $(a_1,\cdots,a_{n+1})\in 
(\Bbbk^{\times})^{n+1}$ such that
\begin{equation}
\label{E6.5.1}\tag{E6.5.1}
\phi(x_i)=a_i x_{\sigma(i)}, \quad 
{\text{for all $i=1,\cdots,n+1$.}}
\end{equation}
Applying $\phi$ to the Poisson bracket \eqref{E3.1.2}, 
we get
\[
\{a_1x_{\sigma(1)},\ldots,\widehat{a_ix_{\sigma(i)}},
\ldots,a_{n+1}x_{\sigma(n+1)}\}=(n+1+d_0)a_i^{n+d_0}x_{\sigma(i)}
\]
for all $1\le i\le n+1$. So, we obtain that 
\begin{equation}
\label{E6.5.2}\tag{E6.5.2}
\prod_{s\neq i} a_s= {\text{sgn}}(\sigma) a_i^{n+d_0}, 
\quad {\text{for all $i=1,\cdots,n+1$.}}
\end{equation}
If $\phi$ preserves $\Bbbk x_i$ individually for each $i$, 
then $\sigma$ is an identity. In this case \eqref{E6.5.2} 
is equivalent to $(a_1,\cdots,a_{n+1})\in G_0$. By an 
elementary computation and the hypothesis on $\Bbbk$ for 
every $\sigma \in S_{n+1}$, \eqref{E6.5.2} has a solution 
$(a_1,\cdots,a_{n+1})$. Indeed, one can show $\prod_{s=1}^{n+1}a_s^{d_0}={\rm sgn}(\sigma)^{n+1}$ and $a_i^{n+1+d_0}={\rm sng}(\sigma)^{(n+1-d_0)/d_0}$ for all $i$. Combining these facts, $G_0$ is a normal subgroup
of $\Aut_{Poi}(P_{\Omega})$ and $\Aut_{Poi}(P_{\Omega})/G_0
\cong S_{n+1}$. Therefore, the main assertion in part (1) 
is proved.

Step 4: Suppose $d$ is odd. For every $\sigma\in S_{n+1}$, 
let $\sigma'$ be the automorphism of $P_{\Omega}$ defined
by $\sigma'(x_i)={\rm{sgn}}(\sigma)x_{\sigma(i)}$ (or $a_i=
{\rm{sgn}}(\sigma)$ for all $i$). One can 
check by \eqref{E6.5.2} that $\sigma'$ is indeed a 
Poisson algebra automorphism of $P_{\Omega}$. Moreover, 
the subgroup of $\Aut_{Poi}(P_{\Omega})$ generated by 
$\{\sigma'\,|\, \sigma\in S_{n+1}\}$ is isomorphic to $S_{n+1}$. 
Therefore $\Aut_{Poi}(P_{\Omega})\cong S_{n+1}\ltimes G_0$.

Conversely, suppose $\Aut_{Poi}(P_{\Omega})\cong S_{n+1}
\ltimes G_0$. Thus it means that $S_{n+1}$ is a subgroup
of $\Aut_{Poi}(P_{\Omega})$. Therefore, there are 
$\varphi_1,\cdots,\varphi_n\in \Aut_{Poi}(P_{\Omega})$ 
satisfying the Coxeter relations in $S_{n+1}$ such that 
\[
\varphi_i^2=1,\ (\varphi_i\varphi_{i+1})^3=1,\ 
(\varphi_i\varphi_j)^2=1,\ \forall |i-j|>1.
\]
Without loss of generality, we can let $\varphi_i$ correspond 
to the transposition $(i,i+1)$ in $S_{n+1}$. So for $\varphi_1$, 
there are $a_1,\cdots,a_{n+1}\in \kk^\times$ such that
$\varphi_1(x_1)=a_1 x_2$ and $\varphi_1(x_2)=a_2 x_1$ and 
$\varphi_1(x_i)=a_i x_i$ for all $i\geq 3$. Since 
$\varphi^2_1=Id$, we have $a_1a_2=1=a_i^2$ for all $i\geq 3$. 

Suppose $d_0$ is even. By using the similar 
expression for $\varphi_i$ as for $\varphi_1$ and relation 
$(\varphi_1\varphi_i)^2=1$ for $i\ge 3$, we obtain 
$a_3=a_4=\cdots=a_{n+1}$. By \eqref{E6.5.2}, we get 
$a_1a_2\widehat{a_3}a_4\cdots a_{n+1}=-a_3^{n+d_0}$. Again since 
$a_1a_2=1$ and $a:=a_3=a_4=\cdots=a_{n+1}=\pm 1$, we have 
$a^{d_0}=-1$. This yields a contradiction. 

(2) Let $\phi$ be an $n$-Lie Poisson algebra automorphism 
of $P_{\Omega-\xi}$ for $\xi\neq 0$. By Proposition 
\ref{xxpro5.2}(3,4), we can write $\phi$ in the form of 
\eqref{E6.5.1} for $(a_1,\cdots,a_{n+1})\in (\kk^\times)^{n+1}$ 
and $\sigma\in S_{n+1}$ such that \eqref{E6.5.2} hold. Since 
$\phi(\Omega)=\xi$, we additionally have 
\[
\prod_{i=1}^{n+1}a_i={\rm sgn}(\sigma). 
\]
Thus if $\phi$ preserves $\kk x_i$ individually for each $i$, 
then $\sigma$ is identity and $(a_1,\ldots,a_{n+1})\in G_1$. 
Now by multiplying \eqref{E6.5.2} together and by using the 
above relation, we get ${\rm sgn}(\sigma)^{n+1+d_0}=1$. If $n+d_0$ is 
even, ${\rm sgn}(\sigma)=1$ and $\sigma\in A_{n+1}$. Then, all 
the rest of the proof follows again the proof of part (1).  

(3) The proof is similar to the proof of part (2), whence it 
is omitted.
\end{proof}

Using the realization lemma, we can determine the 
automorphism group $A$ with fixed quasi-axis $\Omega$.
Let
\begin{equation}
\label{E6.5.3}\tag{E6.5.3}
G_2:=\{(a_1,\cdots,a_{n+1})\in (\Bbbk^{\times})^{n+1} \mid
a_i^{n+d_0+1}=a_j^{n+d_0+1}, \forall \; 1\leq i,j\leq n+1\}.
\end{equation}
and
\begin{equation}
\label{E6.5.4}\tag{E6.5.4}
G_3:=\{(a_1,\cdots,a_{n+1})\in (\Bbbk^{\times})^{n+1} \mid
a_i^{n+d_0+1}=1, \forall \; i=1,\cdots,n+1\}.
\end{equation}
Note that $G_2$ is an infinite group, and $G_3$ is a finite 
group.

\begin{lemma}
\label{xxlem6.6}
Let $\Omega$ be as in \eqref{E6.4.1} where $d_0\geq 2$.
\begin{enumerate}
\item[(1)]
There is a short exact sequence of groups
$$\{1\}\to G_2\to 
\Aut_{\Bbbk^{\times}Poi}(P_{\Omega})\to S_{n+1}\to \{1\}.$$
\item[(2)]
There is a short exact sequence of groups
$$\{1\}\to G_3\to 
\Aut(A\mid \{\Omega\})\to S_{n+1}\to \{1\}.$$
\end{enumerate}
\end{lemma}

\begin{proof}
The proof is similar to the proof of Lemma \ref{xxlem6.5}(1)
where Proposition \ref{xxpro5.2} is replaced by 
Proposition \ref{xxpro6.4}. Details are omitted.
\end{proof}

\section{Partial differential equations in the field
$\Bbbk^{(n)}$}
\label{xxsec7}

Let $\Bbbk^{(n)}$ be the field of rational functions in 
variables $t_1,\cdots,t_n$. In this section, we use $a$ and $b$ 
to denote nonzero elements in $\Bbbk^{(n)}$. When it is 
necessary we write $a(t_1,\cdots,t_n)$ for $a$. Let 
$y_1,\cdots,y_n$ be $n$ elements in $\Bbbk^{(n)}$. The 
Jacobian determinant of $y_1,\cdots, y_n$ is defined to be 
\begin{equation}
\label{E7.0.1}\tag{E7.0.1}
\frac{\partial (y_1,\cdots,y_n)}{\partial (t_1,\cdots,t_n)}
:=
\det \begin{pmatrix} 
(y_1)_{t_1} & (y_1)_{t_2} & \cdots & (y_1)_{t_n} \\
(y_2)_{t_1} & (y_2)_{t_2} & \cdots & (y_2)_{t_n} \\
\cdots      & \cdots      & \cdots & \cdots\\
(y_n)_{t_1} & (y_n)_{t_2} & \cdots & (y_n)_{t_n} \\
\end{pmatrix}\in \Bbbk^{(n)}.
\end{equation}

In this section, we consider a class of ``separable'' 
partial differential equations (PDEs) of the form
\begin{equation}
\label{E7.0.2}\tag{E7.0.2}
\frac{\partial (y_1,\cdots,y_n)}{\partial (t_1,\cdots,t_n)}
=\frac{b(y_1,\cdots,y_2)}{a(t_1,\cdots,t_n)}
\end{equation}
where $a,b$ are nonzero elements in $\Bbbk^{(n)}$. 

\begin{definition}
Given two nonzero elements $a,b\in \Bbbk^{(n)}$. We say 
PDE \eqref{E7.0.2} has a solution if there are 
algebraically independent elements $y_1,\cdots,y_n
\in \Bbbk^{(n)}$ such that \eqref{E7.0.2} holds.
\end{definition}

There are some overlaps between this section and \cite{GZ} 
where PDEs in $\Bbbk^{(2)}$ are discussed.

As in the above definition, we only consider the solutions to 
\eqref{E7.0.2} in the rational function field $\Bbbk^{(n)}$.
Next, we define a class of Nambu-Poisson structure on the field 
$\Bbbk^{(n)}$ of rational functions. Let $a\in\Bbbk^{(n)}$ be 
a nonzero element. We use $\Bbbk^{(n,a)}$ denote the Nambu-Poisson 
algebra $\Bbbk^{(n)}$ with Nambu-Poisson structure determined by
$$\{t_1,\cdots,t_n\}_a= a.$$
If $a=1$, $\Bbbk^{(n,1)}$ is just the Weyl Nambu-Poisson field 
$N_{Weyl}$ [Example \ref{xxex1.3}(3)]. In this case, the 
Nambu-Poisson bracket is denoted by $\{-,\cdots,-\}_w$ instead of 
$\{-,\cdots,-\}_1$. If $a=q t_1\cdots t_n$, then $\Bbbk^{(n,a)}$
is isomorphic to $N_q$ [Example \ref{xxex1.3}(1)]. For any 
$y_1,\cdots, y_n\in \Bbbk^{(n)}$, the Leibniz rule in Definition 
\ref{xxdef1.2}(2) implies that
\begin{equation}
\label{E7.1.1}\tag{E7.1.1}
\{y_1,\cdots,y_n\}_a=
\frac{\partial (y_1,\cdots,y_n)}{\partial (t_1,\cdots,t_n)}
\{t_1,\cdots,t_n\}_a=
a\frac{\partial (y_1,\cdots,y_n)}{\partial (t_1,\cdots,t_n)}.
\end{equation}

The embedding property between the Nambu-Poisson fields 
$\Bbbk^{(n,b)}\to \Bbbk^{(n,a)}$ is closely related to the 
solution of \eqref{E7.0.2}.

\begin{lemma}
\label{xxlem7.2}
Let $a,b\in \Bbbk^{(n)}$. Then \eqref{E7.0.2} has a solution 
if and only if there is a Nambu-Poisson morphism from
$\Bbbk^{(n,b)}\to \Bbbk^{(n,a)}$.
\end{lemma}

\begin{proof}
``$\Longleftarrow$''
Suppose $\phi: \Bbbk^{(n,b)}\to \Bbbk^{(n,a)}$ is a 
Nambu-Poisson morphism between these two Nambu-Poisson
fields. Let $y_i=\phi(t_i)$ for all $i$. Then 
$$\begin{aligned}
b(y_1,\cdots, y_n)& =b(\phi(t_1),\cdots,\phi(t_n))=\phi(b)\\
&=\phi(\{t_1,\cdots,t_n\}_b)=\{\phi(t_1),\cdots,\phi(t_n)\}_a\\
&=\{y_1,\cdots,y_n\}_a\\
&=a(t_1,\cdots, t_n)
\frac{\partial (y_1,\cdots,y_n)}{\partial (t_1,\cdots,t_n)}.
\end{aligned}
$$
Therefore
$$\frac{\partial (y_1,\cdots,y_n)}{\partial (t_1,\cdots,t_n)}
=\frac{b(y_1,\cdots, y_n)}{a(t_1,\cdots,t_n)}.$$
Since $\phi$ is injective, $y_1,\cdots,y_n$ are
algebraically independent. Hence, \eqref{E7.0.2} has a
solution.

``$\Longrightarrow$''
Suppose $(y_1,\cdots,y_n)$ is a solution to \eqref{E7.0.2}. Let 
$\phi$ be the endomorphism of the field $\Bbbk^{(n)}$ determined 
by $\phi(t_i)=y_i$ for all $i$. The above computation shows that
$\phi(\{t_1,\cdots,t_n\}_b)=\{\phi(t_1),\cdots,\phi(t_n)\}_a$ 
assuming \eqref{E7.0.2} holds. Since $\Bbbk^{(n,b)}$ is generated 
by $t_1,\cdots,t_n$, we have 
$\phi(\{f_1,\cdots,f_n\}_b)=\{\phi(f_1),\cdots,\phi(f_n)\}_a$
for all $f_1,\cdots, f_n\in \Bbbk^{(n,b)}$. Therefore 
$\phi$ is a Nambu-Poisson morphism.  
\end{proof}

\begin{lemma}
\label{xxlem7.3} 
Let $a,b,c$ be three nonzero elements in $\Bbbk^{(n)}$. If
PDEs
$$\frac{\partial (y_1,\cdots,y_n)}{\partial (t_1,\cdots,t_n)}
=\frac{b(y_1,\cdots, y_n)}{a(t_1,\cdots,t_n)}
\quad
{\text{and}}\quad
\frac{\partial (y_1,\cdots,y_n)}{\partial (t_1,\cdots,t_n)}
=\frac{c(y_1,\cdots, y_n)}{b(t_1,\cdots,t_n)}
$$
have solutions, then 
the PDE
$$\frac{\partial (y_1,\cdots,y_n)}{\partial (t_1,\cdots,t_n)}
=\frac{c(y_1,\cdots, y_n)}{a(t_1,\cdots,t_n)}$$
has a solution.
\end{lemma}

\begin{proof} By Lemma \ref{xxlem7.2}, there are injective
Nambu-Poisson morphisms $\Bbbk^{(n,b)}\to \Bbbk^{(n,a)}$
and $\Bbbk^{(n,c)}\to \Bbbk^{(n,b)}$. Then, there is an 
Nambu-Poisson morphism $\Bbbk^{(n,c)}\to \Bbbk^{(n,a)}$.
The assertion follows by applying Lemma \ref{xxlem7.2} 
again.
\end{proof}

Applying Lemma \ref{xxlem7.2} to previous non-embedding 
results between Nambu-Poisson algebras $\kk^{(n,a)}$, we have

\begin{corollary}
\label{xxcor7.4} The following PDEs do not have a solution.
\begin{enumerate}
\item[(1)]
$$\frac{\partial (y_1,\cdots,y_n)}{\partial (t_1,\cdots,t_n)}
=\frac{ q y_1\cdots y_n}{q' t_1\cdots t_n}$$
where $q/q'\not\in {\mathbb Z}$.
\item[(2)]
$$\frac{\partial (y_1,\cdots y_n)}{\partial (t_1,\cdots, t_n)}
=q y_1\cdots y_n$$
where $q\in\Bbbk^{\times}$.
\item[(3)]
$$\frac{\partial (y_1,\cdots, y_n)}{\partial (t_1,\cdots, t_n)}
=\frac{ y_1^{1+p_1}\cdots y_n^{1+p_n}}{t_1^{1+q_1}\cdots t_n^{1+q_n}}$$
with $\gcd(p_1,\cdots,p_n)>\gcd(q_1,\cdots,q_n)$.
\end{enumerate}
\end{corollary}

\begin{proof} (1) Let $b=q t_1\cdots t_n$ and $a=q' t_1\cdots t_n$.
Then $\Bbbk^{(n,b)}\cong N_{q}$ and $\Bbbk^{(n,a)}=N_{q'}$. 
Since $q/q'\not\in {\mathbb Z}$, $N_q$ cannot be embedded to 
$N_{q'}$ by Theorem \ref{xxthm4.16}. The assertion follows from
Lemma \ref{xxlem7.2}. 

(2) Let $b=q t_1\cdots t_n$ and $a=1$.
Then $\Bbbk^{(n,b)}\cong N_{q}$ and $\Bbbk^{(n,a)}=N_{Weyl}$. 
By Corollary \ref{xxcor5.9}(2), there is no embedding from
$\Bbbk^{(n,b)}\to \Bbbk^{(n,a)}$. The assertion follows from
Lemma \ref{xxlem7.2}. 

(3) Let $b= t_1^{1+p_1}\cdots t_n^{1+p_n}$ and 
$a=t_1^{1+q_1}\cdots t_n^{1+q_n}$.
By Lemma \ref{xxlem4.7}(2), $\Bbbk^{(n,b)}\cong N(p)$
where $p=\gcd(p_s)$ and $\Bbbk^{(n,a)}=N(q)$ where $q=
\gcd(q_s)$. By Corollary \ref{xxcor4.9}, $N(p)$ cannot be embedded
to $N(q)$ when $p>q$. The assertion follows from
Lemma \ref{xxlem7.2}. 
\end{proof}

We will present two more examples for the rest of
this section. 

\begin{lemma}
\label{xxlem7.5}
Let $b=t_1\cdots t_n p(u)$ where $p(t)\in \Bbbk[t]$ is a 
polynomial of degree $h\geq 1$ and $u\in \Bbbk^{(n)}\setminus
\Bbbk$.
\begin{enumerate}
\item[(1)]
For $1\leq w\leq h-1$,
${^{w}\Gamma}_{0}(\Bbbk^{(n,b)})\supseteq \Bbbk[u]\neq \Bbbk$. 
\item[(2)]
If $u=t_1^{d_1}\cdots t_{n}^{d_n}$ with $\gcd(d_1,\cdots,d_n)=1$,
then ${^{h}\Gamma}_{0}(\Bbbk^{(n,b)})=\Bbbk$. 
\end{enumerate}
\end{lemma}

\begin{proof} Let $N$ be the Nambu-Poisson field $\Bbbk^{(n,b)}$.

(1) If $h=1$, there is nothing to prove. So we assume that 
$h\geq 2$. By Lemma \ref{xxlem2.15}(2), we may assume that 
$w=h-1$. Let $\nu\in {\mathcal V}_{h-1}(N)$ By definition,
$$\nu(p(u))=
\nu(\{t_1,\cdots,t_2\} t_1^{-1}\cdots t_n)\geq 0-(h-1).$$
We claim that $\nu(u)\geq 0$. Suppose to the contrary that
$\nu(u)<0$. Since $p(t)$ is a polynomial of degree $h\geq 2$,
by valuation axioms, 
$\nu(p(u))=\nu(u^h)=h\nu(u)\leq -h$. Then 
$-h\geq h\nu(u)\geq -(h-1)$, yielding a contradiction. Therefore, we proved the claim. Thus $u\in {^{(h-1)}\Gamma}_{0}(N)$ 
as desired.

(2) Up to an automorphism in $GL_{n}({\mathbb Z})$, we may 
assume that $u=t_1$. Define $\deg(t_1)=-1$ and $\deg(t_i)=0$
for all $i\geq 2$. Then, the induced filtration determined by
this degree assignment is a good filtration, and it is denoted
by $\nu$. We claim that $\nu\in {\mathcal V}_{h}(N)$. This 
follows from the following computation
$$\nu(\{t_1,\cdots,t_n\})=\nu(t_1\cdots t_n p(t_1))
=-1+h(-1)=\sum_{s=1}^{n} \nu(t_s)-h$$
and Lemma \ref{xxlem2.9}(3). Since $\nu(t_1)=-1$, $\nu(f(t_1))
<0$ for every polynomial $f(t)$ of positive degree. This 
means that $f(t_1)\not\in {^{h}\Gamma}_{0}(N)$.

By Theorem \ref{xxthm2.16}(2), ${^{h}\Gamma}_{0}(N)\subseteq
\Bbbk[t_1,\cdots,t_n]$. Similarly, one can check that
${^{h}\Gamma}_{0}(N)\subseteq \Bbbk[t_1,t_2^{e_2},\cdots,t_n^{e_n}]$
where $e_s$ is one or $-1$. Thus ${^{h}\Gamma}_{0}(N)\subseteq
\Bbbk[t_1]$. Combining with the result in the previous paragraph, 
we obtain that ${^{h}\Gamma}_{0}(N)\subseteq \Bbbk$. The assertion
follows from Lemma \ref{xxlem2.15}(2).
\end{proof}

\begin{theorem}
\label{xxthm7.6}
Suppose that $b=t_1\cdots t_n p(u)$ where $p(t)$ is a polynomial 
of degree $h\geq 2$ and $u\in \Bbbk^{(n)}\setminus \Bbbk$ and 
that $a=t_1\cdots t_n p_0(t_1)$ where $p_0(t)$ is a polynomial 
of degree $h_0$ strictly between $0$ and $h$. Then \eqref{E7.0.2} 
has no solution.
\end{theorem}

\begin{proof} By Lemma \ref{xxlem7.2} it suffices to show that 
there is no embedding from $\Bbbk^{(n,b)}$ to $\Bbbk^{(n,a)}$.

Suppose to the contrary that there is an embedding 
$\phi: \Bbbk^{(n,b)}\to \Bbbk^{(n,a)}$. Now we compute 
some $\gamma$-type invariants. By Lemma \ref{xxlem7.5}(2), 
${^{h_0}\Gamma}_{0}(\Bbbk^{(n,a)})=\Bbbk$. By Lemma 
\ref{xxlem7.5}(1), ${^{h_0}\Gamma}_{0}(\Bbbk^{(n,b)})\supseteq \Bbbk[u]$.
By Lemma \ref{xxlem2.15}(1), there is an injective 
map $\phi: {^{h_0}\Gamma}_{0}(\Bbbk^{(n,b)})\to
{^{h_0}\Gamma}_{0}(\Bbbk^{(n,a)})$, which yields a contradiction.
Therefore there is no embedding from $\Bbbk^{(n,b)}$ to 
$\Bbbk^{(n,a)}$ as required.
\end{proof}

Theorem \ref{xxthm7.6} provides a lot of examples of 
pairs $(a,b)$ such that \eqref{E7.0.2} has no solution.
The next lemma deals with a special case when $b=qt_1\cdots t_n$. 

\begin{lemma}
\label{xxlem7.7}
Let $a\in \Bbbk^{(n)}$ and $b=qt_1\cdots t_n$ where 
$q\in \Bbbk^{\times}$.
\begin{enumerate}
\item[(1)]
Suppose $\Phi^{-1}\not\in {^{0}\Gamma}_{0}(\Bbbk^{(n,a)})$. 
Then the PDE
\begin{equation}
\label{E7.7.1}\tag{E7.7.1}
\frac{\partial (y_1,\cdots,y_n)}{\partial (t_1,\cdots,t_n)}
=\frac{q y_1\cdots y_2}{a(t_1,\cdots,t_n)\Phi(t_1,\cdots,t_n)}
\end{equation}
has no solution.
\item[(2)]
There is no embedding from $N_{q}$ to $\Bbbk^{(n,a\Phi)}$.
\end{enumerate}
\end{lemma}

\begin{proof}
(1) Suppose to the contrary that \eqref{E7.7.1} has a solution
$(y_1,\cdots,y_n)$. In $\Bbbk^{(n,a)}$, we have
\begin{align}
\notag \{y_1,\cdots,y_n\}_{a}
&=\frac{\partial (y_1,\cdots,y_n)}{\partial (t_1,\cdots,t_n)}a\\
\label{E7.7.2}\tag{E7.7.2}
&= a\frac{q y_1\cdots y_2}{a(t_1,\cdots,t_n)\Phi(t_1,\cdots,t_n)}\\
\notag &=\frac{q y_1\cdots y_2}{\Phi(t_1,\cdots,t_n)}.
\end{align}

By definition, there is a $\nu\in {\mathcal V}_0(\Bbbk^{(n,a)})$ 
such that $\nu(\Phi^{-1})<0$. While, by \eqref{E7.7.2}, we have
$$\nu(\Phi^{-1})=\nu
(\{y_1,\cdots,y_n\}_{a} (qy_1\cdots y_n)^{-1})
\geq 0,$$
yielding a contradiction. Therefore, \eqref{E7.7.1} has no solution.

(2) This follows from part (1) and Lemma \ref{xxlem7.2}.
\end{proof}

\begin{corollary}
\label{xxcor7.8}
Let $b=q t_1\cdots t_n$ for some $q\in \Bbbk^{\times}$ and 
let $a=t_1\cdots t_n$ and $a'=a \Phi^{e}$ where $\Phi$ is in 
$\Bbbk^{[n]}$ with zero constant term and $e$ is $1$ or $-1$. 
\begin{enumerate}
\item[(1)]
The PDE
\begin{equation}
\label{E7.8.1}\tag{E7.8.1}
\frac{\partial (y_1,\cdots,y_n)}{\partial (t_1,\cdots,t_n)}
=\frac{q y_1\cdots y_2}{t_1\cdots t_n(\Phi(t_1,\cdots,t_n))^{e}}\, ,
\end{equation}
where $e=1$ or $-1$, has no solution.
\item[(2)]
There is no embedding from $N_q$ to $\Bbbk^{(n,a')}$.
\end{enumerate}
\end{corollary}

\begin{proof}
(1) Let $f=\Phi^{e}$ where $e$ is $1$ or $-1$. 
By Lemma \ref{xxlem7.7}(1), it suffices to show that
$f^{-1}\not\in {^{0}\Gamma}_{0}(\Bbbk^{(n,a)})
={^{0}\Gamma}_{0}(N_1)$. 

Case 1: $f=\Phi$, which is a polynomial in $t_s$ with 
zero constant term. It remains to show that there is a 
valuation $\nu\in {\mathcal V}_0(N_1)$ such that 
$\nu(f^{-1})<0$, or equivalently $\nu(\Phi)>0$. Since
$f$ has zero constant term, by choosing $\nu(t_i)=1$
for all $i$, we obtain a $0$-valuation $\nu$ on $N_1$ 
[Theorem \ref{xxthm4.3}] such that $\nu(f)>0$ as required. 

Case 2: $f=\Phi^{-1}$ where $\Phi$ is a polynomial in $t_s$ 
with zero constant term. It remains to show that there is a 
valuation $\nu\in {\mathcal V}_0(N_1)$ such that 
$\nu(f^{-1})<0$, or equivalently $\nu(\Phi)<0$. By 
choosing $\nu(t_i)=-1$ for all $i$, we obtain a $0$-valuation 
$\nu$ on $N_1$ [Theorem \ref{xxthm4.3}] such that $\nu(\Phi)>0$ 
as required.

(2) This follows from part (1) and Lemma \ref{xxlem7.2}.
\end{proof}

\subsection*{Acknowledgments.} 
The authors thank Ken Goodearl and Milen Yakimov for many valuable conversations and correspondences on the subject. Wang was partially supported by Simons 
collaboration grant \#688403 and Air Force Office of 
Scientific Research grant FA9550-22-1-0272. Zhang was 
partially supported by the US National Science Foundation 
(No. DMS-2001015 and DMS-2302087). Part of this research 
work was done during the first, second, and third authors' 
visit to the Department of Mathematics at the University of 
Washington in June 2022 and January 2023. They are 
grateful for the fourth author's invitation and wish to 
thank the University of Washington for its hospitality.

\providecommand{\bysame}{\leavevmode\hbox to3em{\hrulefill}\thinspace}
\providecommand{\MR}{\relax\ifhmode\unskip\space\fi MR }
\providecommand{\MRhref}[2]{%

\href{http://www.ams.org/mathscinet-getitem?mr=#1}{#2} }

\end{document}